\documentclass[twocolumn]{autart}    

\usepackage{tikz}
\usepackage{amsmath}
\usepackage{latexsym, amssymb, amsmath, amsbsy,amsopn, amstext}
\usepackage{graphicx}
\usepackage{hyperref}
\usepackage{cancel, color}
\usepackage{epstopdf}
\graphicspath{{images/}}
\DeclareMathOperator{\rank}{rank}
\DeclareMathOperator{\ad}{ad}

\DeclareMathOperator{\Span}{Span}
\DeclareMathOperator{\Col}{Col}
\DeclareMathOperator{\Row}{Row}
\DeclareMathOperator{\diag}{diag}
\DeclareMathOperator{\tr}{trace}
\DeclareMathOperator{\Tr}{Tr}
\DeclareMathOperator{\Dt}{Dt}

\DeclareMathOperator{\I}{{\bf 1}}
\DeclareMathOperator{\lcm}{lcm}
\DeclareMathOperator{\argmin}{argmin}
\def\lplus{\rotatebox[]{-90}{$\pm$}}
\def\rplus{\rotatebox[]{90}{$\pm$}}
\def\lminus{\vdash}
\def\rminus{\dashv}
\def\lvtimes{\vec{\ltimes}}
\def\rvtimes{\vec{\rtimes}}
\def\lvminus{\vec{\vdash}}
\def\rvminus{\vec{\dashv}}
\def\lvplus{\vec{\rotatebox[]{-90}{$\pm$}}}
\def\rvplus{\vec{\rotatebox[]{90}{$\pm$}}}
\def\bd{bd}
\def\pr{pr}
\def\gl{gl}
\def\GL{GL}
\def\PR{Pr}

\def\cal{\mathcal}

\def\pa{\partial}

\def\ra{\rightarrow}

\def\lra{\leftrightarrow}

\def\a{\alpha}
\def\b{\beta}
\def\d{\delta}

\def\D{\Delta}

\def\L{\Lambda}

\def\0{{\bf 0}}
\def\A{\left<A\right>}
\def\B{\left<B\right>}
\def\E{\left<E\right>}

\def\X{\lceil X \rceil}
\def\Y{\lceil Y \rceil}

\newcommand{\R}{{\mathbb R}}
\newcommand{\C}{{\mathbb C}}
\newcommand{\N}{{\mathbb N}}
\newcommand{\Q}{{\mathbb Q}}
\newcommand{\F}{{\mathbb F}}
\newcommand{\Z}{{\mathbb Z}}

\def\dsum{\mathop{\sum}\limits}
\newtheorem{dfn}[thm]{Definition}
\newtheorem{prp}[thm]{Proposition}
\newtheorem{exa}[thm]{Example}

\begin{document}

\begin{frontmatter}

\title{On Equivalence of Matrices\thanksref{footnoteinfo}}

\thanks[footnoteinfo]{This work is supported partly by National Natural Science Foundation of China under Grants 61773371 and 61733018. Corresponding author: Daizhan Cheng. Tel.: +86 10 8254 1232.}

\author{Daizhan Cheng}

\address{Key Laboratory of Systems and Control, AMSS,
Chinese Academy of Sciences, Beijing 100190, P.R.China}


\begin{abstract}  A new matrix product, called the semi-tensor product (STP), is briefly reviewed. The STP extends the classical matrix product to two arbitrary matrices. Under STP the set of matrices becomes a monoid (semi-group with identity). Some related structures and properties are investigated. Then the generalized matrix addition is also introduced, which extends the classical matrix addition to a class of two matrices with different dimensions.
%

Motivated by STP of matrices, two kinds of equivalences of matrices (including vectors) are introduced, which are called matrix equivalence (M-equivalence) and vector equivalence (V-equivalence) respectively. The lattice structure has been established for each equivalence. Under each equivalence, the corresponding quotient space becomes a vector space. Under M-equivalence, many algebraic, geometric, and analytic structures have been posed to the quotient space, which include (i) lattice structure; (ii) inner product and norm (distance); (iii) topology; (iv) a fiber bundle structure, called the discrete bundle; (v) bundled differential manifold; (vi) bundled Lie group and Lie algebra. Under V-equivalence, vectors of different dimensions form a vector space ${\cal V}$, and a matrix $A$ of arbitrary dimension is considered as an operator (linear mapping) on ${\cal V}$. When $A$ is a bounded operator (not necessarily square but includes square matrices as a special case), the generalized characteristic function, eigenvalue and eigenvector etc. are defined.

In one word, this new matrix theory overcomes the dimensional barrier in certain sense. It  provides much more freedom for using matrix approach to practical problems.

\end{abstract}

\end{frontmatter}

\section{Preliminaries}

\subsection{Contents}

For convenience in read, a list of contents is given as follows.

\vskip 2mm

I. Preliminaries

\begin{enumerate}
\item[1.1] Contents
\item[1.2] Introduction
\item[1.3] Symbols
\end{enumerate}

II. M-equivalence and Lattice Structure

\begin{enumerate}
\item[2.1] STP of Matrices
\item[2.2] M-Equivalence of Matrices
\item[2.3] Lattice Structure on ${\cal M}_{\mu}$
\item[2.4] Monoid and Quotient Monoid
\item[2.5] Group Structure of ${\cal M}^{\mu}$
\item[2.6] Semi-tensor Addition and Vector Space Structure of $\Sigma_{\mu}$
\end{enumerate}

III. Topology on M-equivalence Space

\begin{enumerate}
\item[3.1] Topology via Sub-basis
\item[3.2] Fiber Bundle Structure on ${\cal M}_{\mu}$
\item[3.3] Coordinate Frame on ${\cal M}_{\mu}$
\item[3.4] Inner Product on ${\cal M}_{\mu}$
\item[3.5] $\Sigma_{\mu}$ as a Matric Space
\item[3.6] Sub-space of $\Sigma_{\mu}$
\end{enumerate}

IV. Differential Structure on M-equivalence Space

\begin{enumerate}
\item[4.1] Bundled Manifold
\item[4.2] $C^r$ Functions on $\Sigma_{\mu}$
\item[4.3] Generalized Inner Product
\item[4.4] Vector Fields
\item[4.5] Integral Curves
\item[4.6] Forms
\item[4.7] Tensor Fields
\end{enumerate}

V. Lie Algebra Structure on Square M-equivalence Space

\begin{enumerate}
\item[5.1] Lie Algebra on Square M-equivalence Space
\item[5.2] Bundled Lie Algebra
\item[5.3] Bundled Lie Sub-algebra
\item[5.4] Further Properties of $\gl(\F)$
\end{enumerate}

VI. Lie Group on Nonsingular M-equivalence Space

\begin{enumerate}
\item[6.1] Bundled Lie group
\item[6.2] Relationship with $\gl(\F)$
\item[6.3] Lie Subgroup of $\GL(\F)$
\item[6.4] Symmetric Group
\end{enumerate}

VII. V-equivalence

\begin{enumerate}
\item[7.1] Equivalence of Vectors of Different Dimensions
\item[7.2] Vector Space Structure on V-equivalence Space
\item[7.3] Inner Product and Linear Mappings
\item[7.4] Type-1 Invariant Subspace
\item[7.5] Type-2 Invariant Subspace
\item[7.6] Higher Order Linear Mapping
\item[7.7] Invariant Subspace on V-equivalence Space
\item[7.8] Generalized Linear System
\end{enumerate}

VIII.  Concluding Remarks

\subsection{Introduction}

Matrix theory and calculus are two classical and fundamental mathematical tools in modern science and technology. There are two mostly used operators on the set of matrices: conventional matrix product and matrix addition. Roughly speaking, the object of matrix theory is $({\cal M},+,\times)$, where ${\cal M}$ is the set of all matrices. Unfortunately, unlike the arithmetic system $(\R,+,\times)$, in matrix theory both ``$+$" and ``$\times$" are restricted by the matrix dimensions. Precisely speaking: consider two matrices $A\in {\cal M}_{m\times n}$ and $B\in {\cal M}_{p\times q}$. Then the ``product", $A\times B$, is well posed, if and only if, $n=p$; the ``addition" $A+B$, is defined, if and only if, $m=p$ and $n=q$.  Though there are some other matrix products such as Kronecker product, Hadamard product etc., but they are of different stories \cite{hor85}.

The main purpose of this paper is to develop a new matrix theory, which intends to overcome the dimension barrier by extending the matrix product and matrix addition to two matrices which do not meet the classical dimension requirement. As generalizations of the classical ones, they should be consistent with the classical definitions. That is, when the dimension requirements of two argument matrices in classical matrix theory are satisfied, the newly defined operators should coincide with the original one.

Because of the extension of the two fundamental operators, many related concepts can be extended. For instance, the characteristic functions, the eigenvalues and eigenvectors of a square matrix can be extended to certain non-square matrices; Lie algebraic structure can also be extended to dimension-varying square matrices. All these extensions should be consistent with the classical ones. In one word, we are developing the classical matrix theory but not violating any of the original matrix theory.

When we were working on generalizing the fundamental matrix operators we meet a serious problem: Though the extended operators are applicable to certain sets of matrices with different dimensions, they fail to be vector space anymore. This drawback is not acceptable, because it blocked the way to extend many nice algebraic or geometric structures in matrix theory, such as Lie algebraic structure, manifold structure etc.,  to the enlarged set, which includes matrices of different dimensions. To overcome this obstacle, we eventually introduced certain equivalence relations. Then the quotient spaces, called the equivalence spaces, become vector spaces. Two equivalence relations have been proposed. They are \textit{matrix equivalence} (M-equivalence) and \textit{vector equivalence} (V-equivalence).

Then many nice algebraic, analytic, and geometric structures have been developed on the M-equivalence spaces. They are briefly introduced as follows:

\begin{itemize}
\item Lattice structure: The elements in each equivalent class form a lattice. The class of spaces with  different dimensions also form a lattice. The former and the latter are homomorphic. The lattices obtained for M-equivalence and V-equivalence are isomorphic.
\item Topological structure: A topological structure is proposed to the M-equivalence space, which is built by using topological sub-base. It is then proved that under this topology the equivalence space is Hausdorff ($T_2$) space and is second countable.
\item Inner product structure: An inner product is proposed on the M-equivalence space. The norm (distance) is also obtained. It is proved that the topology induced by this norm is the same as the topology produced by using the topological sub-base.
\item Fiber bundle structure. A fiber bundle structure is proposed for the set of matrices (as total space) and the equivalent classes (as base space). The bundle structure is named the discrete bundle, because each fiber has discrete topology.
\item Bundled manifold structure: A (dimension-varying) manifold structure is proposed for the M-equivalence space. Its coordinate charts are constructed via the discrete bundle. Hence it is called a bundled manifold.
\item Bundled Lie algebraic structure: A Lie algebra structure is proposed for the M-equivalence space. The Lie algebra is of infinite dimensional, but almost all the properties of finite dimensional Lie algebras  remain available.
\item Bundled Lie group structure: For the M-equivalence classes of square nonsingular matrices a group structure is proposed. It has also the dimension-varying manifold structure. Both the algebraic and the geometric structures are consistent and hence it becomes a Lie group. The relationship of this Lie group with the bundled Lie algebra is also investigated.
\end{itemize}

Under V-equivalence, all the vectors of varying dimensions form a vector space ${\cal V}$, and any matrix $A$ can be considered as a linear operator on ${\cal V}$. A very important class of $A$, called the bounded operator, is investigated in detail. For a bounded operator $A$, which could be non-square, its characteristic function is proposed. Its eigenvalues and the corresponding eigenvectors are obtained. A generalized $A$-invariant subspace has been discussed in detail.

This work is motivated by the semi-tensor product (STP). The STP of matrices was proposed firstly and formally in 2001 \cite{che01}. Then it has been used to some Newton differential dynamic systems and their control problems \cite{che04}, \cite{xue06}, \cite{mei10}. A basic summarization was given in \cite{che07}.

Since 2008, STP has been used to formulate and analyze Boolean networks as well as general logical dynamic systems, and to solve control design problems for those systems. It has then been developed rapidly. This is witnessed by hundreds of research papers.  A multitude of applications of STP include (i) logical dynamic systems \cite{che11}, \cite{for13}, \cite{las13}; (ii) systems biology \cite{zha13}, \cite{gao13}; (iii) graph theory and formation control \cite{wan12}, \cite{zha13d}; (iv) circuit design and failure detection \cite{li12}, \cite{li12h}, \cite{che13};
(v) game theory \cite{guo13}, \cite{che14}, \cite{che15};  (vi) finite
automata and symbolic dynamics \cite{xu13}, \cite{zhapr}, \cite{hoc13}; (vii) coding and cryptography
\cite{zho15}, \cite{zha14d}; (viii) fuzzy control \cite{che12d}, \cite{fen13}; (ix) some engineering applications \cite{wu15}, \cite{liu13c}; and many other topics \cite{che12}, \cite{lu16}, \cite{yan15}, \cite{zho16}, \cite{zou15}, \cite{liu16}; just to name a few.

As a generalization of conventional matrix product, STP is applicable to two matrices of arbitrary dimensions.
In addition, this generalization keeps all fundamental properties of conventional matrix product available.  Therefore, it becomes a very conventional tool for investigating many matrix expression related problems.

Recently, the journal \textit{IET Control Theory \& Applications} has published a special issue ``Recent Developments in Logical Networks and Its Applications". It provides many up-to-date results of STP and its applications. Particularly, we refer to a survey paper \cite{lu17} in this special issue for a general introduction to STP with applications.

Up to this time, the main effort has been put on its applications. Now when we start to explore the mathematical foundation of STP, we found that the most significant characteristic of STP is that it overcomes the dimension barrier. After serious thought, it can be seen that in fact STP is defined on equivalent classes. Following this thought of train, the matrix theory on equivalence space emerges. The outline of this new matrix theory is presented in this manuscript. The results in this paper are totally new except the concepts and basic properties of STP, which are presented in subsection 2.1.

\subsection{Symbols}

For statement ease, we first give some notations:

\begin{enumerate}
\item $\N$: Set of natural numbers (i.e., $\N=\{1,2,\cdots\}$);

\item $\Z$: set of integers;

\item $\Q$: Set of rational numbers, ($\Q_+$: Set of positive rational numbers);

\item $\R$: Field of real numbers;

\item $\C$: Field of complex numbers;

\item $\F$: certain field of characteristic number $0$ (Particularly, we can understand $\F=\R$ or $\F=\C$).

\item ${\cal M}^{\F}_{m\times n}$: set of $m\times n$ dimensional matrices over field $\F$. When the field $\F$ is obvious or does not affect the discussion, the superscript $\F$ can be omitted, and as a default: $\F=\R$ can be assumed.

\item $\Col(A)$ ($\Row(A)$): the set of columns (rows) of ~$A$; $\Col_i(A)$ ($\Row_i(A)$): the $i$-th column (row) of ~$A$.

\item ${\cal D}_k=\{1,2,\cdots,k\}$,  ${\cal D}:={\cal D}_2$;

\item $\d_n^i$: the $i$-th column of the identity matrix ~$I_n$;

\item $\D_n=\{\d_n^i\,|\,i=1,2,\cdots,n\}$;

\item $L\in {\cal M}_{m\times r}$ is a logical matrix, if ~$\Col(L)\subset \D_m$. The set of  $m\times r$ logical matrices is denoted as ~${\cal L}_{m\times r}$;

\item Assume $A\in {\cal L}_{m\times r}$. Then $L=\left[\d_m^{i_1},\d_m^{i_2},\cdots,\d_m^{i_r}\right]$. It is briefly denoted as
$$
L=\d_m[i_1,i_2,\cdots,i_r].
$$
\item Let $A=(a_{i,j})\in {\cal M}_{m\times n}$, and $a_{i,j}\in \{0,1\}$, $\forall i,j$. Then $A$ is
called a Boolean matrix. Denote the set of $m\times n$ dimensional Boolean matrices by ${\cal B}_{m\times n}$.

\item  Set of probabilistic vectors:
$$
\varUpsilon_k:=\left\{(r_1,r_2,\cdots,r_k)^T\;\left|\; r_i\geq 0,~\dsum_{i=1}^kr_i=1 \right.\right\}.
$$

\item  Set of probabilistic matrices:
$$
\varUpsilon_{m\times n}:=\left\{M\in {\cal M}_{m\times n}\;\big|\; \Col(M)\subset \varUpsilon_m\right\}.
$$

\item
$a|b$: $a$ divides $b$.

\item
$m\wedge n=gcd(m,n)$: The greatest common divisor of $m,n$.

\item
$m\vee n=lcm(m,n)$: The least common multiple of $m,n$.

\item $\ltimes$ ($\rtimes$): The left (right) semi-tensor product of matrices.

\item $\lvtimes$ ($\rvtimes$): The left (right) vector product of matrices.

\item $\sim$ ($\sim_{\ell}$, $\sim_r$): The M-equivalence ((left, right) matrix equivalence).

\item $\lra$ ($\lra_{\ell}$, $\lra_r$): The V-equivalence ((left,right) vector equivalence).


\item The set of all matrices:
$$
{\cal M}=\bigcup_{i=1}^{\infty}\bigcup_{j=1}^{\infty}{\cal M}_{i\times j}.
$$
\item
$$
{\cal M}_{\cdot\times q}:=\left\{A\in {\cal M}\;\big|\; \mbox{column number of}~ A~\mbox{is}~q\right\}.
$$
\item The set of matrices:
$$
{\cal M}_{\mu}:=\left\{A\in {M}_{m\times n}\;\big|\; m/n=\mu \right\}
$$
\item Lattice homomorphism: $\approx$
\item Lattice isomorphism: $\approxeq$
\item Vector order: $\leqslant$
\item Vector space order: $\sqsubseteq$
\item Matrix order: $\prec$
\item Matrix space order: $\sqsubset$
\item The overall matrix quotient space:
$$
\Sigma_{{\cal M}}={\cal M}/\sim.
$$
\item The $\mu$-matrix quotient space:
$$
\Sigma_{\mu}:={\cal M}_{\mu}/\sim.
$$
\item The $\mu=1$-matrix quotient space:
$$
\Sigma:={\cal M}_{1}/\sim.
$$
\item The set of all vectors:
$$
{\cal V}=\bigcup_{n=1}^{\infty} {\cal V}_n, \qquad (\text{Note that real~} {\cal V}_n\sim \R^n).
$$

\item The vector quotient space under V-equivalence:
$$
\Omega_{{\cal V}}:={\cal V}/\lra.
$$
\item The vector quotient subspace under V-equivalence:
$$
\Omega_{{\cal V}}^i:={\cal V}^{[i,\cdot]}/\lra.
$$
\item The matrix quotient space under V-equivalence:
$$
\Omega_{{\cal M}}:={\cal M}/\lra.
$$
\item The matrix quotient subspace under V-equivalence:
$$
\Omega_{{\cal M}}^n:={\cal M}_{\cdot\times n}/\lra.
$$
\item Given a ($C^{r}$) manifold $M$ ($r$ could be $\infty$ or $\omega$),
\begin{itemize}
\item its tangent space is
$T(M)$;
\item its cotangent space is $T^*(M)$;
\item the set of $C^{r}$ functions is
$C^{r}(M)$;
\item the set of vector fields is $V^{r}(M)$; the set of co-vector fields is
$V^{*r}(M)$.
\item the set of tensor fields on $M$ with covariant order $\a$ and contravariant order $\b$ is ${\bf T}^{\a}_{\b}(M)$;
when $\b=0$ it becomes ${\bf T}^{\a}(M)$.

\end{itemize}

\end{enumerate}

\section{M-equivalence and Lattice Structure}

\subsection{STP of Matrices}

This section gives a brief review on STP. We refer to \cite{che07}, \cite{che11}, \cite{che12} for details.

\begin{dfn}\label{d2.2.1} Let $A\in {\cal M}_{m\times n}$, $B\in {\cal M}_{p\times q}$, and $t=n\vee p$ be
the least common multiple of $n$ and $p$. Then
\begin{enumerate}
\item
the left STP of $A$ and $B$, denoted by $A\ltimes B$, is
defined as
\begin{align}\label{2.2.1}
A\ltimes B:=\left(A\otimes I_{t/n}\right)\left(B\otimes I_{t/p}\right),
\end{align}
where $\otimes$ is the Kronecker product \cite{hor85};
\item
the right STP of
$A$ and $B$ is defined as
\begin{align} \label{2.2.2} A\rtimes B:=\left(I_{t/n}\otimes A\right)\left( I_{t/p} \otimes B\right).
\end{align}
\end{enumerate}
\end{dfn}

In the following we mainly discuss the left STP, and briefly call the left STP as STP. Most of the properties of left STP have their corresponding ones for right STP. Please also refer to \cite{che07} or \cite{che11} for their major differences.

\begin{rem}\label{r2.2.2}
If ~$n=p$, both left and right STP defined in Definition \ref{d2.2.1} degenerate to the conventional matrix product. That is, STP is a generalization of the conventional matrix product. Hence, as a default, in most cases the symbol $\ltimes$ can be omitted (but not $\rtimes$). \textbf{That is, unless elsewhere stated throughout this paper}
\begin{align} \label{2.2.201}
AB:=A\ltimes B.
\end{align}

The following proposition shows that this generalization not only keeps the main properties of conventional matrix product available, but also adds some new properties such as certain commutativity.
\end{rem}

Associativity and distribution are two fundamental properties of conventional matrix product. When the product is generalized to STP, these two properties remain available.

\begin{prp}\label{p2.2.3}
\begin{enumerate}
\item[{\rm 1.}] (Distributive Law)
\begin{align} \label{2.3}
\begin{cases} F\ltimes (aG\pm bH)=aF\ltimes G\pm bF\ltimes H,
\\ (aF\pm bG)\ltimes H=a F\ltimes H \pm bG\ltimes H, \quad a, b\in
\F. \end{cases} \end{align}

\item[{\rm 2.}] (Associative Law)
\begin{align} \label{2.2.4} (F\ltimes G)\ltimes H = F\ltimes (G\ltimes H).
\end{align}
\end{enumerate}
\end{prp}

The following proposition is inherited from the conventional matrix product.

\begin{prp}\label{p2.2.4}
\begin{enumerate}
\item[{\rm 1.}]
\begin{align}
\label{2.2.5} (A\ltimes B)^T=B^T\ltimes A^T.
\end{align}

\item[{\rm 2.}] Assume ~$A$ and $B$ are invertible, then
\begin{align}
\label{2.2.6} (A\ltimes B)^{-1}=B^{-1}\ltimes A^{-1}.
\end{align}
\end{enumerate}
\end{prp}
The following proposition shows that the STP has certain commutative property.
\begin{prp}\label{p2.2.5} Given ~$A\in {\cal M}_{m\times n}$.

\begin{enumerate}

\item[{\rm 1.}] Let ~$Z\in \R^t$ be a column vector. Then
\begin{align}
\label{2.2.7} ZA=( I_t\otimes A)Z.
\end{align}

\item[{\rm 2.}] Let ~$Z\in \R^t$ be a row vector. Then
\begin{align}
\label{2.2.8} AZ=Z(I_t\otimes A).
\end{align}
\end{enumerate}
\end{prp}

To explore further commutating properties, we introduce a swap matrix.

\begin{dfn}[\cite{hor85}]\label{d2.2.6} A swap matrix ~$W_{[m,n]}\in {\cal M}_{mn\times
mn}$ is defined as follows:
\begin{align}
\label{2.2.9}
\begin{array}{ccr}
W_{[m,n]}&=&\d_{mn}[1,m+1,\cdots,(n-1)m+1;\\
~&~&2,m+2,\cdots,(n-1)m+2;\\
~&~&\cdots~;~m,2m,\cdots,nm].
\end{array}
\end{align}
\end{dfn}
The following proposition shows that the swap matrix is orthogonal.
\begin{prp}\label{p2.2.7}
\begin{align} \label{2.2.10} W_{[m, n]}^T=W_{[m, n]}^{-1}=W_{[n,m]}.
\end{align}
\end{prp}

The fundamental function of the swap matrix is to ``swap" two factors.

\begin{prp}\label{p2.2.8}
\begin{enumerate}
\item[{\rm 1.}] Let ~$X\in \R^m$, $Y\in \R^n$ be two column vectors. Then
\begin{align} \label{2.2.11} W_{[m, n]}\ltimes X\ltimes Y=Y\ltimes X.
\end{align}
\item[{\rm 2.}] Let ~$X\in \R^m$, $Y\in \R^n$ be two row vectors. Then
\begin{align} \label{2.2.12} X\ltimes Y\ltimes W_{[m, n]}=Y\ltimes X.
\end{align}
\end{enumerate}
\end{prp}

The swap matrix can also swap two factor matrices of the Kronecker product \cite{che07}, \cite{wanpr}.
\begin{prp}\label{p2.2.9} Let $A\in {\cal M}_{m\times n}$ and $B\in {\cal M}_{p\times q}$. Then
\begin{align}
\label{2.2.13}
W_{[m,p]}(A\otimes B)W_{[q,n]}=B\otimes A.
\end{align}
\end{prp}

\begin{rem}\label{r2.2.10} Assume $A\in {\cal M}_{n\times n}$ and $B\in {\cal M}_{p\times p}$ are square matrices. Then (\ref{2.2.13}) becomes
\begin{align}\label{2.2.14}
W_{[n,p]}(A\otimes B)W_{[p,n]}=B\otimes A.
\end{align}
As a consequence, $A\otimes B$ and $B\otimes A$ are similar.
\end{rem}

The following example is an application of Proposition \ref{p2.2.9}.

\begin{exa}\label{e2.2.11} Prove
\begin{align}\label{2.2.15}
e^{A\otimes I_k}=e^A\otimes I_k.
\end{align}

Assume $A\in {\cal M}_{n\times n}$ is a square matrix and $B=A\otimes I_k$. Note that
$$
W(A\otimes I_k)W^{-1}=I_k\otimes A=\diag\underbrace{\left(A,A,\cdots,A\right)}_k,
$$
where $W=W_{[n,k]}$.
Then
$$
\begin{array}{ccl}
e^B&=&e^{W^{-1}(I_k\otimes A)W}\\
~&=&W^{-1}e^{\diag(A,A,\cdots,A)}W\\
~&=&W^{-1}\diag(e^A,e^A,\cdots,e^A)W\\
~&=&W^{-1}\left[I_k\otimes e^A\right]W=e^A\otimes I_k.
\end{array}
$$
\end{exa}

\begin{rem}\label{r2.2.12} Comparing the product of numbers with the product of matrices, two major differences are (i) matrix product has dimension restriction while the scalar product has no restriction; (ii) the matrix product is not commutative while the scalar product is. When the conventional matrix product is extended to STP, these two weaknesses have been eliminated in certain degree. First, the dimension restriction has been removed.  Second, in addition to Proposition ~\ref{p2.2.5}, which shows certain commutativity, the use of swap matrix also provides certain commutativity property. All these improvements make the STP more convenient in use than the conventional matrix product.
\end{rem}

\subsection{M-equivalence of Matrices}

The set of all matrices (over certain field $\F$) is denoted by ${\cal M}$, that is
$$
{\cal M}=\bigcup_{m=1}^{\infty}\bigcup_{n=1}^{\infty}{\cal M}_{m\times n}.
$$

It is obvious that STP is an operator defined as $\ltimes (~\mbox{or}~\rtimes): {\cal M}\times{\cal M}\ra {\cal M}$. Observing the definition of STP carefully, it is not difficult to find that when we use STP to multiply $A$ with $B$, instead of $A$ itself, we modify $A$ by Kronecker multiplying different sizes
of identity to multiply different $B$'s. In fact, STP multiplies an equivalent class of $A$, precisely,
$\A=\{A, A\otimes I_2, A\otimes I_3,\cdots\}$, with an equivalent class of $B$, that is, $\B=\{B,~B\otimes I_2, ~B\otimes I_3,\cdots\}$.

Motivated by this fact, we first propose an equivalence over set of matrices, called the matrix equivalence ($\sim_{\ell}$ or $\sim_r$).  Then STP can be considered as an operator over the equivalent classes. We give a rigorous definition for the equivalence.

\begin{dfn}\label{d2.3.1} Let $A,B\in {\cal M}$ be two matrices.
\begin{enumerate}
\item $A$ and $B$ are said to be left matrix equivalent (LME), denoted by $A\sim_{\ell} B$, if there exist two identity matrices $I_s$, $I_t$, $s,t\in \N$, such that
$$
A\otimes I_s=B\otimes I_t.
$$
\item  $A$ and $B$ are said to be right matrix equivalent (RME), denoted by $A\sim_r B$, if there exist two identity matrices $I_s$, $I_t$, $s,t\in \N$, such that
$$
I_s\otimes A=I_t\otimes B.
$$
\end{enumerate}
\end{dfn}

\begin{rem}\label{r2.3.2} It is easy to verify that the LME $\sim_{\ell}$ (similarly, RME $\sim_r$) is an
 equivalence relation. That is, it is (i) self-reflexive ($A\sim_{\ell} A$); (ii) symmetric (if $A\sim_{\ell} B$,
 then $B\sim_{\ell} A$); and (iii) transitive (if $A\sim_{\ell} B$, and  $B \sim_{\ell} C$, then  $A \sim_{\ell} C$).
\end{rem}

\begin{dfn}\label{d2.3.3} Given $A\in {\cal M}$.
\begin{enumerate}
\item The left equivalent class of $A$ is denoted by
$$
\A_{\ell}:=\left\{B\;|\;B\sim_{\ell} A\right\};
$$
\item The right equivalent class of $A$ is denoted by
$$
\A_{r}:=\left\{B\;|\;B\sim_{r} A\right\}.
$$
\item $A$ is left (right) reducible, if there is an $I_s$, $s\geq 2$, and a matrix $B$, such that
$A=B\otimes I_s$ (correspondingly, $A=I_s\otimes B$). Otherwise, $A$ is left (right) irreducible.
\end{enumerate}
\end{dfn}

\begin{lem}\label{l2.3.4} Assume $A\in {\cal M}_{\b\times \b}$ and $B\in {\cal M}_{\a\times \a}$, where $\a,\b\in \N$, $\a$ and $\b$ are co-prime, and
\begin{align}\label{2.3.1}
A\otimes I_{\a}=B\otimes I_{\b}.
\end{align}
Then there exists a $\lambda\in \F$ such that
\begin{align}\label{2.3.2}
A=\lambda I_{\b},\quad B=\lambda I_{\a}.
\end{align}
\end{lem}

\noindent\textit{Proof}. Split $A\otimes I_{\a}$ into equal size blocks as
$$
A\otimes I_{\a}=\begin{bmatrix}
A_{11}&\cdots&A_{1 \a}\\
\vdots&~&~\\
A_{\a 1}&\cdots&A_{\a\a}
\end{bmatrix}
$$
where $A_{ij}\in {\cal M}_{\b\times \b}$, $i,j=1,\cdots,\a$. Then we have
\begin{align}\label{2.3.3}
A_{i,j}=b_{i,j}I_{\b}.
\end{align}
Note that $\a$ and $\b$ are co-prime. Comparing the entries of both sides of (\ref{2.3.3}), it is clear that (i) the diagonal elements of all $A_{ii}$ are the same; (ii) all other elements ($A_{ij}$, $j\neq i$) are zero. Hence $A=b_{11}I_{\b}$. Similarly, we have $B=a_{11}I_{\a}$. But (\ref{2.3.1}) requires $a_{11}=b_{11}$, which is the required $\lambda$. The conclusion follows.
\hfill $\Box$

\begin{thm}\label{t2.3.5}
\begin{enumerate}
\item If $A\sim_{\ell} B$, then there exists a matrix $\Lambda$ such that
\begin{align}\label{2.3.4}
A=\Lambda\otimes I_{\b},\quad B=\Lambda\otimes I_{\a}.
\end{align}
\item In each class $\A_{\ell}$ there exists a unique $A_1\in \A_{\ell}$, such that $A_1$ is left irreducible. \end{enumerate}
\end{thm}

\noindent\textit{Proof}.
\begin{enumerate}
\item Assume $A\sim_{\ell} B$, that is, there exist $I_{\a}$ and $I_{\b}$ such that
\begin{align}\label{2.3.5}
A\otimes I_{\a}=B\otimes I_{\b}.
\end{align}
Without loss of generality, we assume $\a$ and $\b$ are co-prime. Otherwise, assume their greatest common divisor is $r=\a\wedge \b$, the $\a$ and $\b$ in (\ref{2.3.5}) can be replaced by $\a/r$ and $\b/r$ respectively.

Assume $A\in {\cal M}_{m\times n}$ and $B\in {\cal M}_{p\times q}$. Then
$$
m\a=p\b,\quad n\a=q\b.
$$
Since $\a$ and $\b$ are co-prime, we have
$$
m=s\b,\quad n=t\b,\quad p=s\a,\quad q=t\a.
$$
Split $A$ and $B$ into block forms as
$$
A=\begin{bmatrix}
A_{11}&\cdots&A_{1t}\\
\vdots&~&~\\
A_{s1}&\cdots&A_{st}\\
\end{bmatrix},\quad
B=\begin{bmatrix}
B_{11}&\cdots&B_{1t}\\
\vdots&~&~\\
B_{s1}&\cdots&B_{st}\\
\end{bmatrix},
$$
where $A_{i,j}\in {\cal M}_{\b\times \b}$ and $B_{i,j}\in {\cal M}_{\a\times \a}$, $i=1,\cdots,s$, $j=1,\cdots,t$.
Now (\ref{2.3.5}) is equivalent to
\begin{align}\label{2.3.6}
A_{ij}\otimes I_{\a}=B_{ij}\otimes I_{\b},\quad \forall i,j.
\end{align}
According to Lemma \ref{l2.3.4}, we have $A_{ij}=\lambda_{ij}I_{\b}$ and $B_{ij}=\lambda_{ij}I_{\a}$.
Define
$$
\Lambda:=
\begin{bmatrix}
\lambda_{11}&\cdots&\lambda_{1t}\\
\vdots&~&~\\
\lambda_{s1}&\cdots&\lambda_{st}\\
\end{bmatrix},
$$
equation (\ref{2.3.4}) follows.
\item
For each $A\in \A_{\ell}$ we can find $A_0$ irreducible such that $A=A_0\otimes I_s$. To prove it is unique, let $B\in \A_{\ell}$ and $B_0$ is irreducible and $B=B_0\otimes I_t$. We claim that $A_0=B_0$. Since $A_0\sim_{\ell}B_0$, there exists $\Gamma$ such that
$$
A_0=\Gamma\otimes I_p,\quad B_0=\Gamma\otimes I_q.
$$
Since both $A_0$ and $B_0$ are irreducible, we have $p=q=1$, which proves the claim.
\end{enumerate}
\hfill $\Box$

\begin{figure}[!htbp]
\centering
\includegraphics[width=3cm]{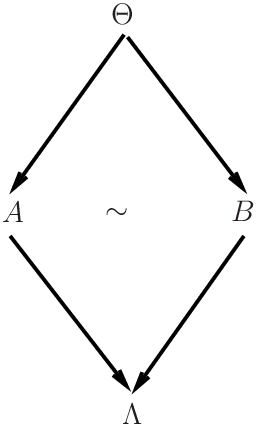}
\caption{$\Theta=lcm(A,B)$ and $\Lambda=gcd(A,B)$}\label{Fig.1}
\end{figure}

\begin{rem}\label{r2.3.6} Theorem \ref{t2.3.5} is also true for $\sim_r$ with obvious modification.
\end{rem}
\begin{rem}\label{r2.3.7}For statement ease, we propose the following terminologies:
\begin{enumerate}
\item  If $A=B\otimes I_s$, then $B$ is called a divisor of $A$ and $A$ is called a multiple of $B$.
\item  If (\ref{2.3.5}) holds and $\a,\b$ are co-prime, then the $\Lambda$ satisfying (\ref{2.3.4}) is called the greatest common divisor of $A$ and $B$. Moreover, $\Lambda=gcd(A,B)$ is unique.
\item  If (\ref{2.3.5}) holds and $\a,\b$ are co-prime, then
\begin{align}\label{2.3.601}
\Theta:=A\otimes I_{\a}=B\otimes I_{\b}
\end{align}
is called the least common multiple  of $A$ and $B$. Moreover, $\Theta=\lcm(A,B)$ is unique. (Refer to Fig.~\ref{Fig.1}.)

\item Consider an equivalent class $\A$, denote the unique irreducible element by $A_1$, which is called the root element. All the elements in $\A$ can be expressed as
 \begin{align}\label{2.3.7}
 A_i=A_1\otimes I_i,\quad i=1,2,\cdots.
 \end{align}
$A_i$ is called the $i$-th element of $\A$. Hence, an equivalent class $\A$ is a well ordered sequence as:
$$
\A=\left\{A_1,~A_2,~A_3,\cdots\right\}.
$$
\end{enumerate}
\end{rem}

Next, we modify some classical matrix functions to make them available for the equivalence class.

\begin{dfn}\label{d4.2.7}
\begin{enumerate}
\item Let $A\in {\cal M}_{n\times n}$. Then a modified determinant is defined as
\begin{align}\label{4.2.801}
\Dt(A)=[|\det(A)|]^{1/n}.
\end{align}
\item
Consider an equivalence of square matrices $\A$, the ``determinant" of $\A$ is defined as
\begin{align}\label{4.2.9}
\Dt(\A)=\Dt(A),\quad A\in \A.
\end{align}
\end{enumerate}
\end{dfn}

\begin{prp}\label{p4.2.8} (\ref{4.2.9}) is well defined, i.e., it is independent of the choice of the representative $A$.
\end{prp}

\noindent\textit{Proof}. To see (\ref{4.2.9}) is well defined, we need to check that $A\sim B$ implies $\Dt(A)=\Dt(B)$. Now assume $A=\Lambda\otimes I_{\b}$, $B=\Lambda\otimes I_{\a}$ and $\Lambda\in {\cal M}_{k\times k}$, then
$$
\Dt(A)=\left[|\det\left(\Lambda\otimes I_{\b}\right)|\right]^{1/k\b}=\left[|\det(\Lambda)\right|]^{1/k},
$$
$$
\Dt(B)=\left[|\det\left(\Lambda\otimes I_{\a}\right)|\right]^{1/k\a}=\left[|\det(\Lambda)|\right]^{1/k}.
$$
It follows that (\ref{4.2.9}) is well defined.
\hfill $\Box$

\begin{rem}\label{r4.2.801}
\begin{enumerate}
\item Intuitively, $\Dt(\A)$ defines only the ``absolute value" of $\A$. Because if there exists an $A\in \A$ such  that $\det(A)<0$, then $\det(A\otimes I_2)>0$. it is not able to define $\det(\A)$ uniquely over the class.
\item When $\det(A)=s$, $\forall A\in \A$, we also use $\det(\A)=s$. But when $\F=\R$, only $\det(\A)=1$ makes sense.
\end{enumerate}
\end{rem}

\begin{dfn}\label{d4.2.9}
\begin{enumerate}
\item  Let $A\in {\cal M}_{n\times n}$. Then a modified trace is defined as
\begin{align}\label{4.2.901}
\Tr(A)=\frac{1}{n}\tr(A).
\end{align}
\item
Consider an equivalence of square matrices $\A$, the ``trace" of $\A$ is defined as
\begin{align}\label{4.2.10}
\Tr(\A)=\Tr(A),\quad A\in \A.
\end{align}
\end{enumerate}
\end{dfn}

Similar to Definition \ref{d4.2.7}, we need and can easily prove (\ref{4.2.10}) is well defined.
These two functions will be used in the sequel.

Next, we show that $\A=\{A_1,A_2,\cdots\}$ has a lattice structure.

\begin{dfn}[\cite{bur81}]\label{d2.3.8} A poset $L$ is a lattice if and only if for every pair $a,b\in L$ both $\sup\{a,b\}$ and $\inf\{a,b\}$ exist.
\end{dfn}

Let $A,B\in \A$. If $B$ is a divisor (multiple) of $A$, then $B$ is said to be proceeding (succeeding) $A$  and denoted by $B\prec A$ ( $B\succ A$). Then $\prec$ is a partial order for $\A$.

\begin{thm}\label{t2.3.9} $(\A,\prec)$ is a lattice.
\end{thm}

\noindent\textit{Proof}. Assume $A,B\in \A$. It is enough to prove that the $\Lambda=gcd(A,B)$ defined in (\ref{2.3.4}) is the $\inf(A,B)$, and the $\Theta=\lcm(A,B)$ defined in (\ref{2.3.601}) is the $\sup(A,B)$.

To prove $\Lambda=\inf(A,B)$ we assume $C\prec A$ and $C\prec B$, then we need only to prove that $C\prec \Lambda$.
Since $C\prec A$ and $C\prec B$, there exist $I_p$ and $I_q$ such that $C\otimes I_p=A$ and $C\otimes I_q=B$. Now
$$
\begin{array}{l}
C\otimes I_p=A=\Lambda\otimes I_{\b},\\
C\otimes I_q=B=\Lambda\otimes I_{\a}.
\end{array}
$$
Hence
$$
\begin{array}{ccl}
C\otimes I_p\otimes I_q&=&\Lambda\otimes I_{\b}\otimes I_{q}\\
~&=&\Lambda\otimes I_{\a}\otimes I_{p}.
\end{array}
$$
It follows that
$$
\b q=\a p.
$$
Since $\a$ and $\b$ are co-prime, we have $p=m\b$ and $q=n\a$, where $m,n\in \N$. Then we have
$$
C\otimes I_p=C\otimes I_m\otimes I_{\b}=\Lambda\otimes I_{\b}.
$$
That is
$$
C\otimes I_m=\Lambda.
$$
Hence, $C\prec \Lambda$.

To prove $\Theta=\sup(A,B)$ assume $D\succ A$ and $D\succ B$. Then we can prove $D\succ \Theta$ in a similar way.
\hfill $\Box$

\begin{dfn}\label{d2.3.10} $\A$ is said to possess a property, if every $A\in \A$ possesses this property. The property is also said to be consistent with the equivalence relation.
\end{dfn}

In the following some easily verifiable consistent properties are collected.

\begin{prp}\label{p2.3.11}
\begin{enumerate}
\item Assume $A\in {\cal M}$ is a square matrix. The following properties are consistent with the matrix equivalence ($\sim_{\ell}$ or $\sim_r$):
\begin{itemize}
\item $A$ is orthogonal, that is $A^{-1}=A^T$;
\item $\det(A)=1$;
\item $\tr(A)=0$;
\item $A$ is upper (lower) triangle;
\item $A$ is strictly upper (lower) triangle;
\item $A$ is symmetric (skew-symmetric);
\item $A$ is diagonal;
\end{itemize}
\item Assume $A\in {\cal M}_{2n\times 2n}$, $n=1,2,\cdots$, and
\begin{align}\label{2.3.8}
J=\begin{bmatrix}
0&1\\
-1&0
\end{bmatrix}
\end{align}
The following property is consistent with the matrix equivalence:
\begin{align}\label{2.3.9}
J\ltimes A+A^T\ltimes J=0.
\end{align}
\end{enumerate}
\end{prp}

\begin{rem}\label{r2.3.12} As long as a property is consistent with an equivalence, then we can say if the equivalent class has the property or not. For instance, because of Proposition \ref{p2.3.11} we can say $\A$ is orthogonal, $\det(\A)=1$, etc.
\end{rem}


\subsection{Lattice Structure on ${\cal M}_{\mu}$}

Denote by
\begin{align}\label{2.5.1}
{\cal M}_{\mu}:=\left\{A\in {\cal M}_{m\times n}\big| m/n=\mu\right\}.
\end{align}
Then it is clear that we have a partition as
\begin{align}\label{2.5.2}
{\cal M}=\bigcup_{\mu\in \Q_+}{\cal M}_{\mu},
\end{align}
where $\Q_+$ is the set of positive rational numbers.

\begin{rem}\label{r2.3.1201} To avoid possible confusion, we assume the fractions in $\Q_+$ are all reduced. Hence for each $\mu\in \Q_+$, there are unique integers $\mu_y$ and $\mu_x$, where $\mu_y$ and $\mu_x$ are co-prime, such that
$$
\mu=\frac{\mu_y}{\mu_x}.
$$
\end{rem}

\begin{dfn}\label{d2.7.0}
\begin{enumerate}
\item Let $\mu\in \Q_+$, $p$ and $q$ are co-prime and $p/q=\mu$. Then we denote by ${\mu}_y=p$ and ${\mu}_x=q$ as $y$ and $x$ are components of $\mu$.
\item Denote the spaces of various dimensions in ${\cal M}_{\mu}$ as
$$
{\cal M}_{\mu}^i:={\cal M}_{i\mu_y\times i\mu_x},\quad i=1,2,\cdots.
$$
\end{enumerate}
\end{dfn}

Assume $A_{\a}\in {\cal M}_{\mu}^{\a}$, $A_{\b}\in {\cal M}_{\mu}^{\b}$, $A_{\a}\sim A_{\b}$, and $\a|\b$, then $A_{\a}\otimes I_k=A_{\b}$, where $k=\b/\a$. One sees easily that we can define an embedding mapping $\bd_k:{\cal M}_{\mu}^{\a}\ra {\cal M}_{\mu}^{\b}$ as
\begin{align}\label{2.7.1}
\bd_k(A):=A\otimes I_k.
\end{align}
In this way, ${\cal M}_{\mu}^{\a}$ can be considered as a subspace of ${\cal M}_{\mu}^{\b}$. The order determined by this space-subspace relation is denoted as
\begin{align}\label{2.7.2}
{\cal M}_{\mu}^{\a}\sqsubset {\cal M}_{\mu}^{\b}.
\end{align}

If (\ref{2.7.2}) holds, ${\cal M}_{\mu}^{\a}$ is called a divisor of ${\cal M}_{\mu}^{\b}$, and ${\cal M}_{\mu}^{\b}$ is called a multiple of ${\cal M}_{\mu}^{\a}$.

Denote by $i\wedge j=gcd(i,j)$ and $i\vee j=lcm(i,j)$. Using the order of (\ref{2.7.2}),  ${\cal M}_{\mu}$  has the following structure.

\begin{figure}[!htbp]
\centering
\includegraphics[width=5cm]{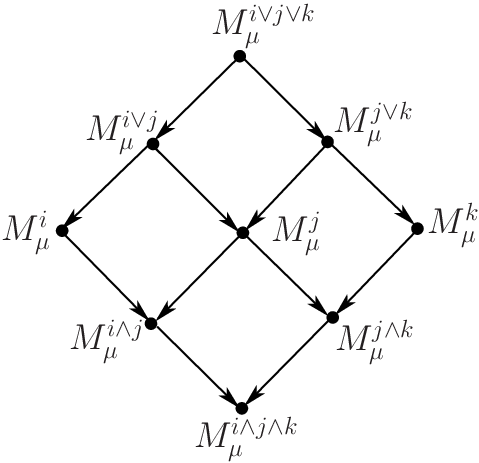}
\caption{Lattice Structure of ${\cal M}_{\mu}$}\label{Fig.2}
\end{figure}

\begin{thm}\label{t2.7.1}
\begin{enumerate}
\item Given ${\cal M}_{\mu}^i$ and ${\cal M}_{\mu}^j$. The greatest common divisor is ${\cal M}_{\mu}^{i\wedge j}$, and the least common multiple is ${\cal M}_{\mu}^{i\vee j}$. (Please refer to Fig.~\ref{Fig.2})

\item Assume $A\sim B$, $A\in {\cal M}_{\mu}^i$ and $B\in {\cal M}_{\mu}^j$. Then their greatest common divisor $\Lambda=gcd(A,B)\in {\cal M}_{\mu}^{i\wedge j}$, and their least common multiple $\Theta=lcm(A,B)\in {\cal M}_{\mu}^{i\vee j}$.
\end{enumerate}
\end{thm}

Next, we define
\begin{align}\label{2.7.3}
\begin{array}{l}
{\cal M}_{\mu}^i\wedge {\cal M}_{\mu}^j={\cal M}_{\mu}^{i\wedge j},\\
{\cal M}_{\mu}^i\vee {\cal M}_{\mu}^j={\cal M}_{\mu}^{i\vee j}.\\
\end{array}
\end{align}

From above discussion, the following result is obvious:

\begin{prp}\label{p2.7.2}
Consider ${\cal M}_{\mu}$. The followings are equivalent:
\begin{enumerate}
\item  ${\cal M}_{\mu}^{\a}$ is a subspace of  ${\cal M}_{\mu}^{\b}$;
\item $\a$ is a factor of $\b$, i.e., $\a|\b$;
\item
${\cal M}_{\mu}^{\a} \wedge {\cal M}_{\mu}^{\b}={\cal M}_{\mu}^{\a}$;
\item
${\cal M}_{\mu}^{\a} \vee {\cal M}_{\mu}^{\b}={\cal M}_{\mu}^{\b}$.
\end{enumerate}
\end{prp}

Using the order $\sqsubset$ defined by (\ref{2.7.2}), it is clear that all the fixed dimension vector spaces ${\cal M}_{\mu}^i$, $i=1,2,\cdots$ form a lattice.

\begin{prp}\label{p2.7.3}  $\left({\cal M}_{\mu},\sqsubset \right)$ is a lattice with
\begin{align}\label{2.7.4}
\begin{array}{ccl}
\sup\left({\cal M}_{\mu}^{\a},{\cal M}_{\mu}^{\b}\right)&=&{\cal M}_{\mu}^{\a\vee \b},\\
\inf\left({\cal M}_{\mu}^{\a},{\cal M}_{\mu}^{\b}\right)&=&{\cal M}_{\mu}^{\a\wedge \b}.\\
\end{array}
\end{align}
\end{prp}

The following properties are easily verifiable.

\begin{prp}\label{p2.7.3} Consider the lattice  $\left({\cal M}_{\mu},\sqsubset \right)$.
\begin{enumerate}
\item It has a smallest (root) subspace $ {\cal M}_{\mu}^1={\cal M}_{p\times q}$, where $p,q$ are co-prime and $p/q=\mu$. That is,
$$
\begin{array}{l}
{\cal M}_{\mu}^i\wedge {\cal M}_{\mu}^1= {\cal M}_{\mu}^1,\\
{\cal M}_{\mu}^i\vee {\cal M}_{\mu}^1={\cal M}_{\mu}^i.
\end{array}
$$
But there is no largest element.

\item The lattice is distributive, i.e.,
$$
\begin{array}{l}
{\cal M}_{\mu}^i\wedge \left({\cal M}_{\mu}^j\vee {\cal M}_{\mu}^k \right)=
\left({\cal M}_{\mu}^i\wedge {\cal M}_{\mu}^j\right) \vee \left({\cal M}_{\mu}^i\wedge {\cal M}_{\mu}^k \right),\\
{\cal M}_{\mu}^i\vee \left({\cal M}_{\mu}^j\wedge {\cal M}_{\mu}^k \right)=
\left({\cal M}_{\mu}^i\vee {\cal M}_{\mu}^j\right) \wedge \left({\cal M}_{\mu}^i\vee {\cal M}_{\mu}^k \right).
\end{array}
$$
\item For any finite set of spaces ${\cal M}_{\mu}^{i_s}$, $s=1,2,\cdots,r$. There exists a smallest supper-space ${\cal M}_{\mu}^u$, $u=\vee_{s=1}^ri_s$, such that
\begin{itemize}
\item
    $$
    {\cal M}_{\mu}^{i_s}\sqsubset {\cal M}_{\mu}^u,\quad s=1,2,\cdots,r;
    $$
\item
If
$${\cal M}_{\mu}^{i_s}\sqsubset {\cal M}_{\mu}^v,\quad s=1,2,\cdots,r,
$$
then
$${\cal M}_{\mu}^{u}\sqsubset {\cal M}_{\mu}^v.
$$
\end{itemize}
\end{enumerate}

\end{prp}

\begin{dfn}[\cite{fan14}]\label{d2.7.4} Let $(L,\prec)$ and $(M,\sqsubset)$ be two lattices.
\begin{enumerate}
\item[(1)] A mapping $\varphi:L\ra M$ is called an order-preserving mapping, if $\ell_1\prec \ell_2$ implies $\varphi(\ell_1)\sqsubset \varphi(\ell_2)$.
\item[(2)] A mapping $\varphi:L\ra M$ is called a homomorphism, and $(L,\prec)$ and $(M,\sqsubset)$ are said to be  lattice homomorphic, denoted by $(L,\prec)\approx (M,\sqsubset)$, if $\varphi$ satisfies the following condition:
\begin{align}\label{7.5}
\varphi \sup(\ell_1,\ell_2)=\sup\left(\varphi(\ell_1),\varphi(\ell_2)\right);
\end{align}
and
\begin{align}\label{7.501}
\varphi \inf(\ell_1,\ell_2)=\inf\left(\varphi(\ell_1),\varphi(\ell_2)\right).
\end{align}
\item[(3)] A homomorphism $\varphi:L\ra M$ is called an isomorphism, and $(L,\prec)$ and $(M,\sqsubset)$ are said to be  lattice isomorphic, denoted by $(L,\prec)\approxeq (M,\sqsubset)$, if $\varphi$ is one to one and onto.
\end{enumerate}
\end{dfn}

Assume $A\in {\cal M}_{\mu}^i$ is irreducible, define $\varphi:\A\ra {\cal M}_{\mu}$ as
\begin{align}\label{2.7.6}
\varphi(A_j):={\cal M}_{\mu}^{ij},\quad j=1,2,\cdots,
\end{align}
then it is easy to verify the following result.

\begin{prp}\label{p2.7.5} The mapping $\varphi:\A\ra {\cal M}_{\mu}$ defined in (\ref{2.7.6}) is a lattice homomorphism from $(\A,\prec)$ to $({\cal M}_{\mu},\sqsubset)$.
\end{prp}

Next, we consider the ${\cal M}_{\mu}$ for different $\mu$'s. It is also easy to verify the following result.

\begin{prp}\label{p2.7.6} Define a mapping $\varphi:{\cal M}_{\mu}\ra {\cal M}_{\lambda}$ as
$$
\varphi\left({\cal M}_{\mu}^i\right):={\cal M}_{\lambda}^i.
$$ The mapping $\varphi:\left({\cal M}_{\mu},\sqsubset\right)\ra \left({\cal M}_{\lambda},\sqsubset\right)$
is a lattice isomorphism.
\end{prp}

\begin{exa}\label{e2.7.601} According to Proposition \ref{p2.7.6}, if we still assume $A\in {\cal M}_{\mu}^i$ and replace $\mu$ in (\ref{2.7.6}) by any $\a\in \Q_+$, that is,  define
$$
\varphi(A_j):={\cal M}_{\a}^{ij},\quad j=1,2,\cdots,
$$
then it is easy to see that  $\varphi:~(\A,\prec) \ra ({\cal M}_{\a},\sqsubset)$ is still a lattice homomorphism.
\end{exa}

\begin{dfn}\label{d2.7.7} Let $(L,\prec)$ be a lattice and $S\subset L$. If $(S,\prec)$ is also a lattice, it is called a sub-lattice of $(L,\prec)$.
\end{dfn}

\begin{rem}\label{r2.7.701} Let $\varphi: (H,\prec)\ra (M,\sqsubset)$ be an injective (i.e., one-to-one) lattice homomorphism. Then $\varphi:H\ra\varphi(H)$ is a lattice isomorphism. Hence $\varphi(H)$ is a sub-lattice of $(M,\sqsubset)$. If we identify $H$ with $\varphi(H)$, we can simply say that $H$ is a sub-lattice of $M$.
\end{rem}

\begin{dfn}\label{d2.7.8}  Let $(L,\prec)$ and $(M,\sqsubset)$ be two lattices. The product order $\subset:=\prec\times \sqsubset$ defined on the product set
$$
L\times M:=\left\{(\ell,m)\;\big|\; \ell\in L,m\in M\right\}
$$
is:
$(\ell_1,m_1)\subset (\ell_2,m_2)$ if and only if $\ell_1\prec \ell_2$ and $m_1\sqsubset m_2$.
\end{dfn}

\begin{thm}\label{t2.7.9} Let $(L,\prec)$ and $(M,\sqsubset)$ be two lattices. Then $\left(L\times M,\prec\times \sqsubset \right)$ is also a lattice, called the product lattice of $(L,\prec)$ and $(M,\sqsubset)$.
\end{thm}

\noindent\textit{Proof}. Let $(\ell_1,~m_1)$ and $(\ell_2,~m_2)$ be two elements in $L\times M$. Denote by $\ell_s=\sup(\ell_1,~\ell_2)$ and $m_s=\sup(m_1,~m_2)$. Then $(\ell_s,~m_s)\supset (\ell_j,~m_j)$, $j=1,2$. To see
$(\ell_s,~m_s)=\sup\left( (\ell_1,~m_1),~(\ell_2,~m_2)\right)$ let $(\ell,~m)\supset (\ell_j,~m_j)$, $j=1,2$. Then
$\ell\succ \ell_j$ and $m\sqsupset m_j$, $j=1,2$. It follows that $\ell\succ \ell_s$ and $m\sqsupset m_s$. That is, $(\ell,~m)\supset (\ell_s,~m_s)$. We conclude that
$$(\ell_s,~m_s)=\sup\left((\ell_1,~m_1),~(\ell_2,~m_2)\right).
$$

Similarly, we set $\ell_i=\inf(\ell_1,~\ell_2)$ and $m_i=\inf(m_1,~m_2)$, then we can prove that
$$(\ell_i,~m_i)=\inf\left((\ell_1,~m_1),~(\ell_2,~m_2)\right).
$$
\hfill $\Box$

Finally, we give an example to show that an order-preserving mapping may not be an lattice homomorphism.

\begin{exa}\label{e2.7.10} Consider the product of two lattices $\left({\cal M}_{\mu},\sqsubset\right)$ and $\left({\cal M}_{\lambda},\sqsubset\right)$. Define a mapping
$$
\varphi:\left({\cal M}_{\mu},\sqsubset\right)\times \left({\cal M}_{\lambda},\sqsubset\right)\ra \left({\cal M}_{\mu\lambda},\sqsubset\right)
$$
as
$$
\varphi\left( {\cal M}_{\mu}^p\times{\cal M}_{\lambda}^q\right):={\cal M}^{pq}_{\mu\lambda}.
$$

Assume ${\cal M}_{\mu}^i\sqsubset {\cal M}_{\mu}^j$ and ${\cal M}_{\lambda}^s\sqsubset {\cal M}_{\lambda}^t$, then $i|j$ and $s|t$, and by the definition of product lattice, we have
$$
{\cal M}_{\mu}^i\times {\cal M}_{\lambda}^s \sqsubset\times \sqsubset {\cal M}_{\mu}^j\times {\cal M}_{\lambda}^t.
$$

Since $is|jt$, we have
\begin{align}\label{2.7.7}
\begin{array}{l}
\varphi\left({\cal M}_{\mu}^i\times {\cal M}_{\lambda}^s\right)= {\cal M}_{\mu\lambda}^{is},\\
\sqsubset {\cal M}_{\mu\lambda}^{jt}=\varphi\left({\cal M}_{\mu}^j\times {\cal M}_{\lambda}^t\right).
\end{array}
\end{align}
That is, $\varphi$ is an order-preserving mapping.

Consider two elements in product lattice as $\a={\cal M}_{\mu}^p\times {\cal M}_{\lambda}^s$ and $\b={\cal M}_{\mu}^q\times {\cal M}_{\lambda}^t$. Following the same arguments in the proof of Theorem \ref{t2.7.9}, one sees easily that
$$
\begin{array}{l}
\lcm(\a,~\b)={\cal M}_{\mu}^{p\vee q}\times {\cal M}_{\lambda}^{s\vee t},\\
\gcd(\a,~\b)={\cal M}_{\mu}^{p\wedge q}\times {\cal M}_{\lambda}^{s\wedge t}.\\
\end{array}
$$
Then
$$
\begin{array}{l}
\varphi(\lcm(\a,~\b))={\cal M}_{\mu\lambda}^{(p\vee q)(s\vee t)},\\
\varphi(\gcd(\a,~\b))={\cal M}_{\mu\lambda}^{(p\wedge q)(s\wedge t)}.\\
\end{array}
$$
Consider
$$
\begin{array}{l}
\varphi(\a)={\cal M}_{\mu\lambda}^{ps},\\
\varphi(\b)={\cal M}_{\mu\lambda}^{qt}.\\
\end{array}
$$
Now
$$
\begin{array}{l}
\lcm(\varphi(\a),\varphi(\b))={\cal M}_{\mu\lambda}^{(ps)\vee(qt)},\\
\gcd(\varphi(\a),\varphi(\b))={\cal M}_{\mu\lambda}^{(ps)\wedge(qt)}.\\
\end{array}
$$
It is obvious that in general
$$(p\vee q)(s\vee t)\neq (ps)\vee(qt),$$
as well as,
$$(p\wedge q)(s\wedge t)\neq (ps)\wedge(qt).
$$
Hence, $\varphi$ is not a homomorphism.
\end{exa}

\subsection{Monoid and Quotient Monoid}

A monoid is a semigroup with identity. We refer readers to \cite{jac85}, \cite{how95}, \cite{fou82} for concepts and some basic properties.

Recall that
$$
{\cal M}:=\bigcup_{m\in \N}\bigcup_{n\in \N}{\cal M}_{m\times n}.
$$
We have the following algebraic structure.

\begin{prp}\label{p2.4.1} The algebraic system $\left( {\cal M},\ltimes \right)$ is a monoid.
\end{prp}

\noindent\textit{Proof}. The associativity comes from the property of $\ltimes$ (refer to (\ref{2.2.4})). The identity element is $1$.
\hfill $\Box$

One sees easily that this monoid covers numbers, vectors, and matrices of arbitrary dimensions.

In the following some of its useful sub-monoids are presented:

\begin{itemize}
\item ${\cal M}(k)$:
$$
{\cal M}(k):=\bigcup_{\a\in \N}\bigcup_{\b\in \N}{\cal M}_{k^{\a}\times k^{\b}},
$$
where $k\in \N$ and $k>1$.

It is obvious that ${\cal M}(k) < {\cal M}$. (In this section $A<B$ means $A$ is a submonoid of $B$). This sub-monoid is useful for calculating the product of tensors over $k$ dimensional vector space \cite{boo79}. It is particularly useful for $k$-valued logical dynamic systems \cite{che11}, \cite{che12}. When $k=2$, it is used for Boolean dynamic systems.

In this sub-monoid the STP can be defined as follows:
\end{itemize}

\begin{dfn}\label{d2.4.2}
\begin{enumerate}
\item Let $X\in \F^n$ be a column vector,  $Y\in \F^m$ a row vector.
\begin{itemize}
\item Assume $n=pm$ (denoted by $X\succ_p Y$):
Split $X$ into $m$ equal blocks as
$$
X=\left[X_1^T,X_2^T,\cdots,X_m^T\right]^T,
$$
where $X_i\in \F^p$, $\forall i$. Define
$$
X\ltimes Y:=\dsum_{s=1}^mX_sy_s\in \F^p.
$$
\item Assume $np=m$ (denoted by $X\prec_p Y$):
Split $Y$ into $n$ equal blocks as
$$
Y=\left[Y_1,Y_2,\cdots,Y_n\right],
$$
where $Y_i\in \F^p$, $\forall i$. Define
$$
X\ltimes Y:=\dsum_{s=1}^mx_sY_s\in \F^p.
$$
\end{itemize}

\item Assume $A\in {\cal M}_{m\times n}$,  $B\in {\cal M}_{p\times q}$, where
$n=tp$ (denoted by $A\succ_t B$), or $nt=p$ (denoted by $A\prec_t B$). Then
$$
A\ltimes B:=C=\left(c_{i,j}\right),
$$
where
$$
c_{i,j}=\Row_i(A)\ltimes \Col_j(B).
$$
\end{enumerate}
\end{dfn}

\begin{rem}\label{r2.4.3}
\begin{enumerate}
\item It is easy to prove that when $A\prec_t B$ or $B\prec_t A$ for some $t\in \N$, this definition of left STP coincides with Definition \ref{d2.2.1}. Though this definition is not as general as Definition \ref{d2.2.1}, it has clear physical meaning.
Particularly, so far this definition covers almost all the applications.
\item Unfortunately, this definition is not suitable for right STP. This is a big difference between left and right STPs.
\end{enumerate}
\end{rem}

\begin{itemize}
\item ${\cal V}$:
$$
{\cal V}:=\bigcup_{k\in \N} {\cal M}_{k\times 1}.
$$

It is obvious that ${\cal V} < {\cal M}$.

This sub-monoid consists of column vectors. In this sub-monoid the STP is degenerated to Kronecker product.
\end{itemize}

We denote by ${\cal V}^T$ the sub-monoid of row vectors. It is also clear that  ${\cal V}^T < {\cal M}$.

\begin{itemize}
\item ${\cal L}$:
$$
{\cal L}:=\left\{A\in {\cal M}\;|\; \Col(A)\subset \D_s,\;s\in \N\right\}.
$$

It is obvious that ${\cal L} < {\cal M}$. This sub-monoid consists of all logical matrices. It is used to express the product of logical mappings.

\end{itemize}

\begin{itemize}
\item ${\cal P}$:
$$
{\cal P}:=\left\{A\in {\cal M}\;|\; \Col(A)\subset \varUpsilon_s,\; \mbox{for some}~s\in \N\right\}.
$$

It is obvious that ${\cal P} < {\cal M}$. This monoid is useful for probabilistic logical mappings.

\end{itemize}

\begin{itemize}
\item ${\cal L}(k)$:
$$
{\cal L}(k):={\cal L} \cap {\cal M}(k).
$$

It is obvious that ${\cal L}(k)< {\cal L} < {\cal M}$. We use it for $k$-valued logical mappings.

\end{itemize}

Next, we define the set of ``short" matrices as
$$
\varXi:=\left\{A\in {\cal M}_{m\times n}\;|\; m\leq n\right\},
$$
and its subset
$$
{\varXi}^r:=\left\{A\in {\cal M}\;|\; A \mbox{~is of full row rank}\right\}.
$$

Then we have the following result.

\begin{prp}\label{p2.4.3}
\begin{align}\label{2.4.1}
{\varXi}^r<{\varXi}<{\cal M}.
\end{align}
\end{prp}

\noindent\textit{Proof}. Assume $A\in {\cal M}_{m\times n}$, $B\in {\cal M}_{p\times q}$ and $A,B\in {\varXi}$, then $m\leq n$ and $p\leq q$. Let $t=n\vee p$. Then $AB\in {\cal M}_{\frac{mt}{n}\times \frac{tq}{p}}$.
It is easy to see that $\frac{mt}{n}\leq \frac{tq}{p}$, so $AB\in {\varXi}$. The second part is proved.

As for the first part, Assume $A,~B\in {\varXi}^r$. Then
$$
\begin{array}{ccl}
\rank(AB)&=&\rank\left[\left(A\otimes I_{t/n}\right)\left(B\otimes I_{t/p}\right)\right]\\
~&\geq&\rank\left[\left(A\otimes I_{t/n}\right)\left(B\otimes I_{t/p}\right)\right.\\
~&~&\left.\left(B^T(BB^T)^{-1}\otimes  I_{t/p}\right)\right]\\
~&=&\rank\left[\left(A\otimes I_{t/n}\right)\left(I_p\otimes I_{t/p}\right)\right]\\
~&=&\rank\left(A\otimes I_{t/n}\right)=mt/n.
\end{array}
$$
Hence, $AB\in {\varXi}^r$.
\hfill $\Box$

Similarly, we can define the set of ``tall" matrices ${\Pi}$ and the set of matrices with full column rank
${\Pi}^c$. We can also prove that
\begin{align}\label{2.4.2}
{\Pi}^c<{\Pi}<{\cal M}.
\end{align}

Next, we consider the quotient space
$$
\Sigma_{{\cal M}}:={\cal M}/\sim.
$$

\begin{dfn}[\cite{rad89}]\label{d2.4.4}
\begin{enumerate}
\item A nonempty set $S$ with a binary operation $\sigma:S\times S\ra S$ is called an algebraic system.
\item Assume $\sim$ is an equivalence relation on an algebraic system $(S,\sigma)$. The equivalence relation is a congruence relation, if for any $A,B,C,D\in S$, $A\sim C$ and $B\sim D$, we have
\begin{align}\label{2.4.3}
A\sigma B \sim C\sigma D.
\end{align}
\end{enumerate}
\end{dfn}

\begin{prp}\label{p2.4.5}
Consider the algebraic system $\left({\cal M}, \ltimes \right)$ with the equivalence relation $\sim=\sim_{\ell}$.  The equivalence relation $\sim$ is congruence.
\end{prp}

\noindent\textit{Proof.} Let $A\sim \tilde{A}$ and $B\sim \tilde{B}$. According to Theorem \ref{t2.3.5}, there exist $U\in {\cal M}_{m\times n}$ and $V\in {\cal M}_{p\times q}$ such that
$$
\begin{array}{ll}
A=U\otimes I_s,& \tilde{A}=U\otimes I_t;\\
B=V\otimes I_{\a},& \tilde{B}=V\otimes I_{\b}.\\
\end{array}
$$
Denote
$$
n\vee p=r,\quad ns\vee \a p=r\xi,\quad nt\vee \b p=r\eta.
$$
Then
$$
\begin{array}{l}
A\ltimes B=\left(U\otimes I_{s}\otimes I_{r\xi/ns}\right)\left(V\otimes I_{\a}\otimes I_{r\xi/\a p}\right)\\
=\left[\left(U\otimes I_{r/n}\right)\left(V\otimes I_{r/p}\right)\right]\otimes I_{\xi}.
\end{array}
$$
Similarly, we have
$$
\tilde{A}\ltimes \tilde{B}
=\left[\left(U\otimes I_{r/n}\right)\left(V\otimes I_{r/p}\right)\right]\otimes I_{\eta}.
$$
Hence we have $A\ltimes B\sim \tilde{A}\ltimes \tilde{B}$.
\hfill $\Box$

According to Proposition \ref{p2.4.5}, we know that $\ltimes$ is well defined on the quotient space $\Sigma_{{\cal M}}$. Moreover, the following result is obvious:

\begin{prp}\label{p2.4.6}
\begin{enumerate}
\item $\left(\Sigma_{{\cal M}}, \ltimes\right)$ is a monoid.

\item Let ${\cal S}< {\cal M}$ be a sub-monoid. Then ${\cal S}/\sim$ is a sub-monoid of $\Sigma_{{\cal M}}$, that is,
$$
{\cal S}/\sim ~<~ \Sigma_{{\cal M}}.
$$
\end{enumerate}
\end{prp}

Since the $S$ in Proposition \ref{p2.4.6} could be any sub-monoid of ${\cal M}$.
All the aforementioned sub-monoids have their corresponding quotient sub-monoids, which are the sub-monoids of $\Sigma_{{\cal M}}$. For instance, ${\cal V}/\sim$, ${\cal L}/\sim$, etc. are the sub-monoids of $\Sigma_{{\cal M}}$.
\subsection{Group Structure on ${\cal M}^{\mu}$}

\begin{prp}\label{p2.5.1} Assume $A\in {\cal M}_{\mu_1}$ and $B\in {\cal M}_{\mu_2}$ then $A\ltimes B\in {\cal M}_{\mu_1\mu_2}$. That is, the operation $\ltimes$ is a mapping
$$
\ltimes:~{\cal M}_{\mu_1}\times {\cal M}_{\mu_2}\ra  {\cal M}_{\mu_1\mu_2}.
$$
\end{prp}

\noindent\textit{Proof}. Assume $A\in {\cal M}_{m\times n}$ and $B\in {\cal M}_{p\times q}$, where $\mu_1=m/n$ and $\mu_2=p/q$, and $t=n\vee p$. Then
$$
\begin{array}{ccl}
A\ltimes B&=&\left(A\otimes I_{t/n}\right)\left(B\otimes I_{t/p}\right)\\
~&\in&{\cal M}_{mt/n\times qt/p}\subset {\cal M}_{\mu_1 \mu_2}.
\end{array}
$$
\hfill $\Box$

\begin{dfn}\label{d2.5.3}
\begin{enumerate}
\item Define
$$
{\cal M}^{\mu}:=\bigcup_{z\in \Z}{\cal M}_{\mu^z}.
$$
Then ${\cal M}^{\mu}$ is closed under operator $\ltimes$.
\item Set
$$
\Sigma_{\mu^z}={\cal M}_{\mu^z}/\sim,
$$
and define
$$
\Sigma^{\mu}:=\bigcup_{z\in \Z}\Sigma_{\mu^z}.
$$
Then $\Sigma^{\mu}$ is also closed under operator $\ltimes$.
\item $\A,\B\in \Sigma^{\mu}$ is said to be power equivalent, denoted by $\A\sim_p \B$, if there exists an integer $z\in \Z$ such that both $\A,~\B\in \Sigma_{\mu^z}$. Denote
\begin{align}\label{2.5.3}
\left<\left<A\right>\right>:=\left\{\B\;|\; \B\sim_p \A \right\}
\end{align}
\end{enumerate}
\end{dfn}

\begin{rem}\label{r2.5.4}
It is obvious that $\ltimes$ is consistent with $\sim_p$. Hence $\ltimes$ is well defined on the set of equivalent classes as
\begin{align}\label{2.5.4}
\left<\left<A\right>\right> \ltimes \left<\left<B\right>\right>:=\left<\left<A\ltimes B\right>\right>.
\end{align}
\end{rem}

Then we have the following group structure.
\begin{thm}\label{t2.5.5} $\left(\Sigma^{\mu}/\sim_p, \ltimes\right)$ is a group, which is isomorphic to $(\Z,+)$. Precisely, assume $A\in {\cal M}_{\mu^z}$ then $\Psi: \Sigma^{\mu}/\sim_p \ra \Z$ is defined as
$$
\Psi\left(\left<\left< A\right>\right>\right):=z,
$$
which is a group isomorphism.
\end{thm}

\subsection{Semi-tensor Addition and Vector Space Structure of $\Sigma_{\mu}$}

\begin{dfn}\label{d2.6.1} Let $A,~B\in {\cal M}_{\mu}$. Precisely, $A\in {\cal M}_{m\times n}$, $B\in {\cal M}_{p\times q}$, and $m/n=p/q=\mu$. Set $t=m\vee p$. Then
\begin{enumerate}
\item the left semi-tensor addition (STA) of $A$ and $B$, denote by
$\lplus$, is defined as
\begin{align}\label{2.6.2}
A~\lplus~ B:=\left(A\otimes I_{t/m}\right)+\left(B\otimes I_{t/p}\right).
\end{align}
Correspondingly, the left semi-tensor subtraction (STS) is defined as
\begin{align}\label{2.6.3}
A\lminus B:=A~\lplus~(-B).
\end{align}
\item The right STA of $A$ and $B$, denote by $\rplus$, is defined as
\begin{align}\label{2.6.4}
A~\rplus ~B:=\left(I_{t/m}\otimes A\right)+\left(I_{t/p}\otimes B\right).
\end{align}
Correspondingly, the right STS is defined as
\begin{align}\label{2.6.5}
A\rminus B:=A\rplus (-B).
\end{align}
\end{enumerate}
\end{dfn}

\begin{rem}\label{r2.6.2} Let $\sigma\in \{\lplus, \lminus, \rplus, \rminus\}$ be one of the four binary operators. Then it is easy to verify that
\begin{enumerate}
\item if $A,~B\in {\cal M}_{\mu}$, then $A\sigma B\in {\cal M}_{\mu}$;
\item If $A$ and $B$ are as in Definition \ref{d2.6.1}, then $A\sigma B\in {\cal M}_{t\times \frac{nt}{m}}$;
\item Set $s=n\vee q$, then $s/n=t/m$ and $s/q=t/p$. So $\sigma$ can also be defined by using column numbers respectively, e.g.,
$$
A\lplus B:=\left(A\otimes I_{s/n}\right)+\left(B\otimes I_{s/q}\right),
$$
etc.
\end{enumerate}
\end{rem}


\begin{thm}\label{t2.6.4}
 Consider the algebraic system $\left({\cal M}_{\mu}, \sigma \right)$, where $\sigma\in \{\lplus, \lminus\}$ and $\sim=\sim_{\ell}$ (or $\sigma\in \{\rplus, \rminus\}$ and  $\sim=\sim_r$). Then the equivalence relation $\sim$ is a congruence relation with respect to $\sigma$.
\end{thm}

\noindent\textit{Proof}. We prove one case, where $\sigma=\lplus$ and $\sim=\sim_{\ell}$. Proofs for other cases are similar.

Assume $\tilde{A}\sim_{\ell} A$ and $\tilde{B}\sim_{\ell} B$. Set $P=gcd(\tilde{A},A)$ and $Q=gcd(\tilde{B},B)$, then
\begin{align}\label{2.6.7}
\tilde{A}=P\otimes I_{\b},\quad
A=P\otimes I_{\a};\\
\label{2.6.8}
\tilde{B}=Q\otimes I_{\gamma},\quad
B=Q\otimes I_{\delta},
\end{align}
 where $P\in {\cal M}_{x\mu\times x}$, $Q\in {\cal M}_{y\mu\times y}$, $x,~y\in \N$ are certain numbers.

Now consider $\tilde{A}\lplus \tilde{B}$. Assume $\eta=x\vee y$, $t=x\b\vee y\gamma=\eta\xi$, $s=x\a\vee y\delta=\eta\zeta$. Then we have
\begin{align}\label{2.6.9}
\begin{array}{ccl}
\tilde{A}\lplus \tilde{B}&=&P\otimes I_{\b}\otimes I_{t/x\b}\\
~&~&+Q\otimes I_{\gamma}\otimes I_{t/y\gamma}\\
~&=&\left[(P\otimes I_{\eta/x})+(Q\otimes I_{\eta/y)}\right]\otimes I_{\xi}.
\end{array}
\end{align}
Similarly, we have
\begin{align}\label{2.6.10}
A\lplus B=\left[(P\otimes I_{\eta/x})+(Q\otimes I_{\eta/y})\right]\otimes I_{\zeta}.
\end{align}
(\ref{2.6.9}) and (\ref{2.6.10}) imply that $
\tilde{A}\lplus \tilde{B}\sim A\lplus B$.
\hfill $\Box$

Define the left and right quotient spaces $\Sigma^{\ell}_{\mu}$ and $\Sigma^{r}_{\mu}$ respectively as
\begin{align}\label{2.6.11}
\Sigma^{\ell}_{\mu}&:={\cal M}_{\mu}/\sim_{\ell};\\
\label{2.6.12}
\Sigma^{r}_{\mu}&:={\cal M}_{\mu}/\sim_{r}.
\end{align}

According to Theorem \ref{t2.6.4}, the operation $\lplus$ (or $\lminus$) can be extended to $\Sigma^{\ell}_{\mu}$
as
\begin{align}\label{2.6.13}
\begin{array}{l}
\A_{\ell} \lplus \B_{\ell}:=<A~\lplus~ B>_{\ell},\\
\A_{\ell} \lminus \B_{\ell}:=<A\lminus B>_{\ell},
\quad \A_{\ell},~\B_{\ell}\in \Sigma^{\ell}_{\mu}.
\end{array}
\end{align}

Similarly,  we can define $\rplus$ (or $\rminus$) on the quotient space$\Sigma^{r}_{\mu}$ as
\begin{align}\label{2.6.14}
\begin{array}{l}
\A_{r} \rplus \B_{r}:=<A\rplus B>_{r},\\
\A_{r} \rminus \B_{r}:=<A\rminus B>_{r},\quad \A_{r},\B_{r}\in \Sigma^{r}_{\mu}.
\end{array}
\end{align}

The following result is important, and the verification is straightforward.

\begin{thm}\label{t2.6.5} Using the definitions in (\ref{2.6.13}) (correspondingly, (\ref{2.6.14})), the quotient space
$\left(\Sigma^{\ell}_{\mu}, \lplus\right)$ (correspondingly, $\left(\Sigma^{r}_{\mu}, \rplus\right)$ )
is a vector space.
\end{thm}

\begin{rem}\label{r2.6.6} As a consequence, $\left( \Sigma^{\ell}_{\mu}, \lplus \right)$ (or $\left( \Sigma^{r}_{\mu}, \rplus \right)$)  is an Abelian group.
\end{rem}

\begin{rem}\label{r2.3.13} Recall Example \ref{e2.2.11}, it shows that the exponential function $\exp$ is well defined on the quotient space $\Sigma:=\Sigma_1$.
\end{rem}

Since each $\A\in \Sigma$ has a unique left (or right) irreducible element $A_1$ (or $\tilde{A}_1$) such that $A\sim_{\ell} A_1$ (or $A\sim_{r} \tilde{A}_1$), in general, we can use the irreducible element, which is also called the root element of an equivalent class, as the representation of this class. But this is not compulsory.

For notational and statement ease, hereafter we consider $\Sigma^{\ell}_{\mu}$ only unless elsewhere stated. As a convention, the omitted script ( ``$\ell$" or``$r$ ") means $\ell$. For instance, $\Sigma_{\mu}=\Sigma^{\ell}_{\mu}$, $\sim=\sim_{\ell}$, $\A=\A_{\ell}$ etc.


\section{Topology on M-equivalence Space}

\subsection{Topology via Sub-basis}

This subsection builds step by step a topology on quotient space $\Sigma_{\mu}$ using a sub-basis.

First, we consider the partition (\ref{2.5.2}),
it is natural to assume that each ${\cal M}_{\mu}$ is a clopen subset in ${\cal M}$, because distinct $\mu$'s correspond to distinct shapes of matrices. Now inside each ${\cal M}_{\mu}$ we assume $\mu_y,~\mu_x\in \N$ are co-prime and $\mu_y/\mu_x=\mu$. Then
$$
{\cal M}_{\mu}=\bigcup_{i=1}^{\infty}{\cal M}_{\mu}^i,
$$
where
$$
{\cal M}_{\mu}^i={\cal M}_{i\mu_y\times i\mu_x},\quad i=1,2,\cdots.
$$
Because of the similar reason, we also assume each ${\cal M}_{\mu}^i$ is clopen.

Overall, we have a set structure on ${\cal M}$ as
\begin{align}\label{0.2.1}
{\cal M}=\bigcup_{\mu\in \Q_+}\bigcup_{i=1}^{\infty}{\cal M}_{\mu}^{i}.
\end{align}

\begin{dfn}\label{d0.2.1} A natural topology on ${\cal M}$, denoted by ${\cal T}_{{\cal M}}$,  consists of
\begin{enumerate}
\item a partition of countable clopen subsets ${\cal M}_{\mu}^i$, $\mu\in \Q_+$, $i\in \N$;
\item the conventional Euclidean  $\R^{i^2\mu_y\mu_x}$ topology for ${\cal M}_{\mu}^i$.
\end{enumerate}
\end{dfn}
Then ${\cal M}$ becomes a topological space. Moreover, it is obvious that $\left({\cal M},{\cal T}_{{\cal M}}\right)$ is a second countable Hausdorff space.

Next, we consider the quotient space
$$
\Sigma_{{\cal M}}:={\cal M}/\sim.
$$
It is clear that
\begin{align}\label{0.2.2}
\Sigma_{{\cal M}}=\bigcup_{\mu\in \Q_+}\Sigma_{\mu}.
\end{align}
Moreover, (\ref{0.2.2}) is also a partition. Hence each $\Sigma_{\mu}$ can be considered as a clopen subset in $\Sigma_{{\cal M}}$. We are, therefore, interested only in constructing a topology on each $\Sigma_{\mu}$.

\begin{dfn}\label{d0.2.2}
\begin{enumerate}
\item Consider ${\cal M}^i_{\mu}$ as an Euclidean space $\R^{i^2\mu_y\mu_x}$ with conventional Euclidean topology. Assume $o_i\neq \emptyset$ is an open set. Define a subset $s_i(o_i)\subset \Sigma_{\mu}$ as follows:
\begin{align}\label{0.2.3}
\A\in s_i(o_i)\Leftrightarrow \A\cap o_i\neq \emptyset.
\end{align}
\item
Let
$$
\begin{array}{ccl}
O_i&=&\{o_i\;|\; o_i ~\mbox{is an open ball in}~ {\cal M}^i_{\mu}\\
 ~&~&~\mbox{with rational center and rational radius}~\}.
\end{array}
$$
\item
Using $O_i$, we construct a set of subsets $S_i\subset 2^{\Sigma_{\mu}}$ as
$$
\begin{array}{ccl}
S_i&:=&\left\{s_i\;|\;s_i=s_i(o_i) ~\mbox{for some}~ o_i\in O_i\right\},\\
~&~&\quad i=1,2,\cdots.
\end{array}
$$
 Taking $S=\cup_{i=1}^{\infty}S_i$ as a topological sub-basis, the topology generated by $S$ is denoted by ${\cal T}$, which makes
 $$
 \left(\Sigma_{\mu},{\cal T}\right)
 $$
 a topological space.  (We refer to \cite{kel75} for a topology produced from a sub-basis.)

 Note that the topological basis consists of the set of finite intersections of $s_i\in S_i$, $i=1,\cdots,t$, $t<\infty$.
\end{enumerate}
\end{dfn}

\begin{figure}[!htbp]
\centering
\includegraphics[width=3cm]{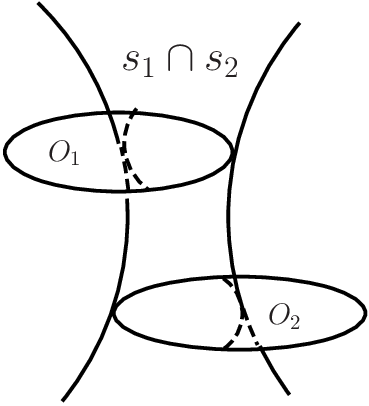}
\caption{$s_i\cap s_j$: An Element In Topological Basis}\label{Fig.0.1}
\end{figure}

\begin{rem}\label{r0.2.3}
\begin{enumerate}
\item It is clear that ${\cal T}$ makes $\left(\Sigma_{\mu},{\cal T}\right)$ a topological space.
\item The topological basis is
\begin{align}\label{0.2.4}
\begin{array}{ccl}
{\cal B}&:=&\left\{s_{i_1}\cap s_{i_2}\cap \cdots\cap s_{i_r}\;|\; s_{i_j}\in S_{i_j};\right.\\
~&~&\left.\;j=1,\cdots,r;\; r<\infty\right\}.
\end{array}
\end{align}
\item Fig.~\ref{Fig.0.1} depicts an element in the topological basis. Here $o_1\in {\cal M}_{\mu}^i$, $o_2\in {\cal M}_{\mu}^j$ are two open discs with  rational center and rational radius. Then
$s_1(o_1)$ and $s_2(o_2)$ are two elements in the sub-basis, and $$
s_1\cap s_2=\left\{\A\;\big|\; \A\cap o_i\neq \emptyset,~i=1,2\right\}
$$
is an element in the basis.
\end{enumerate}
\end{rem}

\begin{thm}\label{t0.2.4} The topological space $\left(\Sigma_{\mu},{\cal T}\right)$ is a second countable, Hausdorff (or $T_2$) space.
\end{thm}

\noindent\textit{Proof}. To see $\left(\Sigma_{\mu},{\cal T}\right)$ is second countable, It is easy to see that
$O_i$ is countable. Then $\{O_i|i=1,2,\cdots\}$, as countable union of countable set, is countable. Finally, ${\cal B}$, as the finite subset of a countable set, is countable.

Next, consider $\A\neq\B\in \Sigma_{\mu}$. Let $A_1\in \A$ and $B_1\in \B$ be their irreducible elements respectively. If $A_1,~B_1\in {\cal M}^i_{\mu}$ for the same $i$, then we can find two open sets $\emptyset \neq o_a,~o_b\subset {\cal M}^i_{\mu}$, $o_a\cap o_b=\emptyset$, such that $A_1\in o_a$ and $B_1\in o_b$. Then by definition, $s_a(o_a)\cap s_b(o_b)=\emptyset$ and $\A\in s_a$, $\B\in s_b$.

Finally, assume $A_1\in {\cal M}^i_{\mu}$, $B_1\in {\cal M}^j_{\mu}$ and $i\neq j$. Let $t=i\vee j$. Then
$$
A_{t/i}=A_1\otimes I_{t/i}\in {\cal M}^t_{\mu},\quad  B_{t/j}=B_1\otimes I_{t/j}\in {\cal M}^t_{\mu}.
$$
Since $A_{t/i}\neq B_{t/j}$, we can find $o_a,~o_b\subset {\cal M}^t_{\mu}$, $o_a\cap o_b=\emptyset$ and $A_{t/i}\in o_a$ and $B_{t/j}\in o_b$. That is,
$s_a(o_a)$ and $s_b(o_b)$ separate $\A$ and $\B$.
\hfill $\Box$


If we consider
\begin{align}\label{0.2.5}
{\cal M}:=\prod_{i=1}^{\infty}\prod_{j=1}^{\infty}{\cal M}_{i\times j}
\end{align}
as a product topological space, then ${\cal T}$ is the quotient topology of the standard product topology on the product space ${\cal M}$ defined by (\ref{0.2.5}). (We refer to \cite{wil70} for product topology.)

\subsection{Bundle Structure on ${\cal M}_{\mu}$}

\begin{dfn}[\cite{hus94}]\label{d0.3.1} A bundle is a triple $(E,p,B)$, where $E$ and $B$ are two topological spaces and $p:E\ra B$ is a continuous map. $E$ and $B$ are  called the total space and base space respectively. For each $b\in B$, $p^{-1}(b)$ is called the fiber of the bundle at $b\in B$.
\end{dfn}

 Observing the two topologies ${\cal T}_{{\cal M}}$ and ${\cal T}$ constructed in previous subsection, the following result is obvious:

\begin{prp}\label{p0.3.2} $\left({\cal M}_{\mu},\PR,\Sigma_{\mu}\right)$, is a bundle, where $\PR$ is the natural projection, i.e.,
$$
\PR(A)=\A.
$$
\end{prp}

\begin{figure}[!htbp]
\centering
\includegraphics[width=5cm]{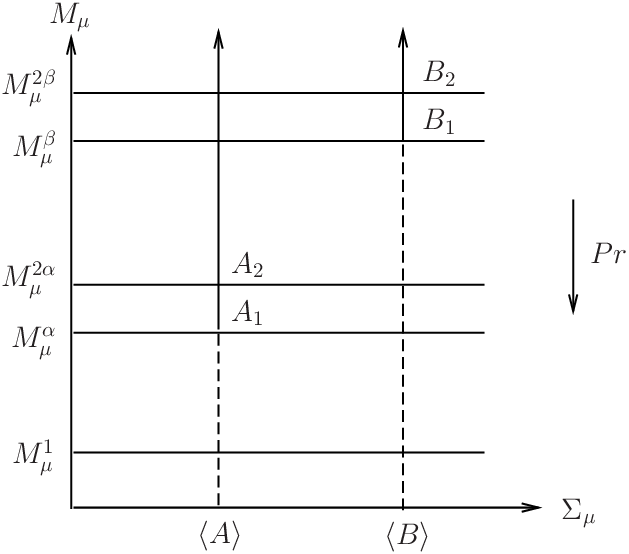}
\caption{Fiber Bundle Structure}\label{Fig.0.2}
\end{figure}

\begin{rem}
\label{r0.3.3}
\begin{enumerate}
\item Of course, $\left({\cal M},\PR,\Sigma_{{\cal M}}\right)$ is also a bundle. But it is of less interest because it is a discrete union of  $\left({\cal M}_{\mu},\PR,\Sigma_{\mu}\right)$, $\mu\in \Q_+$.
\item Consider an equivalent class $\A=\{A_1,~A_2,\cdots\}\in \Sigma_{\mu}$, where $A_1$ is irreducible. Then the fiber over $\A$ is a
discrete set:
$$
\PR^{-1}(\A)=\{A_1,A_2,A_3,\cdots\}.
$$
Hence this fiber bundle is named discrete bundle.

\item Fig.~\ref{Fig.0.2}
illustrates the fiber bundle structure of $\left({\cal M}_{\mu},\PR,\Sigma_{\mu}\right)$. Here $A_1\in \A$ and $B_1\in \B$ are irreducible $A_1\in {\cal M}_{\mu}^{\a}$ and $B_1\in {\cal M}_{\mu}^{\b}$. Their fibers are depicted in Fig.~\ref{Fig.0.2}.
\end{enumerate}
\end{rem}

We can define a set of cross sections \cite{hus94} $c_i:\Sigma_{\mu}\ra {\cal M}_{\mu}$ as:
\begin{align}\label{0.3.1}
c_i(\A):=A_i,\quad i=1,2,\cdots.
\end{align}
It is clear that $\PR\circ c_i=1_{\Sigma_{\mu}}$, where $1_{\Sigma_{\mu}}$ is the identity mapping on $\Sigma_{\mu}$.

Next, we consider some truncated sub-bundles of $\left({\cal M}_{\mu},\PR,\Sigma_{\mu}\right)$.

\begin{dfn}\label{d0.3.3} Assume $k\in \N$.
\begin{enumerate}
\item Set
\begin{align}
{\cal M}^{[\cdot,k]}_{\mu}:=\left\{M_{m\times n}\;|\; M_{m\times n}\in {\cal M}_{\mu}, ~\mbox{and}~m|k\mu_y\right\}.
\end{align}
Then ${\cal M}^{[\cdot,k]}_{\mu}$ is called the $k$-upper bounded subspace of  ${\cal M}_{\mu}$.
\item
$\Sigma^{[\cdot,k]}_{\mu}:={\cal M}^{[\cdot,k]}_{\mu}/\sim$ is called the $k$-upper bounded  subspace of  $\Sigma_{\mu}$.
\end{enumerate}
\end{dfn}

The natural projection $\PR:~{\cal M}^{[\cdot,k]}_{\mu}\ra \Sigma^{[\cdot, k]}_{\mu}$ is obviously defined. Then we have the following bundle structure.

\begin{prp}\label{p0.3.4} $\left({\cal M}_{\mu}^{[\cdot,k]},\PR,\Sigma_{\mu}^{[\cdot,k]}\right)$ is a sub-bundle of $\left({\cal M}_{\mu},\PR,\Sigma_{\mu}\right)$. Precisely speaking, the following graph (\ref{0.3.2}) is commutative, where $\pi$ and $\tilde{\pi}$ are including mappings. (\ref{0.3.2}) is also called the bundle morphism.
\begin{align}\label{0.3.2}
\begin{array}{ccc}
{\cal M}^{[\cdot,k]}_{\mu}&\xrightarrow{~~~\pi~~~}&{\cal M}_{\mu}\\
\PR\downarrow&~&\PR\downarrow\\
\Sigma_{\mu}^{[\cdot,k]}&\xrightarrow{~~~\pi'~~~}&\Sigma_{\mu}
\end{array}
\end{align}
\end{prp}

\begin{rem}\label{r0.3.5}
\begin{enumerate}
\item In a truncated sub-bundle there is a maximum cross section ${\cal M}_{kp\times kq}$ and a minimum cross section (i.e., root leaf) ${\cal M}_{p\times q}$, where $p,~q$ are co-prime and $p/q=\mu$.
\item Let ${\cal M}_{rp\times rq}$, $r=i_1,\cdots,i_t$ be a finite set of cross sections of $\left({\cal M}_{\mu},\PR,\Sigma_{\mu}\right)$. Set $k=i_1\vee i_2\vee \cdots\vee i_t$, then there exists a smallest truncation $\Sigma_{\mu}^{[\cdot,k]}$, which contains
    ${\cal M}_{rp\times rq}$, $r=i_1,\cdots,i_t$ as its cross sections.
\end{enumerate}
\end{rem}

\begin{dfn}\label{d0.3.6}
\begin{enumerate}
\item Define
$$
{\cal M}^{[k,\cdot]}_{\mu}:=\left\{M\in {\cal M}_{\mu}^s\;\big|\; k|s\right\},
$$
which is called the  $k$-lower bounded subspace of  ${\cal M}_{\mu}$.
\item Define the quotient space
$$
\Sigma^{[k,\cdot]}_{\mu}:={\cal M}^{[k,\cdot]}_{\mu}/\sim,
$$
which is called the  $k$-lower bounded subspace of  $\Sigma_{\mu}$.
\item Assume $\a|\b$. Define
$$
{\cal M}^{[\a,\b]}_{\mu}:={\cal M}^{[\a,\cdot]}_{\mu}\bigcap {\cal M}_{\mu}^{[\cdot,\b]},
$$
which is called the  $[\a,\b]$- bounded subspace of  ${\cal M}_{\mu}$.
\item Define the quotient space
$$
\Sigma^{[\a,\b]}_{\mu}:={\cal M}^{[\a,\b]}_{\mu}/\sim,
$$
which is called the  $[\a,\b]$-bounded subspace of  $\Sigma_{\mu}$.
\end{enumerate}
\end{dfn}

\begin{rem}\label{r0.3.7} Proposition \ref{p0.3.4} is also true for the two other truncated sub-bundles: $\left({\cal M}_{\mu}^{[k,\cdot]},\PR,\Sigma_{\mu}^{[k,\cdot]}\right)$ and $\left({\cal M}_{\mu}^{[\a,\b]},\PR,\Sigma_{\mu}^{[\a,\b]}\right)$ respectively. Precisely speaking, in (\ref{0.3.2}) if both ${\cal M}^{[\cdot,k]}_{\mu}$ and $\Sigma_{\mu}^{[\cdot,k]}$ are replaced by ${\cal M}^{[k, \cdot]}_{\mu}$ and $\Sigma_{\mu}^{[k,\cdot]}$ respectively, or by ${\cal M}^{[\a,\b]}_{\mu}$ and $\Sigma_{\mu}^{[\a,\b]}$ respectively, (\ref{0.3.2}) remains commutative.
\end{rem}

\subsection{Coordinate Frame on $\Sigma_{\mu}$}

It is obvious that $\left(\Sigma_{\mu}, \lplus\right)$ is an infinite dimensional vector space. Since each $\A\in \Sigma_{\mu}$ has finite coordinate expression, we may try to avoid using a basis with infinite elements. To this end, we construct a set of ``consistent" coordinate frames as $\{O_1,O_2,\cdots\}$. Then $\A$ can be expressed by $A_1\in \Span\{O_i\}$, $A_2\in \Span\{O_{i+1}\}$, and so on. Moreover, $O_i\subset O_{i+1}$ is a subset (or $O_i$ is part of coordinate elements in $O_{i+1}$). Then it seems that $\A$ can always be expressed in $O_i$ no matter which representative is chosen. The purpose of this section is to build such a set of consistent coordinate sub-frames, which forms an overall coordinate frame.

Assume $A_{\a}\in {\cal M}_{\mu}^{\a}$, $A_{\b}\in {\cal M}_{\mu}^{\b}$, $A_{\a}\sim A_{\b}$, and $\a|\b$, then $A_{\a}\otimes I_k=A_{\b}$, where $k=\b/\a$.
Recall that the order determined by this space-subspace relation is denoted as
\begin{align}\label{0.4.1}
{\cal M}_{\mu}^{\a}\sqsubset {\cal M}_{\mu}^{\b}.
\end{align}

 One sees easily that we can define an embedding mapping $\bd_k:{\cal M}_{\mu}^{\a}\ra {\cal M}_{\mu}^{\b}$ as
\begin{align}\label{0.4.2}
\bd_k(A):=A\otimes I_k.
\end{align}
In this way, ${\cal M}_{\mu}^{\a}$ can be considered as a subspace of ${\cal M}_{\mu}^{\b}$.

  In the following we construct a proper coordinate frame on ${\cal M}_{\mu}^{\b}$, which makes ${\cal M}_{\mu}^{\a}$ its coordinate subspace, that is,  ${\cal M}_{\mu}^{\a}$ is generated by part of coordinate variables of ${\cal M}_{\mu}^{\b}$. To this end,  we build a set of orthonormal basis on ${\cal M}_{\mu}^{\b}$ as follows:

Assume $p=\mu_y$ and $q=\mu_x$. Splitting $C\in {\cal M}_{\mu}^{\b}$ into $\a p\times \a q$ blocks, where each block is of dimension $k\times k$, yields
$$
C=\begin{bmatrix}
C^{1,1}&C^{1,2}&\cdots&C^{1,\a q}\\
\vdots&~&~&~\\
C^{\a p,1}&C^{\a p,2}&\cdots&C^{\a p,\a q}
\end{bmatrix}.
$$
Then for each $C^{I,J}\in {\cal M}_{k\times k}$ we construct a basis, which consists of three classes:
\begin{itemize}
\item Class 1:
\begin{align}\label{0.4.3}
\D^{I,J}_{i,j}=(b_{u,v})\in {\cal M}_{k\times k},\quad i\neq j,
\end{align}
where
$$
b_{u,v}=
\begin{cases}
1,\quad u=i,~v=j\\
0,\quad \mbox{otherwise}.
\end{cases}
$$
That is, for $(I,J)$-th block, at each given non-diagonal position $(i,j)$, set it to be $1$, and all other entries to be $0$.

\item Class 2:
\begin{align}\label{0.4.4}
D^{I,J}:=\frac{1}{\sqrt{k}}I_k^{I,J}.
\end{align}
That is, at each $(I,J)$-th block, set $D^{I,J}=\frac{1}{\sqrt{k}}I_k$ as a basis element.
\item Class 3:
\begin{align}\label{0.4.5}
\begin{array}{ccl}
E^{I,J}_t&=&\frac{1}{\sqrt{t(t-1)}}\diag\left(\underbrace{1,\cdots,1}_{t-1},-(t-1),\underbrace{0,\cdots,0}_{k-t}\right),\\
~&~&\quad t=2,\cdots,k.
\end{array}
\end{align}
That is, set $E^{I,J}_t$ as the rest of basis elements of the diagonal subspace of $(I,J)$-th block, which are orthogonal to $D^{I,J}$.
\end{itemize}

Let $A,B\in {\cal M}_{m\times n}$. Recall that the Frobenius inner product is defined as \cite{hor85}
\begin{align}\label{0.4.6}
(A|~B)_F:=\dsum_{i=1}^m\dsum_{j=1}^n a_{i,j}b_{i,j}.
\end{align}
Correspondingly, the Frobenius norm is defined as
\begin{align}\label{0.4.601}
\|A\|_F:=\sqrt{(A|~A)_F}.
\end{align}
If $(A|~B)_F=0$, then $A$ is said to be orthogonal with $B$.

Using Frobenius inner product, it is easy to verify the following result.

\begin{prp}\label{p0.4.1}
\begin{enumerate}
\item Set
\begin{align}\label{0.4.7}
\begin{array}{ccl}
B^{I,J}&:=&\left\{\D^{I,J}_{i,j},\;1\leq i\neq j\leq k;\;D^{I,J};\right.\\
~&~&\left. E^{I,J}_t,\;t=2,\cdots,k\right\},
\end{array}
\end{align}
then $B^{I,J}$ is an orthonormal basis for $(I,J)$-th block.
\item Set
$$
B:=\left\{B^{I,J}\;|\; I=1,2,\cdots,\a p;\;J=1,2,\cdots,\a q\right\},
$$
then $B$ is an orthonormal basis for ${\cal M}_{\mu}^{\b}$.
\item Set
$$
D:=\left\{D^{I,J}\;|\; I=1,2,\cdots,\a p;\;J=1,2,\cdots,\a q\right\},
$$
then $D$ is an orthonormal basis for subspace ${\cal M}_{\mu}^{\a}\subset {\cal M}_{\mu}^{\b}$.
\end{enumerate}
\end{prp}

\begin{exa}\label{e0.4.2}
Consider ${\cal M}_{1/2}^2 \sqsubset {\cal M}_{1/2}^4$. For any $A\in {\cal M}_{1/2}^4$, we split $A$ as
$$
A=\begin{bmatrix}
A^{1,1}&A^{1,2}&A^{1,3}&A^{1,4}\\
A^{2,1}&A^{2,2}&A^{2,3}&A^{2,4}
\end{bmatrix}.
$$
Then we build the orthonormal basis block-wise as
$$
\begin{array}{ccl}
B^{I,J}&:=&\left\{\D^{I,J}_{1,2}=\begin{bmatrix}0&1\\0&0\end{bmatrix},
\D^{I,J}_{2,1}=\begin{bmatrix}0&0\\1&0\end{bmatrix},\right.\\
~&~&\left.D^{I,J}=\frac{1}{\sqrt{2}}\begin{bmatrix}1&0\\0&1\end{bmatrix},
E^{I,J}_{2}=\frac{1}{\sqrt{2}}\begin{bmatrix}1&0\\0&-1\end{bmatrix}\right\}\\
\end{array}
$$
The orthonormal basis as proposed in Proposition \ref{p0.4.1} is
$$
B=\left\{B^{I,J}\;|\;I=1,2;\;J=1,2,3,4\right\}.
$$

Assume $A\in {\cal M}_{\mu}^{\b}$ and $C\in {\cal M}_{\mu}^{\a}$ and $\a k=\b$. Using (\ref{0.4.7}),  matrix $A$ can be expressed as
$$
A=\dsum_{I}\dsum_{J}\left(\dsum_{i\neq j}c^{I,J}_{i,j}\D^{I,J}_{i,j}+d_{I,J}D^{I,J}+\dsum_{t=2}^ke^{I,J}_tE^{I,J}_t\right).
$$

When ${\cal M}_{\mu}^{\a}$ is merged into ${\cal M}_{\mu}^{\b}$ as a subspace, matrix $C$ can be expressed as
$$
C=\dsum_{I}\dsum_{J}d'_{I,J}D^{I,J}.
$$
\end{exa}

\begin{dfn}\label{d0.4.3} Assume $C=(c_{i,j})\in {\cal M}_{\mu}^{\a}$ and $\b=k\a$. The embedding mapping $\bd_k: C\mapsto C\otimes I_k$ is defined as
\begin{align}\label{0.4.8}
\begin{array}{l}
\bd_k(C)_{I,J}=\sqrt{k}c_{I,J}D^{I,J},\\
~~~~~ I=1,\cdots,\a \mu_y;~J=1,\cdots,\a \mu_x.
\end{array}
\end{align}
\end{dfn}

To be consistent with this, we define the projection as follows:

\begin{dfn}\label{d0.4.4} Assume $A=(A^{I,J})\in {\cal M}_{\mu}^{\b}$, where $\b=k\a$ and $A^{I,J}=\left(A^{I,J}_{i,j}\right)\in {\cal M}_{k\times k}$. The projection
$\pr_k:{\cal M}_{\mu}^{\b}\ra {\cal M}_{\mu}^{\a}$ is defined as
\begin{align}\label{0.4.9}
\pr_k(A)=C=\left(c_{I,J}\right)\in {\cal M}_{\mu}^{\a},
\end{align}
where
$$
c_{I,J}=\frac{1}{k}\dsum_{d=1}^kA^{I,J}_{d,d}.
$$
\end{dfn}

According to the above construction, it is easy to verify the following:

\begin{prp}\label{p0.4.5}
The composed mapping $\pr_k\circ \bd_k$ is an identity mapping. Precisely,
\begin{align}\label{0.4.901}
\pr_k\circ \bd_k(C)=C,\quad \forall C\in {\cal M}^{\a}_{\mu}.
\end{align}
\end{prp}

\subsection{Inner Product on $\Sigma_{\mu}$}

Let $A=\left(a_{i,j}\right),~B=\left(b_{i,j}\right)\in {\cal M}_{m\times n}$. It is well known that the Frobenius inner product of $A$ and $B$ is defined by (\ref{0.4.6}), and
the Frobenius norm is defined by (\ref{0.4.601}).

The following lemma comes from a straightforward computation.

\begin{lem}\label{l0.5.1} Let $A,~B\in {\cal M}_{m\times n}$. Then
\begin{align}\label{0.5.1}
\left(A\otimes I_k\;|\;B\otimes I_k\right)_F=k(A\;|\;B)_F.
\end{align}
\end{lem}

\begin{dfn}\label{d0.5.2} Let $A,~B\in {\cal M}_{\mu}$, where $A\in {\cal M}_{\mu}^{\a}$ and  $B\in {\cal M}_{\mu}^{\b}$. Then the weighted inner product of $A,~B$ is defined as
\begin{align}\label{0.5.2}
(A\;|\;B)_W:=\frac{1}{t}\left(A\otimes I_{t/\a}\;|\;~B\otimes I_{t/\b}\right)_F,
\end{align}
where $t=\a\vee \b$ is the least common multiple of $\a$ and $\b$.
\end{dfn}

Using Lemma \ref{l0.5.1} and Definition \ref{d0.5.2}, we have the following property.
\begin{prp}\label{p0.5.201} Let $A,~B\in {\cal M}_{\mu}$, if $A$ and $B$ are orthogonal, i.e., $(A\;|\;~B)_F=0$, then
$A\otimes I_{\xi}$ and $B\otimes I_{\xi}$ are also orthogonal.
\end{prp}

Now we are ready to define the inner product on $\Sigma_{\mu}$.

\begin{dfn}\label{d0.5.3} Let $\A,~\B\in \Sigma_{\mu}$. Their inner product is defined as
\begin{align}\label{0.5.3}
(\A\;|\;~\B):=(A\;|\;~B)_W.
\end{align}
\end{dfn}

The following proposition shows that (\ref{0.5.3}) is well defined.

\begin{prp}\label{p0.5.4} Definition \ref{d0.5.3} is well defined. That is, (\ref{0.5.3}) is independent of the choice of the representatives $A$ and $B$.
\end{prp}

\noindent\textit{Proof}. Assume $A_1\in \A$ and $B_1\in \B$ are irreducible. Then it is enough to prove that
\begin{align}\label{0.5.4}
(A\;|\;~B)_W=(A_1\;|\;~B_1)_W,\quad A\in \A,~B\in \B.
\end{align}
Assume $A_1\in {\cal M}_{\mu}^{\a}$ and $B_1\in {\cal M}_{\mu}^{\b}$.   Let
$$
\begin{array}{l}
A=A_1\otimes I_{\xi}\in {\cal M}_{\mu}^{\a \xi},\\
B=B_1\otimes I_{\eta}\in {\cal M}_{\mu}^{\b \eta}.\\
\end{array}
$$
Denote by $t=\a\vee \b$, $s=\a\xi\vee \b\eta$, and $s=t\ell$.
Using (\ref{0.5.2}), we have
$$
\begin{array}{ccl}
(A\;|\;~B)_W&=&\frac{1}{s}\left(A\otimes I_{\frac{s}{\a\xi}}\;|\; B\otimes I_{\frac{s}{\b\eta}}\right)_F\\
~&=&\frac{1}{s}\left(A_1\otimes I_{\frac{s}{\a}}\;|\; B_1\otimes I_{\frac{s}{\b}}\right)_F\\
~&=&\frac{1}{t\ell}\left(A_1\otimes I_{\frac{t}{\a}}\otimes I_{\ell}\;|\; B_1\otimes I_{\frac{t}{\b}}\otimes I_{\ell}\right)_F\\
~&=&\frac{1}{t}\left(A_1\otimes I_{\frac{t}{\a}}\;|\; B_1\otimes I_{\frac{t}{\b}}\right)_F\\
~&=&(A_1\;|\;~B_1)_W.
\end{array}
$$
\hfill $\Box$

\begin{dfn}[\cite{tay80}]\label{d0.5.5} A real or complex vector space $X$ is an inner-product space, if there is a mapping $X\times X\ra \R ~(\mbox{or}~\C)$, denoted by $(x|~y)$, satisfying
\begin{enumerate}
\item
$$
(x+y\;|\;z)=(x\;|\;z)+(y\;|\;z),\quad x,y,z\in X.
$$
\item
$$
(x\;|\;y)=\overline{(y\;|\;x)},
$$
where the bar stands for complex conjugate.
\item
$$
(ax\;|\;y)=a(x\;|\;y),\quad a\in \R ~(\mbox{or}~\C).
$$
\item
$$
(x\;|\;x)\geq 0, ~\mbox{and}~ (x\;|\;x)\neq 0 ~\mbox{if}~ x\neq 0.
$$
\end{enumerate}
\end{dfn}

By definition it is easy to verify the following result.

\begin{thm}\label{t0.5.6} The vector space $(\Sigma_{\mu},\lplus)$ with the inner product defined by (\ref{0.5.3}) is an inner product space.
\end{thm}

Then the norm of $\A\in \Sigma_{\mu}$ is defined naturally as:
\begin{align}\label{0.5.5}
\|\A\|:=\sqrt{(\A\;|\;\A)}.
\end{align}

The following is some standard results for inner product space.

\begin{thm}\label{t0.5.7} Assume $\A,\B\in \Sigma_{\mu}$. Then we have the following
\begin{enumerate}
\item (Schwarz Inequality)
\begin{align}\label{0.5.6}
|(\A\;|\;\B)|\leq \|\A\|\|\B\|;
\end{align}
\item (Triangular Inequality)
\begin{align}\label{0.5.7}
\|\A\lplus \B\| \leq \|\A\|+\|\B\|;
\end{align}
\item (Parallelogram Law)
\begin{align}\label{0.5.8}
\begin{array}{l}
\|\A\lplus \B\|^2+\|\A\lminus \B\|^2 \\
= 2\|\A\|^2+2\|\B\|^2.
\end{array}
\end{align}
\end{enumerate}
\end{thm}

Note that the above properties show that $\Sigma_{\mu}$ is a normed space.

Finally, we present the generalized Pythagorean theorem:

\begin{thm}\label{t0.5.8} Let $\A_i\in \Sigma_{\mu}$, $i=1,2,\cdots,n$ be an orthogonal set. Then
\begin{align}\label{0.5.9}
\begin{array}{l}
\|\A_1\lplus \A_2\lplus \cdots\lplus \A_n\|^2 \\
= \|\A_1\|^2+\|\A_2\|^2+\cdots+\|A_n\|^2.
\end{array}
\end{align}
\end{thm}

A natural question is: ``Is $\Sigma_{\mu}$ a Hilbert space?" Unfortunately, this is not true. This fact is shown in the following counter-example.

\begin{exa}\label{e0.5.9} Define a sequence of elements, denoted as
 $\left\{\A_k\;\big|\;k=1,2,\cdots\right\}$,  as follows:
$A_1\in {\cal M}_{\mu}^1$ is arbitrary. Define $A_k$ inductively as
$$
A_{k+1}=A_k\otimes I_2+E_{k+1}\in {\cal M}_{\mu}^{2^{k}},\quad k=1,2,\cdots,
$$
where $E_{s}=\left(e^{s}_{i,j}\right)\in {\cal M}_{\mu}^{2^{s-1}}$ ($s\geq 2$) is defined as
$$
e^s_{i,j}=\begin{cases}
\frac{1}{2^s},& i=1,\;j=2,\\
0,            & \mbox{Otherwise}.
\end{cases}
$$

First, we claim that $\left\{\A_k:=\left<A_k\right>\;\big|\;k=1,2,\cdots\right\}$ is a Cauchy sequence. Let $n>m$. Then
\begin{align}\label{0.5.10}
\begin{array}{l}
\left\|\A_m\lminus \A_n\right\|\\
\leq \left\|\A_m\lminus\A_{m+1}\right\|+\cdots+\left\|\A_{n-1}\lminus\A_{n}\right\|\\
\leq \frac{1}{2^{m+1}}+\cdots+\frac{1}{2^{n}}\leq \frac{1}{2^{m}}.
\end{array}
\end{align}

Then we prove by contradiction that it does not converge to any element. Assume it converges to $\left<A_0\right>$, it is enough to consider the following three cases:

\begin{itemize}
\item Case 1, assume $A_0\in {\cal M}_{\mu}^{2^s}$ and $A_0=A_{s+1}$.
Then
$$
\left\|\left<A_0\right> \lminus \left<A_{s+2}\right>\right\|
=\left\|\left<A_{s+1}\right> \lminus \left<A_{s+2}\right>\right\|
=\frac{1}{2^{s+2}}.
$$
Similar to (\ref{0.5.10}) we can prove that
$$
\left\|\left<A_0\right> \lminus \left<A_{t}\right>\right\|
>\frac{1}{2^{s+2}},\quad t>s+2.
$$
Hence the sequence can not converge to $\left< A_0 \right>$.

\item Case 2, assume $A_0\in {\cal M}_{\mu}^{2^s}$ and $A_0\neq A_{s+1}$.
Note that $A_0 \lminus A_{s+1}$ is orthogonal to $E_{s+2}$, then it is clear that
$$
\left\|\left<A_0 \lminus A_{s+2}\right>\right\|>\left\|\left<A_0 \lminus A_{s+1}\right>\right\|.
$$
Note that by construction we have that as long as $t\geq s+2$ then
$A_{t}-A_{s+1}$ and $A_0-A_{s+1}$ are orthogonal. Using generalized Pythagorean theory, we have
$$
\begin{array}{l}
\left\|\left<A_0 \lminus A_{t}\right>\right\|\\
=\sqrt{\left[\left<A_0 \lminus A_{s+1}\right>\lplus \left<A_{s+1} \lminus A_{t}\right>\right]^2}\\
=\sqrt{\left[\left<A_0 \lminus A_{s+1}\right>\right]^2+\left[\left<A_{s+1} \lminus A_{t}\right>\right]^2}\\
>\|\left<A_0 \lminus A_{s+1}\right>\|>0,\quad t\geq s+2.
\end{array}
$$
Hence the sequence can not converge to $\left<A_0\right>$.

\item Case 3, $A_0\in {\cal M}_{\mu}^{2^s\xi}$, where ~$\xi>1$ is odd. Corresponding to Case 1, we assume
$A_0=A_{s+1}\otimes I_{\xi}$. Then we have
$$
\begin{array}{l}
\left\|\left<A_0 \lminus A_{s+2}\right>\right\|\\
=\left\|\left<A_{s+1}\otimes I_{\xi} \lminus (A_{s+1}\otimes I_2+E_{s+2})\otimes I_{\xi}\right>\right\|\\
=\left\|\left<E_{s+2})\otimes I_{\xi}\right>\right\|\\
=\frac{1}{2^{s+2}}.
\end{array}
$$
and
$$
\left\|\left<A_0 \lminus A_{t}\right>\right\|>\frac{1}{2^{s+2}},\quad t>s+2.
$$
So the sequence can not converge to $\left< A_0 \right>$.

Corresponding to Case 2, assume $A_0\neq A_{s+1}\otimes I_{\xi}$. Using Proposition \ref{p0.5.201}, a similar argument shows that the sequence cannot converge to $\left< A_0 \right>$ too.

\end{itemize}
\end{exa}

\subsection{$\Sigma_{\mu}$ as a Matric Space}

Using the norm defined in previous section one sees easily that $\Sigma_{\mu}$ is a matric space:

\begin{thm}\label{t0.6.1} $\Sigma_{\mu}$ with distance
\begin{align}\label{0.6.1}
d(\A,\B):=\|\A\lminus \B\|, \quad \A,~\B\in \Sigma_{\mu}
\end{align}
is a matric space.
\end{thm}

\begin{thm}\label{t0.6.2} Consider $\Sigma_{\mu}$. The topology deduced by the distance $d$, denoted by ${\cal T}_d$ is exactly the same as the topology ${\cal T}$ defined in Definition \ref{d0.2.2}.
\end{thm}

\noindent\textit{Proof}. Assume $U\in {\cal T}_d$ and $p\in U$. Then there exists a ball $B_{\epsilon}(p)$ such that
$B_{\epsilon}(p)\subset U$, where $\epsilon>0$. Assume
$
p=\left<A_0\right>$ and $A_0\in {\cal M}_{\mu}^s={\cal M}_{s\mu_y\times s\mu_x}$.

Now we can construct a ball ${\cal B}_{\d}(A_0)\subset {\cal M}_{\mu}^s$, where $\d>0$. Note that ${\cal B}_{\d}(A_0)$ is a sub-basis element of ${\cal T}$, and hence is an open set in $\left(\Sigma_{\mu}, {\cal T}\right)$. By continuity, as $\d>0$ small enough, $q\in {\cal B}_{\d}(A_0)$ implies $d(q,A_0)<\epsilon$. That is,
$$
{\cal B}_{\d}\left(\left<A_0\right>\right)\subset B_{\epsilon}(p)\subset U,
$$
which means $U\in {\cal T}$. Hence, ${\cal T}_d\subset {\cal T}$.

Conversely, assume $q\in V\in {\cal T}$. Then there exists a basic open set $s_1 \cap \cdots \cap s_r\in {\cal T}$ such that $q\in s_1 \cap \cdots \cap s_r\in {\cal T}$. Express
$$
q=\left\{A_1,~A_2,~\cdots,~A_r\cdots\right\},
$$
where $A_i\in s_i\subset {\cal M}_{\mu}^{\xi_i}$, $i=1,\cdots,r$. For each $A_i$, we can find ${\cal B}_{\epsilon_i}^i(\left<A_i\right>)\subset s_i$, where $\epsilon_i>0$, $i=1,\cdots,r$. Choosing $\d_i>0$ small enough such that
$$
B_{\d_i}^i\left(\left<A_i\right>\right)\subset {\cal B}_{\epsilon_i}^i(\left<A_i\right>)
\subset s_i,\quad i=1,\cdots,r.
$$
Then we have
$$
q\in \bigcap_{i=1}^rB_{\d_i}^i\left(\left<A_i\right>\right)\in {\cal T}_d.
$$
That is, $V\in {\cal T}_d$. Hence, ${\cal T}\subset {\cal T}_d$.

We conclude that  ${\cal T}_d = {\cal T}$.
\hfill $\Box$

\begin{dfn}[\cite{dug66}]\label{d0.6.3}
\begin{enumerate}
\item A topological space is regular (or $T_3$) if for each closed set $X$ and $x\not\in X$ there exist open neighborhoods $U_x$ of $x$ and $U_X$ of $X$, such that $U_x\cap U_X=\emptyset$.
\item A topological space is normal (or $T_4$) if for each pair of closed sets $X$ and $Y$ there exist open neighborhoods $U_X$ of $X$ and $U_Y$ of $Y$, such that $U_X\cap U_Y=\emptyset$.
\end{enumerate}
\end{dfn}

Since a matric space is regular and normal, as a corollary of Theorem \ref{t0.6.2}, we have the following result.

\begin{cor}\label{c0.6.4} The topological space $\left(\Sigma_{\mu},{\cal T}\right)$, defined in Definition \ref{d0.2.2}, is both regular and normal.
\end{cor}

Note that
$$
T_4\Rightarrow T_3\Rightarrow T_2.
$$

Finally, we show some properties of $\Sigma_{\mu}$.

\begin{prp}\label{p0.6.5}
$\Sigma_{\mu}$ is convex. Hence it is arcwise connected.
\end{prp}
\noindent\textit{Proof}. Assume $\A,~\B\in \Sigma_{\mu}$. Then it is clear that
$$
\lambda\A\lplus (1-\lambda)\B=\left<\lambda A\lplus (1-\lambda) B\right>\in \Sigma_{\mu},\quad \lambda\in [0,1].
$$
So $\Sigma_{\mu}$ is convex. Let $\lambda$ go from $1$ to $0$, we have a path connecting $\A$ and $\B$.
\hfill $\Box$

\begin{prp}\label{p0.6.6}
$\Sigma_{\mu}$ and $\Sigma_{1/\mu}$ are isometric spaces.
\end{prp}

\noindent\textit{Proof}. Consider the transpose:
$$
\A\mapsto \left<A^T\right>.
$$
Then it is obvious that
$$
d(\A,\B)=d\left(\left<A^T\right>, \left<B^T\right>\right).
$$
Hence the transpose is an isometry. Moreover, its inverse is itself.
\hfill $\Box$

\subsection{Subspaces of $\Sigma_{\mu}$}

Consider the $k$-upper bounded subspace $\Sigma_{\mu}^{[\cdot,k]}$. We have

\begin{prp}\label{p0.7.1} $\Sigma_{\mu}^{[\cdot,k]}$ is a Hilbert space.
\end{prp}

\noindent\textit{Proof}. Since $\Sigma_{\mu}^{[\cdot,k]}$ is a finite dimensional vector space and any finite dimensional inner product space is a Hilbert space \cite{die69}, the conclusion follows.
\hfill $\Box$

\begin{prp}[\cite{die69}]\label{p0.7.2} Let $E$ be an inner product space, $\{0\}\neq F\subset E$ be a Hilbert subspace.
\begin{enumerate}
\item For each $x\in E$ there exists a unique $y:=P_F(x)\in F$, called the projection of $x$ on $F$, such that
\begin{align}\label{0.7.1}
\|x-y\|=\min_{z\in F}\|x-z\|.
\end{align}
\item
\begin{align}\label{0.7.2}
F^{\perp}:=P_F^{-1}(0)
\end{align}
is the subspace orthogonal to $F$.
\item
\begin{align}\label{0.7.3}
E=F\oplus F^{\perp},
\end{align}
where $\oplus$ stands for orthogonal sum.
\end{enumerate}
\end{prp}

Using above proposition, we consider the projection: $P_F:\Sigma_{\mu}\ra \Sigma_{\mu}^{[\cdot,\a]}$. Let $\A\in \Sigma_{\mu}^{\b}$. Assume $\left<X\right>\in \Sigma_{\mu}^{\a}$, $t=\a\vee \b$. Then the norm of $\A\lminus \left<X\right>$ is:
\begin{align}\label{0.7.4}
\left\|\A\lminus \left<X\right>\right\|=\frac{1}{\sqrt{t}}\left\|A\otimes I_{t/\b}-X\otimes I_{t/\a}\right\|_F.
\end{align}
Set $p=\mu_y$ $q=\mu_x$, and $k:=t/\a$. We split $A$ as
$$
A\otimes I_{t/\b}=\begin{bmatrix}
A_{1,1}&A_{1,2}&\cdots&A_{1,q\a}\\
A_{2,1}&A_{2,2}&\cdots&A_{2,q\a}\\
\vdots&~&~&~\\
A_{p\a,1}&A_{p\a,2}&\cdots&A_{p\a,q\a}\\
\end{bmatrix},
$$
where $A_{i,j}\in {\cal M}_{k\times k}$, $i=1,\cdots,p\a;~j=1,\cdots,q\a$.
set
\begin{align}\label{0.7.5}
C:=\argmin_{x\in \Sigma_{\mu}^{\a}}\left\|A\otimes I_{t/\b}-X\otimes I_{t/\a}\right\|.
\end{align}
Then it is easy to calculate that
\begin{align}\label{0.7.6}
c_{i,j}=\frac{1}{k}\tr(A_{i,j}),\quad i=1,\cdots,p\a;~j=1,\cdots,q\a.
\end{align}
We conclude that
\begin{prp}\label{p0.7.3} Let $\A\in \Sigma_{\mu}$, $P_F: \Sigma_{\mu}\ra \Sigma_{\mu}^{[\cdot,\a]}$. Precisely speaking, $\A\in \Sigma_{\mu}^{\b}$. Using the above notations, the projection of $\A$ is
\begin{align}\label{0.7.7}
P_F(\A)=\left<C\right>,
\end{align}
where $C$ is defined in (\ref{0.7.6}).
\end{prp}

We give an example to depict this.

\begin{exa}\label{e0.7.4}

Given
$$
A=\begin{bmatrix}
1&2&-3&0&2&1\\
2&1&-2&-1&1&0\\
0&-1&-1&3&1&-2
\end{bmatrix}\in \Sigma_{0.5}^3.
$$
We consider the projection of $\A$ onto $\Sigma_{0.5}^{[\cdot,2]}$. Denote $t=2\vee 3=6$. Using formulas (\ref{0.7.6})--(\ref{0.7.7}), we have
$$
P_F(\A)=\left<\begin{bmatrix}
1&0&1/3&0\\
0&-1/3&0&-1\\
\end{bmatrix}\right>.
$$
Then we have
$$
\left<E\right>=\A \lminus P_F(\A),
$$
where
$$
\begin{array}{l}
E=\\
\left[
\begin{array}{cccccccccccc}
0&0&2&0&-3&0&-1/3&0&2&0&1&0\\
0&0&0&2&0&-3&0&-1/3&0&2&0&1\\
2&0&0&0&-2&0&-1&0&2/3&0&0&0\\
0&2&0&4/3&0&-2&0&-1&0&2&0&0\\
0&0&-1&0&-2/3&0&3&0&1&0&-1&0\\
0&0&0&-1&0&-2/3&0&3&0&1&0&-1\\
\end{array}\right].
\end{array}
$$

It is easy to verify that $\left< E \right>$ and $\A$ are mutually orthogonal.

\end{exa}

 We also have $\Sigma_{\mu}^{[k,\cdot]}$ and $\Sigma_{\mu}^{[\a,\b]}$ (where $\a|\b$) as metric subspaces of $\Sigma_{\mu}$.

 Finally, we would like to point out that since $\Sigma_{\mu}$ is an infinity dimensional vector space, it is possible that $\Sigma_{\mu}$ is isometric to its proper subspace. For instance, consider the following example.

 \begin{exa}\label{e0.7.5} Consider a mapping $\varphi:\Sigma_{\mu}\ra \Sigma_{\mu}^{[k,\cdot]}$ defined by
 $\A\mapsto \left<A\otimes I_k\right>$. It is clear that this mapping satisfies
 $$
 \|\A \lminus \B\|=\|\varphi(\A) \lminus \varphi(\B)\|,\quad \A,\B\in \Sigma_{\mu}.
 $$
 That is, ${\cal M}_{\mu}$ can be isometrically embedded into its proper subspace.
 \end{exa}


\section{Differential Structure on M-equivalence Space}

\subsection{Bundled Manifold}

Unlike the conventional manifolds, which have fixed dimensions, this section explores a new kinds of manifolds, called the bundled manifold. Intuitively speaking, it is a fiber bundle, which has fibers belonging to manifolds of different dimensions. To begin with, the following definition is proposed, which is a mimic to the definition of an $n$-dimensional manifold \cite{boo79}.
\begin{dfn}\label{d3.2.1} Let $\{M,~{\cal T}\}$ be a topological space.
\begin{enumerate}
\item An open set $U\neq \emptyset$ is said to be a simple coordinate chart, if there is an open set $\varTheta\subset \R^s$ and a homeomorphism $\phi:U\ra \varTheta\subset \R^s$. The integer $s$ is said to be the dimension of $U$.
\item An open set $U\neq \emptyset$ is said to be a bundled coordinate chart, if there exist finite simple coordinate charts $U_i$ with homeomorphisms $\phi_i:U_i\ra \varTheta_i\subset \R^{s_i}$, $i=1,\cdots,k$, $k<\infty$ is set $U$ depending, such that
    $$
    U=\bigcap_{i=1}^{k}U_i.
    $$
\item Let $U$, $V$ be two bundled coordinate charts. $U=\cap_{i=1}^{k_1}U_i$, $V=\cap_{j=1}^{k_2}V_j$.
$U$ and $V$ are said to be $C^r$ comparable if for any $U_i$ and $V_j$, as long as their dimensions are equal, they are $C^r$ comparable. (Where $r$ could be ${\infty}$, that is they are $C^{\infty}$ comparable;  or $\omega$, which means they are analytically comparable.)
\end{enumerate}
\end{dfn}

\begin{rem}\label{r3.2.2} In Definition \ref{d3.2.1} for a bundled coordinate chart $U=\bigcap_{i=1}^{k}U_i$, we can, without loss of generality, assume $\dim(U_i)$, $i=1,\cdots,k$ are distinct. Because simple coordinate charts of same dimension can be put together by set union $\cup$. Hereafter, this is assumed.
\end{rem}

\begin{figure}[!htbp]
\centering\includegraphics[width=7cm]{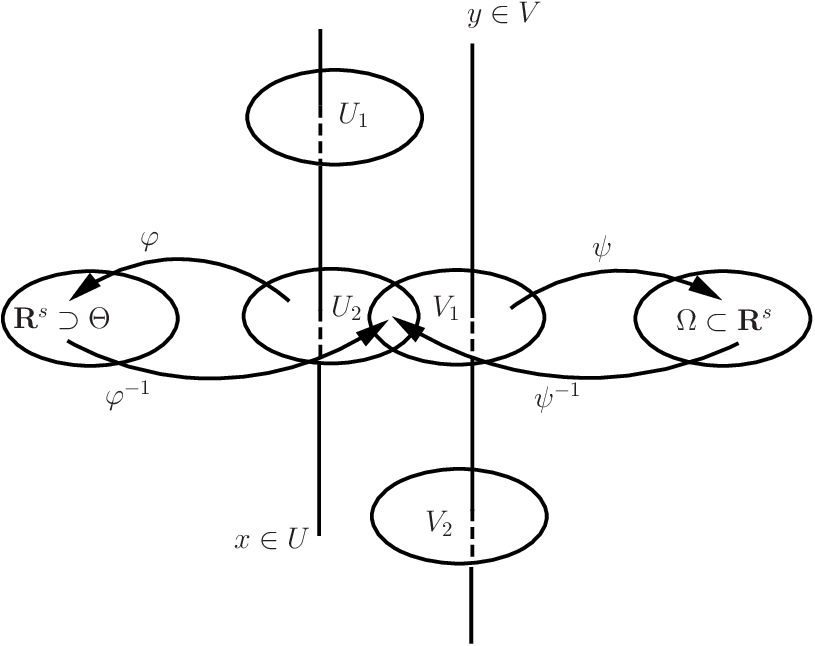}
\caption{Multiple Coordinate Charts}\label{Fig.3.1}
\end{figure}

To depict the bundle coordinate chart, we refer to Fig.~\ref{Fig.3.1}, where two bundled coordinate charts $U$ and $V$ are described. Note that $U=U_1\cap U_2$ and $V=V_1\cap V_2$, where $U_1,~U_2,~V_1,~V_2$ are simple coordinate charts. Also, $U$ and $V$ are coordinate neighborhoods of $x$ and $y$ respectively. Next, we assume that $\dim(U_2)=\dim(V_1)=s$, $U_2\cap V_1\neq \emptyset$, and  $U_1\cap V_2 = \emptyset$. $U$ and $V$ are $C^r$  comparable, if and only if the two mappings
$$
\begin{array}{l}
\psi\circ\phi^{-1}:\phi(U_2\cap V_1)\ra \psi(U_2\cap V_1)~\mbox{and}\\
\phi\circ\psi^{-1}:\psi(U_2\cap V_1)\ra \phi(U_2\cap V_1)
\end{array}
$$
are $C^r$.
As a convention, we assume $\dim(U_1)\neq \dim(U_2)$, $\dim(V_1)\neq \dim(V_2)$.

\begin{dfn}\label{d3.2.2} A topological space $M$ is a bundled $C^r$ (or $C^{\infty}$, or analytic, denoted by $C^{\omega}$) manifold, if the following conditions are satisfied.
\begin{enumerate}
\item $M$ is second countable and Hausdorff.
\item There exists an open cover of $M$, described as
$${\cal C}=\{U_{\lambda}\;|\; \lambda\in \Lambda\},$$
where each $U_{\lambda}$ is a  bundled coordinate chart. Moreover, any two bundled coordinate charts in ${\cal C}$ are   $C^r$  comparable.
\item If a bundled coordinate chart $V$ is comparable with $U_{\lambda}$, $\forall \lambda\in \Lambda$, then $V\in {\cal C}$.
\end{enumerate}
\end{dfn}

It is obvious that the  topological structure of $\Sigma_{\mu}$ with natural ($\R^{i^2\mu_y\mu_x}$) coordinates on each cross section (or leaf) meets the above requirements for a bundled manifold. Hence, we have the following result.

 \begin{thm}\label{t3.2.3} $\Sigma_{\mu}$ is a bundled analytic manifold.
\end{thm}

\noindent\textit{Proof.} Condition 1 has been proved in Theorem \ref{t0.2.4}. For condition $2$, set $p=\mu_y$, $q=\mu_x$ and
$$
O_k:={\cal M}_{\mu}^{k},\quad k=1,2,\cdots.
$$
Choosing any finite open subset $o_{i_s}\subset O_{i_s}$, $s=1,2,\cdots, t$, $t<\infty$, and constructing corresponding $s(o_{i_1}),\cdots,s(o_{i_t})$.
Set $U_I:=s(o_{i_1})\cap\cdots\cap s(o_{i_t})$, where $I=\{i_1,i_2,\cdots,i_t\}$.
Define
$$
W:=\left\{U_I\;|; I ~\mbox{is a finite subset of}~ \N \right\}.
$$
Then  $W$ is an open cover of $M$. Identity mappings from $s(o_{i})\ra {\cal M}_{\mu}^i\simeq \R^{ip\times iq}$ makes any two $U_I$ and $U_J$ being $C^{\omega}$ comparable.
As for condition 3, just add all bundled coordinate charts which are comparable with $W$ into $W$, the condition is satisfied.
\hfill $\Box$

Next, we consider the lattice-related coordinates on $\Sigma_{\mu}$.

Consider $\Sigma_{\mu}$ and assume $p=\mu_y$ and $q=\mu_x$. Then $\Sigma_{\mu}$ has leafs
 $$
 \Sigma_{\mu}=\left\{{\cal M}_{\mu}^i\;|\;i=1,2,\cdots\right\},
 $$
 where ${\cal M}_{\mu}^i={\cal M}_{ip\times iq}$.

 Consider an element $\left<x\right> \in \Sigma_{\mu}$, then there exists a unique irreducible $x_1\in \left<x\right>$ such that $\left<x\right>=\left\{x_j=x_1\otimes I_j\;|\;j=1,2,\cdots \right\}$. Now assume $x_1\in {\cal M}_{\mu}^s$. As defined above, ${\cal M}_{\mu}^s$ is the root leaf of $\left<x\right>$.

\begin{figure}[!htbp]
\centering
\includegraphics[width=8cm]{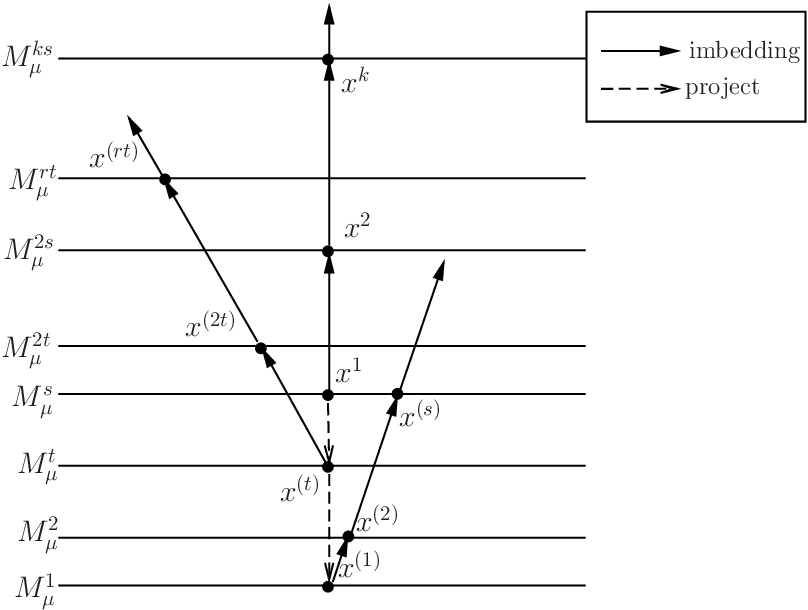}
\caption{Lattice-related Coordinates}\label{Fig.3.2}
\end{figure}

It is obvious that $\left<x\right>$ has different coordinate representations on different leafs. But because of the subspace lattice structure, they must be consistent. Fig.~\ref{Fig.3.2} shows the lattice-related subspaces. Any geometric objects defined on its root leaf must be consistent with its representations on all embedded spaces and projected spaces.

 As shown in Fig. 6 the following subspaces are related:
\begin{itemize}
\item Class 1 (Embedded Elements): Starting from $x^1\in {\cal M}_{\mu}^s$, we have
$$
x^1\prec x^2\prec x^3\prec \cdots.
$$
\item Class 2 (Projected Elements): Let $x^{(t)}\in {\cal M}_{\mu}^t$, where $t|s$. Then
$$
x^{(t)}\prec x^1.
$$
Particularly, $x^{(1)}\in {\cal M}_{\mu}^1$ satisfies
$$
x^{(1)}\prec x^1.
$$
\item Class 3 (Embedded Elements from Projected Elements): Starting from any $x^{(t)}$, $t|s$, we have
$$
x^{(t)}\prec x^{(2t)}\prec x^{(3t)}\prec \cdots.
$$
\end{itemize}

\begin{rem}\label{r3.2.4}
\begin{enumerate}

\item Classes 1--3 are the set of coordinates, which are related to a given irreducible element $x_1$.
\item Elements in Class 1 are particularly important. Say, we may firstly define a root element on ${\cal M}_{\mu}^s$, such as $A_1\in {\cal M}_{\mu}^s$. Then we use it to get an equivalent class, such as $\A$, and use the elements in this class to perform certain calculation, such as STP. All the elements in the equivalent class, such as $\A$, are of Class 1.
\item The elements in subspace and their equivalent classes are less important. Sometimes we may concern only the elements of Class 1, say for STP etc.
\item The subspace elements obtained by project mapping may not be ``uniformed" with the object obtained from the real subspaces of the original space. More discussion will be seen in the sequel.
\item Because of the above argument, sometimes we may consider only the equivalent classes which have their root elements defined on their root leaf. Therefore, the objects may only be defined on the multiple of the root leaf (root space).
\end{enumerate}
\end{rem}

\subsection{$C^{r}$ Functions on $\Sigma_{\mu}$}

\begin{dfn}\label{d3.3.1} Let $M$ be a bundled manifold, $f:~M\ra \R$ is called a $C^r$ function, if for each simple coordinate chart $U\subset {\cal M}_{\mu}^s$ $f|_U:=f_s$ is $C^r$.
The set of $C^r$ functions on $M$ is denoted by $C^r(M)$.
\end{dfn}

Assume $f\in C^r(\Sigma_{\mu})$, $A,~B\in {\cal M}_{\mu}$ and $A\sim B$. Then $f$ is well defined on $\Sigma_{\mu}$ means it is defined on different leafs consistently, and hence on leafs corresponding to $A$ and $B$ we  have $f(A)=f(B)$. To this end, the $f$ can be constructed as follows:

\begin{dfn}\label{d3.3.2} Assume $f$ is firstly defined on root leaf ${\cal M}_{\mu}^s$ as $f_s(x)$. Then we extend it to other leafs as:
\begin{itemize}
\item Step 1. Let $Q=\{t\in \N\;\big|\; t|s\}$. Then
\begin{align}\label{3.3.1}
f_t(y):=f_s(x=\bd_k(y)),
\end{align}
where $k=\frac{s}{t}\in \N$.
\item Step 2. Assume $gcd(\ell,s)=t$. If $\ell=t$, $f_{\ell}=f_t$ has already been defined in Step 1. So we assume $\ell=kt$. Then
\begin{align}\label{3.2}
f_{\ell}(y):=f_t(x=\pr_k(y)).
\end{align}

Note that in Step 2, $t=s$ is allowed.
\end{itemize}
\end{dfn}

Then it is easy to verify the following:

\begin{prp}\label{p3.3.3} The function $f$ defined in Definition \ref{d3.3.2} is consistent with the equivalence $\sim$. Hence it is well defined on $\Sigma_{\mu}$.
\end{prp}

\begin{exa}\label{e3.3.4} Consider $\Sigma_2$, and assume $f$ is defined on its root leaf ${\cal M}_{2}^2$ firstly as
$$
f_2\left(\begin{bmatrix}
a_{11}&a_{12}\\
a_{21}&a_{22}\\
a_{31}&a_{32}\\
a_{41}&a_{42}\\
\end{bmatrix}
\right):=a_{11}a_{22}-a_{32}a_{41}.
$$
Then we can determine the other expressions of $f$ as follows:
\begin{itemize}
\item Consider $f_1$:
$$
f_1\left(\begin{bmatrix}
a_{1}\\
a_{2}
\end{bmatrix}
\right)=f_2\left(\begin{bmatrix}
a_{1}\\
a_{2}
\end{bmatrix}\otimes I_2\right)=a_1^2.
$$
\item Consider $f_3$. Let $A=\left(a_{i,j}\right)\in {\cal M}_{2}^3$. Then
$$
\begin{array}{l}
f_3(A)=f_1(\pr_3(A))=f_1\left(\begin{bmatrix}
\frac{1}{3}(a_{11}+a_{22}+a_{33})\\
\frac{1}{3}(a_{41}+a_{52}+a_{63})
\end{bmatrix}
\right)\\
=\frac{1}{9}(a_{11}+a_{22}+a_{33})^2.
\end{array}
$$

Similarly, for any $n=2k-1$ we have
$$
f_n(A)=\frac{1}{n^2}\left(a_{11}+a_{22}+\cdots+a_{nn}\right)^2.
$$

\item Consider $f_4$. Let $A=\left(a_{i,j}\right)\in {\cal M}_{2}^4$. Then
$$
\begin{array}{l}
f_4(A)=f_2(\pr_2(A))\\
=f_2\left(\begin{bmatrix}
\frac{1}{2}(a_{11}+a_{22})& \frac{1}{2}(a_{13}+a_{24})\\
\frac{1}{2}(a_{31}+a_{42})& \frac{1}{2}(a_{33}+a_{44})\\
\frac{1}{2}(a_{51}+a_{62})& \frac{1}{2}(a_{53}+a_{64})\\
\frac{1}{2}(a_{71}+a_{82})& \frac{1}{2}(a_{73}+a_{84})\\
\end{bmatrix}
\right)\\
=\frac{1}{4}\left[(a_{11}+a_{22})(a_{33}+a_{44})-
(a_{53}+a_{64})(a_{71}+a_{82})\right]
\end{array}
$$

Similarly, for $n=2k$ ($k\geq 2$) we have
$$
\begin{array}{ccl}
f_n(A)&=\frac{1}{k^2}&\left[\left(a_{11}+a_{22}+\cdots+a_{kk}\right)\right.\\
~&~&\left(a_{k+1,k+1}+a_{k+2,k+2}+\cdots+a_{2k,2k}\right)\\
~&~&-\left(a_{2k+1,k+1}+a_{2k+2,k+2}+\cdots+a_{3k,2k}\right)\\
~&~&\left.\left(a_{3k+1,1}+a_{3k+2,2}+\cdots+a_{4k,k}\right)\right].
\end{array}
$$
\end{itemize}
\end{exa}

\begin{rem}\label{r3.3.5} For a smooth function $f$ defined firstly on ${\cal M}_{\mu}^s$, its extensions to both ${\cal M}_{\mu}^{[\cdot,s]}$ and ${\cal M}_{\mu}^{[s, \cdot]}$ are consistently defined. Hence, $f$ is completely well posed on $\Sigma_{\mu}$.
\end{rem}

\subsection{Generalized Inner Products}

We define the generalized Frobenius inner product as follows.

\begin{dfn}\label{d3.4.1} Given $A\in {\cal M}_{m\times n}$ and $B\in {\cal M}_{p\times q}$.
\begin{enumerate}
\item[Case 1] (Special Case): Assume $p=rm$ and $q=sn$.
Split $B$ into equal blocks as
$$
B=\begin{bmatrix}
B_{1,1}&B_{1,2}&\cdots&B_{1,s}\\
B_{2,1}&B_{2,2}&\cdots&B_{2,s}\\
\vdots&~&~&~\\
B_{r,1}&B_{r,2}&\cdots&B_{r,s}\\
\end{bmatrix},
$$
where $B_{i,j}\in {\cal M}_{m\times n}$, $i=1,\cdots,r$; $j=1,\cdots,s$.
Then the generalized Frobenius inner product of $A$ and $B$ is defined as
\begin{align}\label{3.4.1}
\begin{array}{l}
(A\;\big|\;
 B)_F:=\\
\begin{bmatrix}
\left(A|B_{1,1}\right)_F&\left(A|B_{1,2}\right)_F&\cdots&\left(A|B_{1,s}\right)_F\\
\left(A|B_{2,1}\right)_F&\left(A|B_{2,2}\right)_F&\cdots&\left(A|B_{2,s}\right)_F\\
\vdots&~&~&~\\
\left(A|B_{r,1}\right)_F&\left(A|B_{r,2}\right)_F&\cdots&\left(A|B_{r,s}\right)_F\\
\end{bmatrix}.
\end{array}
\end{align}
Note that here $(A|B_{i,j})_F$ is the standard Frobenius inner product defined in (\ref{0.4.6}).

\item[Case 2] (General Case): Assume $A\in {\cal M}_{m\times n}$ and $B\in {\cal M}_{p\times q}$ and let the great common divisor of $m,~p$ be $\a=m\wedge p$, and the great common divisor of $n,~q$ be $\b=n\wedge q$. Denote by $\xi=m/\a$ and $\eta=n/\b$, $r=p/\a$ and $s=q/\b$. Then we split $A$ into $\xi\times \eta$ blocks as
$$
A=\begin{bmatrix}
A_{1,1}&A_{1,2}&\cdots&A_{1,\eta}\\
A_{2,1}&A_{2,2}&\cdots&A_{2,\eta}\\
\vdots&~&~&~\\
A_{\xi,1}&A_{\xi,2}&\cdots&A_{\xi,\eta}\\
\end{bmatrix},
$$
where $A_{i,j}\in {\cal M}_{\a\times \b}$, $i=1,\cdots,\xi$; $j=1,\cdots,\eta$.
Then the generalized Frobenius inner product of $A$ and $B$ is defined as
\begin{align}\label{3.4.2}
\begin{array}{l}
(A\;\big|\;
 B)_F:=\\
\begin{bmatrix}
\left(A_{1,1}|B\right)_F&\left(A_{1,2}|B\right)_F&\cdots&\left(A_{1,\eta}|B\right)_F\\
\left(A_{2,1}|B\right)_F&\left(A_{2,2}|B\right)_F&\cdots&\left(A_{2,\eta}|B\right)_F\\
\vdots&~&~&~\\
\left(A_{\xi,1}|B\right)_F&\left(A_{\xi,2}|B\right)_F&\cdots&\left(A_{\xi,\eta}|B\right)_F\\
\end{bmatrix}.
\end{array}
\end{align}
Note that here $(A_{i,j}|B)_F$ is the (Case 1) generalized Frobenius inner product defined in (\ref{3.4.1}).
\end{enumerate}
\end{dfn}

\begin{exa}\label{e3.4.2}
Let
$$
A=\begin{bmatrix}
1&-1&1&0\\
1&2&0&1
\end{bmatrix}\in {\cal M}_{0.5}^2,
$$
and
$$
B=\begin{bmatrix}
1&0\\
-1&2\\
-1&0\\
1&-1
\end{bmatrix}\in {\cal M}_{3/2}^1.
$$
Note that $m=2$, $n=4$, $p=4$, $q=2$, $\a=\gcd(m,p)=2$, $\b=\gcd(n,q)=2$.
Then we split $A$ and $B$ as follows
$$
A=\begin{bmatrix}
A_{1,1}&A_{1,2}
\end{bmatrix},\quad
B=\begin{bmatrix}
B_{1,1}\\
B_{2,1}
\end{bmatrix},
$$
were $A_{i,j},~B_{k,\ell}\in {\cal M}_{2\times 2}$, $i=1,2;~j=1,2;~k=1,2;~\ell=1$.

Finally, we have
$$
(A\;|\; B)_F=
\begin{bmatrix}
(A_{1,1}|B_{1,1})&(A_{1,2}|B_{1,1})\\
(A_{1,1}|B_{2,1})&(A_{1,2}|B_{2,1})\\
\end{bmatrix}=\begin{bmatrix}
4&3\\
-2&-2
\end{bmatrix}.
$$
\end{exa}

\begin{dfn}\label{d3.4.3} Assume $A\in {\cal M}_{\mu}^{\a}$ and $B\in {\cal M}_{\lambda}^{\b}$. $\mu_y\wedge \lambda_y=s$, $\mu_x\wedge \lambda_x=t$, $\frac{\mu_y}{s}=m$, $\frac{\mu_x}{t}=n$, $\frac{\lambda_y}{s}=p$, $\frac{\lambda_x}{t}=q$. Since $s,~t$ are co-prime, denote $\sigma=s/t$, then $\sigma_y=s$ and $\sigma_x=t$.

Split $A$ as
$$
A=\begin{bmatrix}
A_{1,1}&A_{1,2}&\cdots&A_{1,n}\\
A_{2,1}&A_{2,2}&\cdots&A_{2,n}\\
\cdots&~&~&~\\
A_{m,1}&A_{m,2}&\cdots&A_{m,n}\\
\end{bmatrix},
$$
where $A_{i,j}\in {\cal M}_{\sigma}^{\a}$;
and split $B$ as
$$
B=\begin{bmatrix}
B_{1,1}&B_{1,2}&\cdots&B_{1,q}\\
B_{2,1}&B_{2,2}&\cdots&B_{2,q}\\
\cdots&~&~&~\\
B_{p,1}&B_{p,2}&\cdots&B_{p,q}\\
\end{bmatrix},
$$
where $B_{i,j}\in {\cal M}_{\sigma}^{\b}$.
Then the generalized weighted inner product is defined as
\begin{align}\label{3.4.3}
\begin{array}{l}
(A\;\big|\;B)_W:=\\
\left[
\begin{array}{llll}
\left(A_{1,1}|B_{1,1}\right)_W&\cdots&\left(A_{1,1}|B_{1,q}\right)_W&\cdots\\
\left(A_{2,1}|B_{2,1}\right)_W&\cdots&\left(A_{2,1}|B_{2,q}\right)_W&\cdots\\
\vdots&~&~&~\\
\left(A_{m,1}|B_{p,1}\right)_W&\cdots&\left(A_{m,1}|B_{p,q}\right)_W&\cdots
\end{array}\right.\\
\left. \begin{array}{lll}
\left(A_{1,n}|B_{1,1}\right)_W&\cdots&\left(A_{1,n}|B_{1,q}\right)_W\\
\left(A_{2,n}|B_{2,1}\right)_W&\cdots&\left(A_{2,n}|B_{2,q}\right)_W\\
\vdots&~&~\\
\left(A_{m,n}|B_{p,1}\right)_W&\cdots&\left(A_{m,n}|B_{p,q}\right)_W
\end{array}\right],
\end{array}
\end{align}
where $\left(A_{i,j}|B_{r,s}\right)$ are defined in (\ref{0.5.2}).
\end{dfn}

\begin{dfn}\label{d3.4.4} Assume $\A\in \Sigma_{\mu}$ and $\B\in \Sigma_{\lambda}$. Then the generalized inner product of $\A$ and $\B$, denoted by $(\A\;|\; \B)$, is defined as
\begin{align}\label{3.4.4}
(\A\;|\; \B):=\left<(A\;|\; B)_W\right>.
\end{align}
\end{dfn}

Of course, we need to prove that (\ref{3.4.4}) is independent of the choice of representatives $A$ and $B$. This is verified by a straightforward computation.

Next, we would like to define another ``inner product" called the $\d$-inner product, where $\d\in \Q_+$.
First we introduce a new notation:

\begin{dfn}\label{3.4.401}
Let $\mu,\d\in \Q_+$, $\mu$ is said to be superior to $\d$, denoted by
$$
\mu\gg \d,
$$
if $\d_y|\mu_y$ and $\d_x|\mu_x$.
\end{dfn}

The $\d$ inner product is a mapping $(\cdot|\cdot): \bigcup_{\mu\gg \d}\Sigma_{\mu}\times \bigcup_{\mu\gg \d}\Sigma_{\mu} \ra \Sigma_{\d}$.

\begin{dfn}\label{d3.4.5} Assume $A\in {\cal M}^{\a}_{\mu}$, $B\in {\cal M}^{\b}_{\lambda}$ and $ \mu\gg \d$, $\lambda\gg \d$. Denote $\mu_y/\d_y=\xi$, $\mu_x/\d_x=\eta$, $\lambda_y/\d_y=\zeta$, and $\lambda_x/\d_x=\ell$, then the $\d$-inner product of $A$ and $B$ is defined as follows:
Split
$$
A=\begin{bmatrix}
A_{1,1}&A_{1,2}&\cdots&A_{1,\eta}\\
A_{2,1}&A_{2,2}&\cdots&A_{2,\eta}\\
\vdots&~&~&~\\
A_{\xi,1}&A_{\xi,2}&\cdots&A_{\xi,\eta}\\
\end{bmatrix},
$$
where $A_{i,j}\in {\cal M}_{\d}^{\a}$,
and
$$
B=\begin{bmatrix}
B_{1,1}&B_{1,2}&\cdots&B_{1,\ell}\\
B_{2,1}&B_{2,2}&\cdots&B_{2,\ell}\\
\vdots&~&~&~\\
B_{\zeta,1}&B_{\zeta,2}&\cdots&B_{\zeta,\ell}\\
\end{bmatrix},
$$
where $B_{i,j}\in {\cal M}_{\d}^{\b}$.
Then the $\d$-inner product of $A$ and $B$ is defined as
\begin{align}\label{3.4.5}
(A|B)_{\d}:=
\begin{bmatrix}
C_{1,1}&C_{1,2}&\cdots&C_{1,\eta}\\
C_{2,1}&C_{2,2}&\cdots&C_{2,\eta}\\
\vdots\\
C_{\xi,1}&C_{\xi,2}&\cdots&C_{\xi,\eta}\\
\end{bmatrix},
\end{align}
where
\begin{align}\label{3.4.6}
C_{i,j}:=
\begin{bmatrix}
(A_{i,j}|B_{1,1})_W&(A_{i,j}|B_{1,2})_W&\cdots&(A_{i,j}|B_{1,\ell})_W\\
(A_{i,j}|B_{2,1})_W&(A_{i,j}|B_{2,2})_W&\cdots&(A_{i,j}|B_{2,\ell})_W\\
\vdots\\
(A_{i,j}|B_{\zeta,1})_W&(A_{i,j}|B_{\zeta,2})_W&\cdots&(A_{i,j}|B_{\zeta,\ell})_W\\
\end{bmatrix}.
\end{align}
\end{dfn}

\begin{dfn}\label{d3.4.6} Assume $\A\in \Sigma_{\mu}$ and $\B\in \Sigma_{\lambda}$, where  $\mu\gg \d$ and $\lambda\gg \d$. Then the $\d$-inner product of $\A$ and $\B$ is defined as
\begin{align}\label{3.4.7}
(\A\;|\;\B)_{\d}:=(A\;|\;B)_{\d},\quad A\in \A,~B\in \B.
\end{align}
\end{dfn}

\begin{rem}\label{r3.4.7}
\begin{enumerate}
\item It is easy to verify that (\ref{3.4.7}) is independent of the choice of $A$ and $B$. Hence the $\d$-inner product is well defined.
\item Definition \ref{d3.4.3} ( or Definition \ref{d3.4.4} ) is a special case of  Definition \ref{d3.4.5} (correspondingly, Definition \ref{d3.4.6} ).
\item Definition \ref{d3.4.1} cannot be extended to the equivalence space, because it depends on the choice of representatives.
\item Unlike the generalized inner product defined in Definition \ref{d3.4.3} (as well as Definition \ref{d3.4.4} ), the $\d$-inner product is defined on a subset of ${\cal M}$ (or $\Sigma$).
\end{enumerate}
\end{rem}

Using $\d$-inner product, we have the following.

\begin{prp}\label{p3.4.8} Assume $\varphi:\Sigma_{\mu}\ra \Sigma_{\lambda}$ is a linear mapping. Then there exists a matrix $\Lambda\in {\cal M}_{r\mu_y\lambda_y\times r\mu_x\lambda_x}$, called the structure matrix of $\varphi$, such that
\begin{align}\label{3.4.8}
\varphi(\A)=\left(A\;|\;\L\right)_{\mu}.
\end{align}
\end{prp}

\subsection{Vector Fields}

\begin{dfn}\label{d3.5.1}  Let $M$ be a bundled manifold and $T(M)$ the tangent space of $M$. $V:M\ra T(M)$ is called a $C^r$ vector field, if for each simple coordinate chart $U$, $V|_U$ is $C^r$.
The set of $C^r$ vector fields on $M$ is denoted by $V^r(M)$.
\end{dfn}

We express a vector field in a matrix form. That is, let $X\in V^r\left({\cal M}_{m\times n}\right)$. Then
$$
X=\dsum_{i=1}^m\dsum_{j=1}^nf_{i,j}(x)\frac{\pa}{\pa x_{i,j}}:=\left[f_{i,j}(x)\right]\in {\cal M}_{m\times n}.
$$
Similar to smooth functions, the vector fields on ${\cal M}_{\mu}$ can be defined as follows:

\begin{dfn}\label{d3.5.2} Assume $\left<X\right>$ is firstly defined on $T\left({\cal M}_{\mu}^s\right)$ as $X_s(x)$, i.e., ${\cal M}_{\mu}^s$ is the root leaf of $\left<X\right>$. Then we extend it to other leafs as:
\begin{itemize}
\item Step 1. Let $Q=\{t\in \N\;\big|\; t|s\}$. Then for
\begin{align}\label{3.5.1}
X_t(y):=\left(\pr_k\right)_*(X_s)(x=\bd_k(y)),
\end{align}
where $k=\frac{s}{t}\in \N$.
\item Step 2. Assume $\ell\wedge s=t$. If $\ell=t$, $X_{\ell}=X_t$ has already been defined in Step 1. So we assume $\ell=kt$. Then\footnote{Let $M$ and $N$ be two manifolds, and $\varphi:M\ra N$  a smooth mapping, $x\in M$, $y\in N$, and $\varphi(x)=y$. Then
\begin{itemize} \item
$\varphi_*:~T_y(N)\ra T_x(M)$, satisfying
$$
\varphi_*(X)h(y)=L_X(h\circ \varphi)(x), \quad \forall h(y)\in C^{\infty}(N), \;\forall X\in T_x(M),
$$
where $L_X$ is the Lie derivative with respect to $X$;
\item
$\varphi^*:~T^*_x(M)\ra T^*_y(N)$, satisfying
$$
\varphi^*(\omega)Y=\omega(\varphi_*(Y)|_x), \quad \forall Y\in T_y(N), \;\forall \omega\in T^*_x(M).
$$
\end{itemize}
We refer readers to \cite{boo79} for concepts, and to \cite{spi79} for computations.
}
\begin{align}\label{3.5.2}
X_{\ell}(y):=
\left(\bd_k\right)_*(X_t)(x=\pr_k(y))\otimes I_k.
\end{align}
\end{itemize}
\end{dfn}

Next, we consider the computation of the related expressions of a vector field, which is originally defined on its root leaf.
\begin{itemize}
\item To calculate (\ref{3.5.1})  we first set
\begin{align}\label{3.5.3}
x=y\otimes I_k
\end{align}
to get $X_s(x(y)):=X_s(y)$. Then we split $X_s(y)$ into $tp\times tq$ blocks as $X_s=[X_s^{i,j}]$, where each $X_s^{i,j}\in {\cal M}_{k\times k}$. Then $X_t=[X_{i,j}]\in T_x\left({\cal M}_{tp\times tq}\right)$, and
\begin{align}\label{3.5.4}
X_{i,j}=\Tr\left([X_s^{i,j}]\right).
\end{align}

\item To calculate (\ref{3.5.2})  we split $y$ into $tp\times tq$ blocks as $y=[y^{i,j}]$, where each $y^{i,j}\in {\cal M}_{k\times k}$. Then $X_t(y)$ is obtained by replacing $x$ by $x=[x_{i,j}]\in {\cal M}_{tp\times tq}$ as
\begin{align}\label{3.5.5}
x_{i,j}=\Tr\left([y^{i,j}]\right).
\end{align}
It follows that
\begin{align}\label{3.5.6}
X_{\ell}(y)= X_t(y)\otimes I_k.
\end{align}
\end{itemize}

Then it is easy to verify the following:

\begin{prp}\label{p3.5.3} The vector field $X$ defined in Definition \ref{d3.5.2} is consistent with the equivalence $\sim\big|_{{\cal M}_{\mu}^{[s,\cdot]}}$. Hence the equivalent class $\left<X\right>$ is well defined on $T\left(\Sigma_{\mu}^{[s,\cdot]}\right)$.
\end{prp}

\begin{rem}\label{r3.5.301} In fact, Proposition \ref{p3.5.3} only claims that the representations on embedded supper spaces are consistent, which is obviously weaker than Proposition \ref{p3.3.3}. Please refer to Remark \ref{r3.2.4} and the latter Remark \ref{r3.6.5.4} for the extension of $\left<X\right>$ to the projected subspace.
\end{rem}

\begin{exa}\label{e3.5.4} Consider $\Sigma_{1/2}$. Assume $X$ is firstly defined on $T\left({\cal M}_{1/2}^2\right)$ as
$$
X_2(x)=F(x)=\begin{bmatrix}F^{11}(x)&F^{12}(x)\end{bmatrix},
$$
where $x=(x_{i,j})\in {\cal M}_{1/2}^2$ and
$$
F^{1,1}(x)=\begin{bmatrix}
x_{1,1}&0\\
0&x_{1,3}
\end{bmatrix};\quad
F^{1,2}(x)=\begin{bmatrix}
x_{2,2}&0\\
x_{1,1}&0
\end{bmatrix}.
$$
Then we consider the expression of $X$ on the other cross sections:
\begin{itemize}
\item Consider $X_1\in T\left({\cal M}_{1/2}^1\right)$. Set
$$
X_1(y)=\begin{bmatrix}f_{1}(y)&f_{2}(y)\end{bmatrix},
$$
where $y=(y_1,y_2)\in {\cal M}_{1/2}^1$.
According to (\ref{3.5.1}),
$$
\begin{array}{ccl}
f_{1}(y)&=&\frac{1}{2}\left(F^{1,1}_{1,1}(\bd_2(y))+F^{1,1}_{2,2}(\bd_2(y))\right)\\
~&=&\frac{y_1+y_2}{2};
\end{array}
$$
and
$$
\begin{array}{ccl}
f_{2}(y)&=&\frac{1}{2}\left(F^{1,2}_{1,1}(\bd_2(y))+F^{1,2}_{2,2}(\bd_2(y))\right)\\
~&=&\frac{1}{2}y_1.
\end{array}
$$
\item
Consider $X_3\in T\left({\cal M}_{1/2}^3\right)$. Set
$$
X_3(z)=\left(g_{i,j}(z)\right)\in {\cal M}_{1/2}^3,
$$
where $z=\left(z_{i,j}\right)\in {\cal M}_{1/2}^3$.

Consider the projection $\pr_3:~{\cal M}_{1/2}^3\ra {\cal M}_{1/2}^1$:
$$
\pr_3(z)=\begin{bmatrix}\frac{z_{1,1}+z_{2,2}+z_{3,3}}{3}& \frac{z_{1,4}+z_{2,5}+z_{3,6}}{3}\end{bmatrix}.
$$
According to (\ref{3.5.2}),
$$
X_3(z)=[G_{1}(z),G_{2}(z)]\otimes I_3,
$$
where
$$
\begin{array}{ccll}
G_1(z)&=&f_{1}(\pr_3(z))=&\frac{1}{6}\left(z_{1,1}+z_{2,2}+z_{3,3}\right.\\
~&~&~&\left.+z_{1,4}+z_{2,5}+z_{3,6}\right)\\
G_2(z)&=&f_{2}(\pr_3(z))=&\frac{1}{6}\left(z_{1,1}+z_{2,2}+z_{3,3}\right).
\end{array}
$$

Similarly, for $n=2k-1$ we have
$$
X_n(z)=[G_{1}(z),G_{2}(z)]\otimes I_n\in T\left({\cal M}_{1/2}^n\right),
$$
where
$$
\begin{array}{ccl}
G_1(z)&=&\frac{1}{2n}\left(z_{1,1}+z_{2,2}+\cdots+z_{n,n}\right.\\
~&~&+\left.z_{1,n+1}+z_{2,n+2}+\cdots+z_{n,2n}\right)\\
G_2(z)&=&=\frac{1}{2n}\left(z_{1,1}+z_{2,2}+\cdots+z_{n,n}\right).
\end{array}
$$

\item
Consider $X_4\in {\cal M}_{1/2}^4$. Set
$$
X_4(z)=\left(g_{i,j}(z)\right)\in T\left({\cal M}_{1/2}^4\right),
$$
where $z=\left(z_{i,j}\right)\in {\cal M}_{1/2}^4$.
Consider the projection $\pr_2:~{\cal M}_{1/2}^4\ra {\cal M}_{1/2}^2$:
$$
\pr_2(z)=\begin{bmatrix}
\frac{z_{1,1}+z_{2,2}}{2}&\frac{z_{1,3}+z_{2,4}}{2}&\frac{z_{1,5}+z_{2,6}}{2}&\frac{z_{1,7}+z_{2,8}}{2}\\
\frac{z_{3,1}+z_{4,2}}{2}&\frac{z_{3,3}+z_{4,4}}{2}&\frac{z_{3,5}+z_{4,6}}{2}&\frac{z_{3,7}+z_{4,8}}{2}\\
\end{bmatrix}.
$$
According to (\ref{3.5.2}),
$$
X_4(z)=\begin{bmatrix}
G_{1,1}(z)&G_{1,2}(z)&G_{1,3}(z)&G_{1,4}(z)\\
G_{2,1}(z)&G_{2,2}(z)&G_{2,3}(z)&G_{2,4}(z)\\
\end{bmatrix}\otimes I_2,
$$
where
$$
\begin{array}{cclccl}
G_{1,1}(z)&=&\frac{z_{1,1}+z_{2,2}}{2},&
G_{1,2}(z)&=&0,\\
G_{1,3}(z)&=&\frac{z_{3,3}+z_{4,4}}{2},&
G_{1,4}(z)&=&0,\\
G_{2,1}(z)&=&0,&
G_{2,2}(z)&=&\frac{z_{1,5}+z_{2,6}}{2},\\
G_{2,3}(z)&=&\frac{z_{1,1}+z_{2,2}}{2}&
G_{2,4}(z)&=&0.\\
\end{array}
$$

Similarly, for $n=2k$ we have $X_n\in T\left({\cal M}_{1/2}^n\right)$ as
$$
X_n(z)=\begin{bmatrix}
G_{1,1}(z)&G_{1,2}(z)&G_{1,3}(z)&G_{1,4}(z)\\
G_{2,1}(z)&G_{2,2}(z)&G_{2,3}(z)&G_{2,4}(z)\\
\end{bmatrix}\otimes I_k,
$$
where
$$
\begin{array}{ccl}
G_{1,1}(z)&=&\frac{z_{1,1}+z_{2,2}+\cdots+z_{k,k}}{k},\\
G_{1,2}(z)&=&0,\\
G_{1,3}(z)&=&\frac{z_{k+1,k+1}+z_{k+2,k+2}+\cdots+z_{2k,2k}}{k},\\
G_{1,4}(z)&=&0,\\
G_{2,1}(z)&=&0,\\
G_{2,2}(z)&=&\frac{z_{1,2k+1}+z_{2,2k+2}+\cdots+z_{k,3k}}{k},\\
G_{2,3}(z)&=&\frac{z_{1,1}+z_{2,2}+\cdots+z_{k,k}}{k},\\
G_{2,4}(z)&=&0.\\
\end{array}.
$$
\end{itemize}
\end{exa}

\subsection{Integral Curves}

\begin{dfn}\label{d3.6.5.1} Let $\left<\xi\right>\in V^r(\Sigma_{\mu})$ be a vector field. Then for each $\A\in \Sigma_{\mu}$ there exists a curve $\left<X(t)\right>$ such that $\left<X(0)\right>=\A$ and
\begin{align}\label{3.6.5.1}
\left<\dot{X}(t)\right>=\left<\xi\right>\left(\left<X(t)\right>\right),
\end{align}
which is called the integral curve of $\left<\xi\right>$, starting from $\A$.
\end{dfn}

In fact, the integral curve of $\left<\xi\right>$ is a bundled integral curve.
Assume $\A=\left<A_0\right>$, where $A_0\in {\cal M}_{m\times n}$ is irreducible. Here $m$, $n$ may not be co-prime but $m/n=\mu$. So we denote $A_0\in {\cal M}_{\mu}^s$. Then on this root leaf we denote
$$
\xi^s=\left<\xi\right>\cap T\left({\cal M}_{m\times n}\right),
$$
and the integral curve of $\xi^s$, starting from $A_0$ is a standard one, which is the cross section of the bundled integral curve $\left<X\right>$ passing through $A_0$. That is, it is the solution of
\begin{align}\label{3.6.5.2}
\begin{cases}
\dot{X}_s(t)=\xi^s(X_s),\\
X_s(0)=A_0.
\end{cases}
\end{align}
We may denote the solution as
\begin{align}\label{3.6.5.2}
X_s(t)=\Phi^{\xi^s}_t(A_0).
\end{align}

Next, we consider the other cross sections of the integral curve, which correspond to the cross sections of $\left<\xi\right>$ respectively.

Recall Definition \ref{d3.5.2}, we can get the cross sections of the bundled integral curve on each leafs by using the cross sections of $\left<\xi\right>$ on corresponding leafs. The following result is then obvious:

\begin{thm}\label{t3.6.5.2}
Assume $s|\tau$. The corresponding cross section of the bundled integral curve is the integral curve, denoted by
$$
X_{\tau}(t)=\Phi^{\xi^{\tau}}_t(\bd_k(A_0)),\quad k=\tau/s,
$$
satisfying
\begin{align}\label{3.6.5.3}
\begin{cases}
\dot{X}_{\tau}(t)=\xi^{\tau}\left(X_{\tau}\right),\\
X_{\tau}(0)=\bd_k(A_0),
\end{cases}
\end{align}
where  $\xi^{\tau}$ is defined by (\ref{3.5.2}).
\end{thm}

Note that a cross section here is a mapping from base space to a leave such that  (\ref{0.3.1}) holds.

\begin{figure}[!htbp]
\centering
\includegraphics[width=5cm]{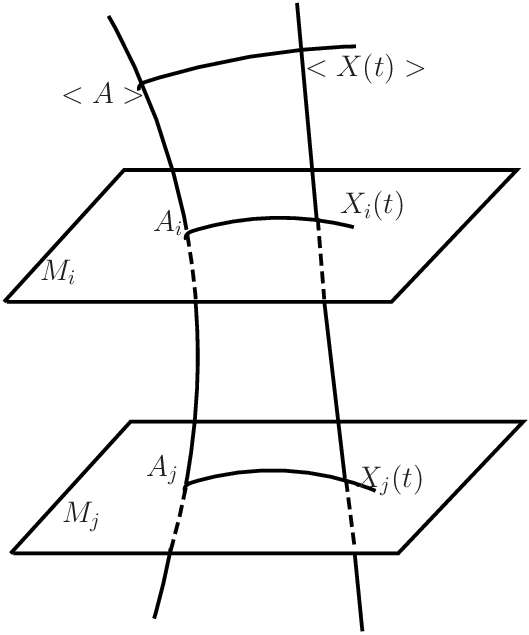}
\caption{Integral Curve of Vector Field}\label{Fig.7}
\end{figure}

 Fig.~\ref{Fig.7} shows the integral curve and its projections on each leaf.

The following result comes from the construction directly.

\begin{thm}\label{t3.6.5.3} Assume $A_{\a}=\A\cap {\cal M}_{\mu}^{\a}$, $A_{\b}=\A\cap {\cal M}_{\mu}^{\b}$,
$\xi^{\a}=\left<\xi\right>\cap T\left({\cal M}_{\mu}^{\a}\right)$, $\xi^{\b}=\left<\xi\right>\cap T\left({\cal M}_{\mu}^{\b}\right)$, and $\a=k\b$, $k\geq 2$. $\xi$ is defined firstly on $T\left({\cal M}_{\mu}^{s}\right)$, where $s|\b$.
Then there exists a one-to-one correspondence between the two cross sections (or corresponding integral curves).
Precisely,
\begin{align}\label{3.6.5.5}
\begin{array}{l}
\Phi^{\xi^{\a}}_t\left(A_{\a}\right)=\pr_k\left[\Phi^{\xi^{\b}}_t\left(\pr_k(A_{\a})\right)\right],\\
\Phi^{\xi^{\b}}_t\left(A_{\b}\right)=\bd_k\left[\Phi^{\xi^{\a}}_t\left(\bd_k(A_{\b})\right)\right].
\end{array}
\end{align}
\end{thm}

Particularly, assume the vector field is a linear vector field and $X_0\in {\cal M}_{\d}^1$, then we have
\begin{align}\label{3.6.5.6}
\Phi^{\xi^1}_{i,j}=\left(X^{k}_{i,j}\;\big|\; X_0\right)_{\d},\quad i=1,\cdots,m;j=1,\cdots,n.
\end{align}
Moreover, on ${\cal M}_{sm\times sn}$, the cross section of the integral curve can be expressed by modifying (\ref{3.6.5.6}) as
\begin{align}\label{3.6.5.7}
\Phi^{\xi}_{i,j}&=&\Phi^{\xi^1}_{i,j}\otimes I_s,\quad i=1,\cdots,m;j=1,\cdots,n.
\end{align}

\begin{rem}\label{r3.6.5.4}
\begin{itemize}
\item Note that if $\left<\xi\right>$ is firstly defined on $T\left({\cal M}_{\mu}^s\right)$, $t<s$ and $k=s/t\in \N$. Then $\xi^{(t)}$ is obtained through the following two steps: (i) Restrict $\xi^s$ on ${\cal M}_{\mu}^t\otimes I_k$ as
$$
\xi^s(y=\bd_k(x))=\xi^s\big|_{{\cal M}_{\mu}^t\otimes I_k}.
$$
(ii) Project $\xi^s(y=\bd_k(x))$ onto the tangent space of the subspace $T\left({\cal M}_{\mu}^t\otimes I_k\right)$.

 According to (ii), $\xi^{(t)}(A)$ does not correspond to $\xi^s(A\otimes I_k)$. Hence, the relationship demonstrated in Theorems \ref{t3.6.5.2} and \ref{t3.6.5.3} are not available for the integral curves of $\xi^{(t)}$ and $\xi^s$ respectively.

Because of this argument, if $\left<\xi\right>$ is firstly defined on $T\left({\cal M}_{\mu}^s\right)$, where ${\cal M}_{\mu}^s$ is the root leaf for $\left<\xi\right>$,  then we are mainly concerning the integral curves of $\xi^{\tau}$ on $T\left({\cal M}_{\mu}^{\tau}\right)$, where $s|\tau$.

\item The projections of $\left<\xi\right>=\left<\xi^s\right>$ onto subspaces $\xi^{(t)}$ are useful in some other problems. For instance, assume $\left<\xi\right>$ is a constant vector field. Using notations in Theorem \ref{t3.6.5.3}, we set
    $$
    {\cal D}:=\Span\left\{D^{I,J}\;\big|\; I=1,\cdots,\a p;J=1,\cdots,\a q\right\},
    $$
    where $D^{I,J}$ is defined in (\ref{0.4.4}). Then the projection of $\xi^{\b}$ on ${\cal D}$, denoted by
    $\xi^{\b}_{{\cal D}}$, is well defined. Moreover, (\ref{3.6.5.5}) becomes
\begin{align}\label{3.6.5.701}
\begin{array}{l}
\Phi^{\xi^{\a}}_t\left(A_{\a}\right)=\pr_k\left[\Phi^{\xi^{\b}_{{\cal D}}}_t\left(\pr_k(A_{\a})\right)\right],\\
\Phi^{\xi^{\b}_{{\cal D}}}_t\left(A_{\b}\right)=\bd_k\left[\Phi^{\xi^{\a}}_t\left(\bd_k(A_{\b})\right)\right].
\end{array}
\end{align}
\end{itemize}
\end{rem}

\begin{exa}\label{e3.6.5.5} Recall Example \ref{e3.3.4}. Since it is a linear vector field, it is easy to calculate that the cross section on ${\cal M}_{2\times 4}$ can be expressed as in (\ref{3.6.5.6}), where
$X_0\in {\cal M}_{2\times 4}$ and
$$
\begin{array}{l}
X^2_{1,1}=
\begin{bmatrix}
e^t&0&0&0\\
0&0&0&0
\end{bmatrix};~
X^2_{1,2}=
\begin{bmatrix}
0&1&0&0\\
0&0&0&0
\end{bmatrix};\\
X^2_{1,3}=
\begin{bmatrix}
0&0&\frac{e^t+e^{-t}}{2}&0\\
0&\frac{e^t-e^{-t}}{2}&0&0
\end{bmatrix};~
X^2_{1,4}=
\begin{bmatrix}
0&0&0&1\\
0&0&0&0
\end{bmatrix};\\
X^2_{2,1}=
\begin{bmatrix}
0&0&0&0\\
1&0&0&0
\end{bmatrix};~
X^2_{2,2}=
\begin{bmatrix}
0&0&\frac{e^t-e^{-t}}{2}&0\\
0&\frac{e^t+e^{-t}}{2}&0&0
\end{bmatrix};\\
X^2_{2,3}=
\begin{bmatrix}
0&0&0&0\\
t&0&1&0
\end{bmatrix};~
X^2_{2,4}=
\begin{bmatrix}
0&0&0&0\\
0&0&0&1
\end{bmatrix}.
\end{array}
$$

Consider the cross section on ${\cal M}_{2k\times 4k}$ with $X^k_0=X_0\otimes I_k\in \left<X_0\right>$. Then
$$
X^{2k}_{i,j}=X^2_{i,j}\otimes I_k,\quad i=1,2;~j=1,2,3,4.
$$
\end{exa}

An integral manifold of an involutive distribution on $T(\Sigma_{\mu})$ can be defined and calculated in a similar way.

\subsection{Forms}

\begin{dfn}\label{d3.6.6.1}  Let $M$ be a bundled manifold. ~$\omega:~M\ra T^*(M)$ is called a $C^r$ co-vector
field (or one form), if for each simple coordinate chart $U$ $\omega|_U$ is a $C^r$ co-vector field.
The set of $C^r$ co-vector fields on $M$ is denoted by ${V^*}^r(M)$.
\end{dfn}

We express a co-vector field in matrix form. That is, let $\omega\in {V^*}\left({\cal M}_{m\times n}\right)$. Then
$$
\omega=\dsum_{i=1}^m\dsum_{j=1}^n\omega_{i,j}(x)d x_{i,j}:=\left[\omega_{i,j}(x)\right]\in {\cal M}_{m\times n}.
$$

Similar to the construction of vector fields, the co-vector fields on $\Sigma_{\mu}$ can be established as follows:

\begin{dfn}\label{d3.6.6.2} Assume $\left<\omega\right>$ is firstly defined on $T^*\left({\cal M}_{\mu}^s\right)$ as $\omega_s(x)$, where ${\cal M}_{\mu}^s$ is the root leaf of $\left<\omega\right>$.  Then we extend it to other leafs as:
\begin{itemize}
\item Step 1. Let $Q=\{t\in \N\;\big|\; t|s\}$.
\begin{align}\label{3.6.6.1}
\omega_t(y):=k\left[\left(\bd_k\right)^*(\omega_s)(x=\bd_k(y))\right],
\end{align}
where $k=\frac{s}{t}\in \N$.
\item Step 2. Assume $\ell\wedge s=t$. If $\ell=t$, $\omega_{\ell}=\omega_t$ has already been defined in Step 1. So we assume $\ell=kt$, $k\geq 2$. Then
\begin{align}\label{3.6.6.2}
\omega_{\ell}(y):=\frac{1}{k}\left[\left(\pr_k\right)^*(\omega_t)(x=\pr_k(y))\otimes I_k\right].
\end{align}
\end{itemize}
\end{dfn}

Similar to the calculation of vector fields, to calculate (\ref{3.6.6.1})  we first set
\begin{align}\label{3.6.6.201}
x=y\otimes I_k
\end{align}
to get $\omega_s(x(y)):=\omega_s(y)$. Then we split $\omega_s(y)$ into $tp\times tq$ blocks as $\omega_s=[\omega_s^{i,j}]$, where each $\omega_s^{i,j}\in {\cal M}_{k\times k}$. Then $\omega_t=[\omega_{i,j}]\in T^*_x\left({\cal M}_{tp\times tq}\right)$, where the entries are
\begin{align}\label{3.6.6.202}
\omega_{i,j}=\Tr\left([\omega_s^{i,j}]\right),\quad i=1,\cdots,tp; j=1,\cdots,tq.
\end{align}

To calculate (\ref{3.6.6.2})  we split $y$ into $tp\times tq$ blocks as $y=[y^{i,j}]$, where each $y^{i,j}\in {\cal M}_{k\times k}$. Then $\omega_t(y)$ is obtained by replacing $x$ by $x=[x_{i,j}]\in {\cal M}_{tp\times tq}$ as
\begin{align}\label{3.6.6.203}
x_{i,j}=\Tr\left([y^{i,j}]\right).
\end{align}
It follows that
\begin{align}\label{3.6.6.204}
\omega_{\ell}(y)=\frac{1}{k}\left(\omega_t(y)\otimes I_k\right).
\end{align}

Then it is easy to verify the following:

\begin{prp}\label{p3.6.6.3} The co-vector field $\left<\omega\right>$ defined by Definition \ref{d3.6.6.2} is consistent with the equivalence $\sim\big|_{{\cal M}_{\mu}^{[s,\cdot]}}$. Hence it is well defined on $\Sigma_{\mu}^{[s,\cdot]}$.
\end{prp}

\begin{dfn}\label{p3.6.6.4} Let $\left<\omega\right>\in V^*\left(\Sigma_{\mu}\right)$ and $\left<X\right>\in V\left(\Sigma_{\mu}\right)$. Then the action of  $\left<\omega\right>$ on $\left<X\right>$ is defined as
\begin{align}\label{3.6.6.3}
 \left<\omega\right>(\left<X\right>):= \left(\left<\omega\right>\;\big|\; (\left<X\right>\right)_W.
\end{align}
\end{dfn}

Similar to vector field case, if the co-vector field is firstly defined on ${\cal M}_{\mu}^s$, then we can assume $\left<\omega\right>$ is only defined on ${\cal M}_{\mu}^{\tau}$ satisfying
$$
{\cal M}_{\mu}^{s}\sqsubset {\cal M}_{\mu}^{\eta}.
$$

\subsection{Tensor Fields}

The set of tensor fields on $\Sigma_{\mu}$ of covariant order $\a$ and contravariant order $\b$ is denoted by
$ {\bf T}^{\a}_{\b}(\Sigma_{\mu})$. To avoid complexity, we consider only the covariant tensor, $\left<t\right>\in {\bf T}^{\a}(\Sigma_{\mu})$.

\begin{dfn}\label{d3.6.7.1}
A covariant tensor field $\left<t\right>\in {\bf T}^{\a}\left(\Sigma_{\mu}\right)$ is a multi-linear mapping
$$
\left<t\right>: \underbrace{V^r(\Sigma_{\mu})\times\cdots\times V^r(\Sigma_{\mu})}_{\a}\ra C^r(\Sigma_{\mu}).
$$
\end{dfn}

Assume $\left<t\right>$ is defined on root leaf at $x\in {\cal M}_{\mu}^{s}$ as $t_s\in {\bf T}^{\a}\left({\cal M}_{\mu}^{s}\right)$,  $p,~q$ are co-prime and $p/q=\mu$.  The calculation is performed as follows:
 Construct the structure matrix of $t_s$ as
\begin{align}\label{3.6.7.1}
\begin{array}{l}
M_s(x):=\\
\begin{bmatrix}
t^{1,\cdots,1}_{1,\cdots,1}(x)&t^{1,\cdots,1}_{1,\cdots,2}(x)&\cdots&t^{1,\cdots,1}_{sq,\cdots,sq}(x)\\
t^{1,\cdots,2}_{1,\cdots,1}(x)&t^{1,\cdots,2}_{1,\cdots,2}(x)&\cdots&t^{1,\cdots,2}_{sq,\cdots,sq}(x)\\
\vdots&~&~&~\\
t^{sp,\cdots,sp}_{1,\cdots,1}(x)&t^{sp,\cdots,sp}_{1,\cdots,2}(x)&\cdots&t^{sp,\cdots,sp}_{sq,\cdots,sq}(x)\\
\end{bmatrix}\\
\in {\cal M}_{(sp)^{\a}\times (sq)^{\a}},
\end{array}
\end{align}
where
$$
\begin{array}{l}
t^{i_1,\cdots,i_{\a}}_{j_1,\cdots,j_{\a}}(x)
=\left.t_s\left(\frac{\pa}{\pa x_{i_1,j_1}},\frac{\pa}{\pa x_{i_2,j_2}},\cdots,\frac{\pa}{\pa x_{i_{\a},j_{\a}}}\right)\right|_x,\\
i_d=1,\cdots,sp;~j_d=1,\cdots,sq;~d=1,\cdots,\a;\\
x\in {\cal M}_{sp\times sq}.
\end{array}
$$

Consider $\left<X^k\right>\in V\left(\Sigma_{\mu}\right)$ with its irreducible element $X_s^k\in V\left({\cal M}_{\mu}^s\right)$, where ${\cal M}_{\mu}^s$ is the root leaf of $\left<X^k\right>$.   $X^k_s$ is expressed in matrix form as
$$
X_s^k:=\dsum_{i=1}^{sp}\dsum_{j=1}^{sq}v^k_{i,j}\frac{\pa}{\pa x_{i,j}}
:=\left[v^k_{ij}\right]:=V_s^k,\quad k=1,\cdots.\a.
$$

Then we have

\begin{prp}\label{p3.6.7.2}  \begin{align}\label{3.6.7.2}
 \begin{array}{l}
t_s\left(\left<X\right>^1,\cdots,\left<X\right>^{\a}\right)\\
~=\left(M_s(x)\;\big|\; V^1\otimes V^2\otimes \cdots \otimes V^{\a}\right)\\
~=\left(\cdots \left(\cdots \left(M_s\;\big|\; V^1\right)_W\;\big|\; V^2\right)_W\cdots\;\big| V^{\a}\right)_W.
\end{array}
\end{align}
\end{prp}

Next, we calculate the expressions of $\left<t\right>$ on other leafs. The following algorithm can be verified to be consistent on ${\cal M}_{\mu}^{[s,\cdot]}$, where ${\cal M}_{\mu}^{s}$ is the root leaf of $\left<t\right>$.

\begin{alg}\label{a3.6.7.3} Assume $\left<t\right>$ is firstly defined on $T^{\a}\left({\cal M}_{\mu}^s\right)$ as $t_s(x)$. Then we extend it to other leafs as:
\begin{itemize}
\item Step 1. Denote $Q:=\left\{r\;\big|\; r|s\right\}$.
Let $\tau\in Q$ and $k=\frac{s}{\tau}\in N$. Then
\begin{align}\label{3.6.7.201}
t_{\tau}(y)=\left(\bd_k\right)^*(t_s)(\bd_k(y))
\end{align}
can be calculated by
constructing $M_{\tau}(y)$ ($y\in {\bf T}^{\a}\left({\cal M}_{\tau p\times \tau q}\right)$) as follows: First, split $M_s(\bd(y))$ into $(\tau p)^{\a}\times (\tau q)^{\a}$ blocks
$$
\begin{array}{l}
M_s(\bd_r(y))=\\
\begin{bmatrix}
T^{1,\cdots,1}_{1,\cdots,1}(y)&T^{1,\cdots,1}_{1,\cdots,2}(y)&\cdots&T^{1,\cdots,1}_{rq,\cdots,\tau q}(y)\\
T^{1,\cdots,2}_{1,\cdots,1}(y)&T^{1,\cdots,2}_{1,\cdots,2}(y)&\cdots&T^{1,\cdots,2}_{rq,\cdots,\tau q}(y)\\
\vdots&~&~&~\\
T^{\tau p,\cdots,\tau p}_{1,\cdots,1}(y)&T^{\tau p,\cdots,\tau p}_{1,\cdots,2}(y)&\cdots&t^{\tau p,\cdots,\tau p}_{\tau q,\cdots,\tau q}(y)\\
\end{bmatrix},
\end{array}
$$
where each block $T^{i_1,\cdots,i_{\a}}_{j_1,\cdots,j_{\a}}\in {\cal M}_{k\times k}$.
Then we set
\begin{align}\label{3.6.7.3}
M_r(y):=k^{\a}[\xi^{i_1,\cdots,i_{\a}}_{j_1,\cdots,j_{\a}}(y)]\in {\cal M}_{(\tau p)^{\a}\times (\tau q)^{\a}},
\end{align}
where
$$
\xi^{i_1,\cdots,i_{\a}}_{j_1,\cdots,j_{\a}}(y)=\Tr\left(T^{i_1,\cdots,i_{\a}}_{j_1,\cdots,j_{\a}}\right).
$$
\item Step 2. Assume $\ell\wedge s=\tau$. If $\ell=\tau$, $t_{\ell}=t_{\tau}$ has already been defined in Step 1. So we assume $\ell=k\tau$, where $k\geq 2$. Then
\begin{align}\label{3.6.7.301}
t_{\ell}(z)=\left(\pr_k\right)^*(t_{\tau})(\pr_k(y))
\end{align}
can be calculated by constructing
\begin{align}\label{3.6.7.4}
M_{\ell}(z):=\frac{1}{k^{\a}}\left[M_{\tau}(y=\pr_k(z))\otimes I_{k^{\a}}\right],\quad z\in {\cal M}_{\ell p\times \ell q}.
\end{align}
\end{itemize}
\end{alg}

\begin{exa}\label{e3.6.7.4} Consider a covariant tensor field $\left<t\right>\in {\cal T}^2\left(\Sigma_{\mu}\right)$, where $\mu=\frac{2}{3}$. Moreover, $\left<t\right>$ is firstly defined at $x=(x_{i,j})\in {\cal M}_{2\times 3}$ with its structure matrix as
\begin{align}\label{3.6.7.5}
M_1(x)=\begin{bmatrix}
1&x_{12}&0&0&0&1&0&0&x_{22}\\
0&1&0&-1&0&0&0&0&0\\
1&0&-1&0&0&0&1&-1&0\\
0&1&0&0&0&x_{22}&0&0&x_{23}
\end{bmatrix}.
\end{align}
\begin{enumerate}
\item Let $X(\left<x\right>),Y(\left<x\right>)\in V\left(\Sigma_{\mu}\right)$ be defined at $x\in {\cal M}_{2\times 3}$ as
$$
X(x)=\begin{bmatrix}
x_{13}&0&0\\
0&0&x_{21}
\end{bmatrix},\quad
Y(x)=\begin{bmatrix}
1&0&1\\
0&x_{11}^2&0
\end{bmatrix}.
$$
Evaluate $\left<t\right>(X,Y)$.

First, using (\ref{3.4.3}), we calculate that
$$
\begin{array}{ccl}
t_1(X,\cdot)&=&\left(M_1(x)\;\big|\; X\right)_W\\
~&=&\begin{bmatrix}
x_{13}+x_{21}&x_{12}x_{13}-x_{21}&0\\
0&x_{13}&x_{21}x_{23}
\end{bmatrix}.
\end{array}
$$
Then
$$
t_1(X,Y)=x_{13}\left(1+x_{11}^2\right)+x_{21}.
$$

\item Expressing $\left<t\right>$ on leaf ${\cal M}_{4\times 6}$:
Note that
$$
x=\pr_2(y)=\begin{bmatrix}
\frac{y_{11}+y_{22}}{2}&\frac{y_{13}+y_{24}}{2}&\frac{y_{15}+y_{26}}{2}\\
\frac{y_{31}+y_{42}}{2}&\frac{y_{33}+y_{44}}{2}&\frac{y_{35}+y_{46}}{2}
\end{bmatrix},
$$
where $y=\left[y_{i,j}\right]\in {\cal M}_{4\times 6}$.

Using (\ref{3.6.7.4}), we have
$$
M_2(y)=\frac{1}{2^2}
\begin{bmatrix}
t^{11}&t^{12}&t^{13}\\
t^{21}&t^{22}&t^{23}\\
\end{bmatrix},
$$
where
$$
\begin{array}{ccl}
t^{11}&=&\begin{bmatrix}
1&\frac{y_{13}+y_{24}}{2}&0\\
0&1&0\\
\end{bmatrix}\otimes I_2;
\\
t^{12}&=&\begin{bmatrix}
0&0&1\\
-1&0&0\\
\end{bmatrix}\otimes I_2;
\\
t^{13}&=&\begin{bmatrix}
0&0&\frac{y_{33}+y_{44}}{2}\\
0&0&0\\
\end{bmatrix}\otimes I_2;
\\
t^{21}&=&\begin{bmatrix}
1&0&-1\\
0&1&0\\
\end{bmatrix}\otimes I_2;
\\
t^{22}&=&\begin{bmatrix}
0&0&0\\
0&0&\frac{y_{33}+y_{44}}{2}\\
\end{bmatrix}\otimes I_2;
\\
t^{23}&=&\begin{bmatrix}
1&-1&0\\
0&0&\frac{y_{35}+y_{46}}{2}\\
\end{bmatrix}\otimes I_2.
\end{array}
$$
\end{enumerate}

\end{exa}
%


\section{Lie Algebra on Square M-equivalence Space}

\subsection{Ring Structure on ~$\Sigma$}

Consider the vector space of the equivalent classes of square matrices $\Sigma:=\Sigma_{1}$. Since $\Sigma$ is closed under the STP, more algebraic structures may be posed on it. First, polynomials; Second, Lie algebra structure.

To begin with, we extend some fundamental concepts of matrices to their equivalent classes.

\begin{dfn}\label{d4.1.1}
\begin{enumerate}
\item $\A$  is nonsingular (symmetric, skew symmetric, positive/negative (semi-)definite, upper/lower (strictly) triangular, diagonal, etc.) if its irreducible element $A_1$ is (equivalently, every $A_i$ is).

\item $\A$ and $\B$ are similar, denoted by $\A\sim \B$, if there exists a nonsingular $\left<P\right>$ such that
    \begin{align}\label{4.1.1}
    \left<P^{-1}\right>\A\left<P\right>=\B.
    \end{align}
\item $\A$ and $\B$ are congruent, denoted by $\A\simeq \B$,  if there exists a nonsingular $\left<P\right>$ such that
    \begin{align}\label{4.1.2}
    \left<P^{T}\right>\A\left<P\right>=\B.
    \end{align}
\item $\left<J\right>$ is called the Jordan normal form of $\A$, if the irreducible element $J_1\in \left<J\right>$ is the Jordan normal form of $A_1$.
\end{enumerate}
\end{dfn}

\begin{dfn}\label{d4.2.1} \cite{lan02} A set $R$ with two operators $+,~\times$ is a ring. If the followings hold:
\begin{enumerate}
\item $(R,+)$ is an Abelian group;
\item $(R,\times)$ is a monoid;
\item (Distributive Rule)
$$
\begin{array}{l}
(a+b)\times c=a\times c+b\times c\\
c\times (a+b)=c\times a+c\times b,\quad a,b,c\in R.
\end{array}
$$
\end{enumerate}
\end{dfn}

Observing ${\cal M}_1$, which consists of all square matrices, both $\lplus$ (including $\lminus$) and $\ltimes$ are well defined. Unfortunately, $\left( {\cal M}_1, \lplus   \right)$ is not a group because there is no unit element. Since both $\lplus$ and $\ltimes$ are consistent with the equivalence $\sim$, we consider
$\Sigma:={\cal M}_1/\sim$. Then it is easy to verify the following:

\begin{prp}\label{p4.2.2} $\left(\Sigma,\lplus,\ltimes\right)$ is a ring.
\end{prp}

Consider a polynomial on a ring $R$, $p:R\ra R$, as
\begin{align}\label{4.2.1}
\begin{array}{ccl}
p(x)&=&a_nx^n+a_{n-1}x^{n-1}+\cdots+a_0,\\
~&~&~~~ a_i\in R,\;i=0,1,\cdots,n.
\end{array}
\end{align}
It is obvious that this is well defined. Set $R=(\Sigma,\lplus,\ltimes)$, then $p(x)$ defined in (\ref{4.2.1}) is also well defined on $R$. Particularly, the coefficients $a_i$ can be chosen as $\left<a_i\right>$ for $a_i\in \F$, $i=1,\cdots,n$. Then $p(x)$ is as a ``standard" polynomial. Unless elsewhere is stated, in this paper only such standard polynomials are considered.

For any $\A\in \Sigma$, the polynomial $p(\A)$ is well defined, and it is clear that
\begin{align}\label{4.2.2}
p(\A)=\left<p(A)\right>,\quad \mbox{for any}~A\in \A.
\end{align}

Using Taylor expansion, we can consider general matrix functions. For instance, we have the following result:
\begin{thm}\label{t4.2.3} Let $f(x)$ be an analytic function. Then $f(\A)$ is well defined provided $f(A)$ is well defined. Moreover,
\begin{align}\label{4.2.3}
f(\A)=\left<f(A)\right>.
\end{align}
\end{thm}

In fact, the above result can be extended to multi-variable case.
\begin{dfn}\label{d4.2.4} Let $F(x_1,\cdots,x_k)$ be a $k$-variable analytic function. Then $F(\A_1,\cdots,\A_k)$ is a well posed expression, where $\A_i\in \Sigma$, $i=1,\cdots,k$. Assume $A_i\in \A_i$, $i=1,\cdots,k$, then $F(A_1,\cdots,A_k)$ is a realization of $F(\A_1,\cdots,\A_k)$. Particularly, if $A_i\in {\cal M}_{r\times r}$, $\forall i$, we call  $F(A_1,\cdots,A_k)$ a realization of $F(\A_1,\cdots,\A_k)$ on $r$-th leaf.
\end{dfn}
Similar to (\ref{4.2.3}), we also have
\begin{align}\label{4.2.4}
F(\A_1,\cdots,\A_k)=\left<F(A_1,\cdots,A_k)\right>,
\end{align}
provided $F(A_1,\cdots,A_k)$ is well defined.

In the following we consider some fundamental matrix functions for $\Sigma$. We refer to \cite{cur84} for the definitions and basic properties of some fundamental matrix functions. Using these acknowledges, the following results are obvious:

\begin{thm}\label{t4.2.5} Let $\A,~\B \in \Sigma$ (i.e., $A$, $B$ are square matrices). Then the followings hold:
\begin{enumerate}
\item Assume $A\ltimes B=B\ltimes A$, then
\begin{align}\label{4.2.5}
e^{\A}\ltimes e^{\B}=e^{\left<A\lplus B\right>}.
\end{align}
\item If $\A$ is real skew symmetric, then $e^{\A}$ is orthogonal.
\item Assume $B$ is invertible, we denote $\B^{-1}=\left<B^{-1}\right>$. Then
\begin{align}\label{4.2.6}
e^{\B^{-1}\A \B}=\B^{-1}\ltimes e^{\A}\ltimes \B.
\end{align}
\item Let $A,~B$ be closed enough to identity so that $\log(A)$ and $\log(B)$ are defined, and $A\ltimes B=B\ltimes A$. Then
\begin{align}\label{4.2.7}
\log(\A \ltimes \B)=\log(\A)\lplus \log(\B).
\end{align}
\end{enumerate}
\end{thm}

Many known results for matrix functions can be extended to $\Sigma$. For instance, it is easy to prove the following Euler fromula:

\begin{prp}\label{p4.2.6} Consider $\F=\R$ and let $\A\in \Sigma$. Then the Euler formula holds. That is,
\begin{align}\label{4.2.8}
e^{i\A}=\cos(\A)~\lplus~ i\sin(\A).
\end{align}
\end{prp}

%
%
%
%
%

Recall the modification of trace and determinant in Definitions \ref{d4.2.7} and \ref{d4.2.9}. The following proposition shows that for the modifications the relationship between $\tr(A)$ and $\det(A)$ \cite{cur84} remains available.

\begin{prp}\label{p4.2.10} Assume $\F=\R$ (or $\F=\C$), $\A\in \Sigma$, then
\begin{align}\label{4.2.11}
e^{\Tr(\A)}=\Dt\left( e^{\A}\right).
\end{align}
\end{prp}

\noindent\textit{Proof.}
Let $A_0\in \A$ and $A_0\in {\cal M}_{n\times n}$. Then
$$
\begin{array}{ccl}
e^{\Tr(\A)}&=&e^{\frac{1}{n}\tr(A_0)}=\left(e^{\tr(A_0)}\right)^{\frac{1}{n}}\\
~&=&\left(|e^{\tr(A_0)}|\right)^{\frac{1}{n}}=\Dt(e^{A_0})\\
~&=&\Dt\left(e^{\A}\right).
\end{array}
$$
\hfill $\Box$

Next, we consider the characteristic polynomial of an equivalent class, they comes from standard matrix theory~\cite{hor85}.

\begin{dfn}\label{d4.2.11} Let $\A\in \Sigma$. $A_1\in \A$ is its irreducible element. Then
\begin{align}\label{4.2.12}
p_{\A}(\lambda):=\det(\lambda \lminus A_1)
\end{align}
is called the characteristic polynomial of $\A$.
\end{dfn}

The following result is an immediate consequence of the definition.
\begin{thm}[Cayley-Hamilton]\label{t4.2.12} Let $p_{\A}$ be the characteristic polynomial of $\A$. Then
\begin{align}\label{4.2.13}
p_{\A}(\A)=0.
\end{align}
\end{thm}

\begin{rem}\label{r4.2.13}
\begin{enumerate}
\item If we choose $A_k=A_1\otimes I_k$ and calculate the characteristic polynomial of $p_{A^k}$, then $p_{A_k}(\lambda)=\left(p_{A_1}(\lambda)\right)^k$. So $p_{A_k}(\A)=0$ is equivalent to $p_{A_1}(\A)=0$.
\item Choosing any $A_i\in \A$, the corresponding minimal polynomials $q_{A_i}(\lambda)$ are the same. So we have unique minimal polynomial as $q_{\A}(\lambda)=q_{A_i}(\lambda)$.
\end{enumerate}
\end{rem}

\subsection{Bundled Lie Algebra}

Consider the vector space of the equivalent classes of square matrices $\Sigma:=\Sigma_{1}$,  this section gives a Lie algebraic structure to it.

\begin{dfn}[\cite{hal03}]\label{d4.3.1} A Lie algebra is a vector space $g$ over some field $\F$ with a binary operation $[\cdot,\cdot]:~g\times g\ra g$, satisfying
\begin{enumerate}
\item (bi-linearity)
\begin{align}\label{4.3.1}
\begin{array}{ccl}
[\alpha A+\beta B,C]&=&\alpha [A,C]+\beta [B,C];\cr
[C,\alpha A+\beta B]&=&\alpha [C,A]+\beta [C,B],
\end{array}
\end{align}
where $\a,~\b\in \F$.
\item (skew-symmetry)
\begin{align}\label{4.3.2}
[A,B]=-[B,A];
\end{align}
\item (Jacobi Identity)
\begin{align}\label{4.3.3}
\begin{array}{r}
[A,[B,C]]+[B,[C,A]]+[C,[A,B]]=0,\\
\forall A,B,C\in g.
\end{array}
\end{align}
\end{enumerate}
\end{dfn}

\begin{dfn}\label{d4.3.2} Let $(E,\PR, B)$ be a discrete bundle with leaves $E_i$, $i=1,2,\cdots$. If
\begin{enumerate}
\item $(B,\oplus,\otimes)$ is a Lie algebra;
\item $(E_i,+,\times)$ is a Lie algebra, $i=1,2,\cdots$;
\item The restriction $\PR|_{E_i}: ~(E_i,+,\times) \ra (B,\oplus,\otimes)$ is a Lie algebra homomorphism, $i=1,2,\cdots$,
\end{enumerate}
then $(B,\oplus,\otimes)$ is called a bundled Lie algebra.
\end{dfn}
(We refer to Definition \ref{d0.3.1} and Remark \ref{d0.3.3} for the concept of discrete bundle.)

On vector space $\Sigma$ we define an operation $[\cdot,\cdot]:~\Sigma\times \Sigma\ra \Sigma$ as
\begin{align}\label{4.3.4}
[\A,\B]:=\A\ltimes \B\lminus \B\ltimes \A.
\end{align}
Then we have the following Lie algebra:

\begin{thm}\label{t4.3.3} The vector space ${\Sigma}$ with Lie bracket $[\cdot,\cdot]$ defined in (\ref{4.3.4}), is  a bundled Lie algebra, denoted by $\gl(\F)$.
\end{thm}

\noindent\textit{Proof}. Let ${\cal M}_1=\cup_{i=1}^{\infty}{\cal M}^i$, where ${\cal M}^i={\cal M}_1^i$. Then it is clear that $\left({\cal M}_1,\PR, \Sigma\right)$ is a discrete bundle.

Next, we prove $\left(\Sigma, \lplus, [\cdot,\cdot]\right)$ is a Lie algebra.  Equations (\ref{4.3.1}) and (\ref{4.3.2}) are obvious. We prove (\ref{4.3.3}) only.

Assume $A_1\in \A$, $B_1\in \B$ and $C_1\in \left<C\right>$ are irreducible, and $A_1\in {\cal M}_{m\times m}$, $B_1\in {\cal M}_{n\times n}$, and $C_1\in {\cal M}_{r\times r}$. Let $t=n\vee m\vee r$. Then it is easy to verify that
\begin{align}\label{4.3.5}
\begin{array}{l}
[\A,[\B,\left<C\right>]]=\\
~\left< \left[(A_1\otimes I_{t/m}), [(B_1\otimes I_{t/n}),(C_1\otimes I_{t/r})]\right]\right>.
\end{array}
\end{align}
Similarly, we have
\begin{align}\label{4.3.6}
\begin{array}{l}
[\B,[\left<C\right>,\A]]=\\
~\left< \left[(B_1\otimes I_{t/n}), [(C_1\otimes I_{t/r}),(A_1\otimes I_{t/m})]\right]\right>.
\end{array}\\
\label{4.3.7}
\begin{array}{l}
[\left<C\right>,[\A,\B]]=\\
~\left< \left[(C_1\otimes I_{t/r}), [(A_1\otimes I_{t/m}),(B_1\otimes I_{t/n})]\right]\right>.
\end{array}
\end{align}
Since (\ref{4.3.3}) is true for any $A,B,C\in gl(t,\F)$, it is true for $A=A_1\otimes I_{t/m}$, $B=B_1\otimes I_{t/n}$, and $C=C_1\otimes I_{t/r}$. Using this fact and equations (\ref{4.3.5})--(\ref{4.3.7}), we have
$$
\begin{array}{l}
[\A,[\B,\left<C\right>]]\lplus [\B,[\left<C\right>,\A]]\lplus [\left<C\right>,[\A,\B]]\\
=\left< \left[(A_1\otimes I_{t/m}),[(B_1\otimes I_{t/n}),(C_1\otimes I_{t/r})]\right]\right.\\
+~~\left[(B_1\otimes I_{t/n}), [(C_1\otimes I_{t/r}),(A_1\otimes I_{t/m})]\right]~~\\
+\left. \left[(C_1\otimes I_{t/r}), [(A_1\otimes I_{t/m}),(B_1\otimes I_{t/n})]\right]\right>\\
=\left<0\right>=0.
\end{array}
$$

 Let the Lie algebraic structure on ${\cal M}^n$ be $\gl(n,\F)$ (where $\F=\R$ or $\F=\C$). It follows from the consistence of $\lplus$ and $\ltimes$ with the equivalence that $\PR: \gl(n,\F) \ra \gl(\F)$ is a Lie algebra homomorphism.
\hfill $\Box$

\subsection{Bundled Lie Sub-algebra}

This section considers some useful Lie sub-algebras of Lie algebra $\gl(\F)$. Assume $g$ is a Lie algebra, $h\subset g$ is a vector subspace. Then $h$ is called a Lie sub-algebra if and only if $[h,h]\subset h$.

\begin{dfn}\label{d4.4.1} Let $(E,\PR, B)$ be a bundled Lie algebra. if
\begin{enumerate}
\item $(H,\oplus,\otimes)$ is a Lie sub-algebra of $(B,\oplus,\otimes)$;
\item  $(F_i,+,\times)$ is a Lie sub-algebra of  $(E_i,+,\times)$, $i=1,2,\cdots$;
\item The restriction $\PR|_{F_i}: ~(F_i,+,\times) \ra (H,\oplus,\otimes)$ is a Lie algebra homomorphism, $i=1,2,\cdots$,
\end{enumerate}
then $(F,\PR,H)$ is called a bundled Lie sub-algebra of $(E,\PR,B)$.
\end{dfn}

It is well known that there are some useful Lie sub-algebras of Lie algebra $\gl(n,\F)$. When $\gl(n,\F)$, $\forall n$, are generalized to the bundled Lie algebra $\gl(\F)$ (over $\Sigma$), the corresponding Lie sub-algebras are investigated one-by-one in this section.

\begin{itemize}

\vskip 2mm
\item

Bundled orthogonal Lie sub-algebra

\vskip 2mm

\begin{dfn}\label{d4.1.2}
$\A\in \Sigma$ is said to be symmetric (skew symmetric) if $A^T=A$ ($A^T=-A$), $\forall A\in \A$.
\end{dfn}

 The symmetric (skew symmetric) $\A$ is well defined because if $A\sim B$ and $A^T=A$ (or $A^T=-A$), then so is $B$. It is also easy to verify the following:

\begin{prp}\label{p4.4.2.2} Assume $\A$ and $\B$ are skew symmetric, then so is $[\A,\B]$.
\end{prp}

We, therefore, can define the following bundled Lie sub-algebra.

\begin{dfn}\label{d4.4.2.3}
$$
o(\F):=\left\{\A\in \gl(\F)\;|\;  \A^T=-\A\right\}
$$
is called the bundled orthogonal algebra.
\end{dfn}

\vskip 2mm
\item

Bundled special linear algebra

\vskip 2mm

\begin{dfn}\label{d4.4.2.4}
$$
sl(\F):=\left\{\A\in \gl(\F)\;|\;  \Tr(\A)=0\right\}
$$
is called the bundled special linear algebra.
\end{dfn}

Similar to the case of orthogonal algebra, it is easy to verify that $sl(\F)$ is a Lie sub-algebra of $\gl(\F)$.

\vskip 2mm

\item

Bundled upper triangular algebra

\vskip 2mm

\begin{dfn}\label{d4.4.2.5}
$$
t(\F):=\left\{\A\in \gl(\F)\;|\;  \A~\mbox{is upper triangular}\right\}
$$
is called the bundled upper triangular algebra.
\end{dfn}

Similarly, we can define bundled lower triangular algebras.

\vskip 2mm

\item

Bundled strictly upper triangular algebra

\vskip 2mm

\begin{dfn}\label{d4.4.2.6}
$$
n(\F):=\left\{\A\in \gl(\F)\;|\;  \A~\mbox{is strictly upper triangular}\right\}
$$
is called the bundled strictly upper triangular algebra.
\end{dfn}

\vskip 2mm

\item

Bundled diagonal algebra

\vskip 2mm

\begin{dfn}\label{d4.4.2.7}
$$
d(\F):=\left\{\A\in \gl(\F)\;|\;  \A~\mbox{is diagonal}\right\}
$$
is called the bundled diagonal algebra.
\end{dfn}

\vskip 2mm
\item

Bundled symplectic algebra

\vskip 2mm

\begin{dfn}\label{d4.4.2.8}
$$
\begin{array}{ccl}
sp(\F)&:=&\left\{\A\in gl(\F)\;|\; \A~\mbox{satisfies (\ref{4.4.2.0})}\right.\\
~&~&~\left.\mbox{and}~A_1\in {\cal M}_{2n\times 2n},~n\in N\right\},
\end{array}
$$
 is called the bundled symplectic algebra.
\begin{align}\label{4.4.2.0}
\left<J\right>\ltimes \A\lplus \A^T\ltimes  \left<J\right>=0,
\end{align}
where
$$
J=\begin{bmatrix}0&1\\-1&0\end{bmatrix}.
$$
\end{dfn}

\end{itemize}

\begin{dfn}\label{d4.4.2.9} A Lie sub-algebra ${\cal J}\subset {\cal G}$ is called an idea, if
\begin{align}\label{4.4.2.1}
\left[g, {\cal J}\right]\in {\cal J}.
\end{align}
\end{dfn}

\begin{exa}\label{e4.4.2.10}
$sl(\F)$ is an idea of $gl(\F)$. Because
$$
\Tr[g,h]=\Tr(g\ltimes h\lminus h\ltimes g)=0,\quad  \forall g\in \gl(\F), \forall h\in sl(\F).
$$
Hence, $\left[\gl(\F), sl(\F)\right]\subset sl(\F)$.
\end{exa}

Many properties of the sub-algebra of $\gl(n,\F)$ can be extended to the sub-algebra of $\gl(\F)$. The following is an example.

\begin{prp}\label{p4.4.2.11}
\begin{align}\label{4.4.2.2}
\gl(\F)=sl(\F)\lplus r\left<1\right>,\quad r\in \F.
\end{align}
\end{prp}

\noindent\textit{Proof.} It  is obvious that
$$
\A=\left(\A\lminus \Tr(\A)\right)\lplus \Tr(\A)\left<1\right>.
$$
Since $\Tr\left(\A\lminus \Tr(\A)\right)=0$, which means
$$
\left(\A\lminus \Tr(\A)\right)\in sl(\F).
$$
The conclusion follows.
\hfill $\Box$

\begin{exa}\label{e4.4.2.12}
\begin{enumerate}
\item Denote by
\begin{align}\label{4.4.2.3}
gl(\A,\F):=\left\{\left<X\right> \;\big|\;\left<X\right>\A\lplus \A\left<X\right>=0\right\}.
\end{align}
It is obvious that $\gl(\A,\F)$ is a vector sub-space of $\gl(\F)$.
Let $\left<X\right>,~\left<Y\right>\in \gl(\A,\F)$. Then
$$
\begin{array}{l}
\left[\left<X\right>,\left<Y\right>\right]\A\lplus \A \left[\left<X\right>,\left<Y\right>\right]^T\\
=\left<X\right>\left<Y\right>\A\lminus \left<Y\right>\left<X\right>\A\\
~+\A \left<Y^T\right>\left<X^T\right>\lminus \left<X^T\right>\left<Y^T\right>\A\\
=\left<X\right>\left<Y\right>\A\lminus \left<Y\right>\left<X\right>\A\\
~+\A \left<Y\right>\left<X\right>\A \lminus \left<X\right>\left<Y\right>\A\\
=0.
\end{array}
$$
Hence, $\gl(\A,\F)\subset \gl(\F)$ is a bundled Lie sub-algebra.

\item Assume $\A$ and $\B$ are congruent. That is, there exists a non-singular $\left<P\right>$ such that $\A=\left<P^T\right>\B\left<P\right>$. Then it is easy to verify that $\pi:~\gl(\A,\F)\ra ~\gl(\B,\F)$ is  an isomorphism, where
    $$\pi(\left<X\right>)=\left<P^{-T}\right>\left<X\right> \left<P\right>.
$$
\end{enumerate}
\end{exa}

\subsection{Further Properties of $gl(\F)$}

\begin{dfn}\label{d4.5.1}Let $p\in \N$.
\begin{enumerate}
\item The $p$-truncated equivalent class is defined as
\begin{align}\label{4.5.1}
\A^{[\cdot,p]}:=\left\{A_i\in \A\;\big|\;i|p\right\}.
\end{align}
\item The $p$-truncated square matrices is defined as
\begin{align}\label{4.5.2}
{\cal M}^{[\cdot,p]}:=\left\{A\in {\cal M}_{n\times n}\;\big|\; n|p\right\}.
\end{align}
\item The $p$-truncated equivalence space is defined as
\begin{align}\label{4.5.3}
\Sigma^{[\cdot,p]}:= {\cal M}^{[\cdot,p]}/\sim.
\end{align}
\end{enumerate}
\end{dfn}

\begin{rem}\label{r4.5.2}
\begin{enumerate}
\item If $A_1\in \A$ is irreducible, $A_1\in {\cal M}_{n\times n}$ and $n$ is not a divisor of $p$,  then $\A^{[\cdot,p]}=\emptyset$.
\item It is obvious that $\left({\cal M}^{[\cdot,p]},\PR, \Sigma^{[\cdot,p]}\right)$ is a discrete bundle, which is a sub-bundle of $\left({\cal M},\PR,\Sigma\right)$. That is, the following (\ref{4.5.4}) is commutative:
 \begin{align}\label{4.5.4}
\begin{array}{ccc}
{\cal M}^{[\cdot,k]}&\xrightarrow{~~~\pi~~~}&{\cal M}\\
\!\!\!\!\!\!\!\!\!\!\!\!\!\!\!\!\PR\downarrow&~&\!\!\!\!\!\!\!\!\!\!\PR\downarrow\\
\Sigma^{[\cdot,k]}&\xrightarrow{~~~\pi'~~~}&\Sigma
\end{array}
\end{align}
Where $\pi$ and $\pi'$ are including mappings.

\item It is ready to verify that  $\Sigma^{[\cdot,p]}$ is closed with $\lplus$ and $[\cdot,\cdot]$, defined in (\ref{4.3.4}), hence the including mapping $\pi': \Sigma^{[\cdot,p]}\ra \Sigma$ is a Lie algebra homomorphism. Identifying $\A^{[\cdot,p]}$ with its image $\A=\pi'\left(\A^{[\cdot,p]}\right)$, then $\Sigma^{[\cdot,p]}$ becomes a Lie sub-algebra of $gl(\F)$. We denote this as
  \begin{align}\label{4.5.5}
\gl^{[\cdot,p]}(\F):=\left(\Sigma^{[\cdot,p]},~\lplus,~[\cdot,\cdot]\right),
\end{align}
and call $\gl^{[\cdot,p]}(\F)$ the $p$-truncated Lie sub-algebra of $\gl(\F)$.

\item Let $\Gamma \subset \gl(\F)$ be a Lie sub-algebra. Then its $p$-truncated sub-algebra $\Gamma^{[\cdot,p]}$ is defined in a similar way as for $\gl^{[\cdot,p]}$. Alternatively, it can be considered as
\begin{align}\label{4.5.6}
\Gamma^{[\cdot,p]}:=\Gamma \bigcap \Sigma^{[\cdot,p]}.
\end{align}
\end{enumerate}
\end{rem}

\begin{dfn}[\cite{hum72}]\label{d4.5.3} Let $g$ be a Lie algebra.
\begin{enumerate}
\item  Denote the derived serious as ${\cal D}(g):=[g,g]$, and
$$
{\cal D}^{(k+1)}(g):={\cal D}\left({\cal D}^k(g)\right),\quad k=1,2,\cdots.
$$
Then $g$ is solvable, if there exists an $n\in \N$ such that ${\cal D}^{(n)}(g)=\{0\}$.

\item  Denote the descending central series as ${\cal C}(g):=[g,g]$, and
$$
{\cal C}^{(k+1)}(g):=\left[g,{\cal C}^{(k)}(g)\right],\quad k=1,2,\cdots.
$$
The $g$ is nilpotent, if there exists an $n\in \N$ such that ${\cal C}^{(n)}(g)=\{0\}$.
\end{enumerate}
\end{dfn}

\begin{dfn}\label{d4.5.4} Let $\Gamma\subset \gl(\F)$ be a sub-algebra of $\gl(\F)$.

\begin{enumerate}
\item $\Gamma$ is solvable, if for any $p\in \N$, the truncated sub-algebra $
\Gamma^{[\cdot, p]}$
is solvable.
\item $\Gamma$ is nilpotent, if for any $p\in \N$, the truncated sub-algebra
$\Gamma^{[\cdot, p]}$ is nilpotent.
\end{enumerate}
\end{dfn}

\begin{dfn}[\cite{hum72}]\label{d4.5.5} Let $g$ be a Lie algebra.
\begin{enumerate}
\item $g$ is simple if it has no non-trivial idea, that is, the only ideas are $\{0\}$ and $g$ itself.

\item $g$ is semi-simple if it has no solvable idea except $\{0\}$.
\end{enumerate}
\end{dfn}

Though the following proposition is simple, it is also fundamental.

\begin{prp}\label{p4.5.6} The Lie algebra $\gl^{[\cdot,p]}(\F)$ is isomorphic to the classical linear algebra $\gl(p,\F)$.
\end{prp}

\noindent\textit{Proof}. First we construct a mapping $\pi: gl^{[\cdot,p]}(\F)\ra \gl(p,\F)$ as follows: Assume $\A\in gl^{[\cdot,p]}(\F)$ and $A_1\in \A$ is irreducible. Say, $A_1\in {\cal M}_{n\times n}$, then by definition, $n|p$. Denote by $s=p/n$, then define
$$
\pi(\A):=A_1\otimes I_s\in \gl(p,\F).
$$
Set $\pi': \gl(p,\F)\ra \gl^{[\cdot,p]}(\F)$ as
$$
\pi'(B):=\B\in \gl^{[\cdot,p]}(\F),
$$
then it is ready to verify that $\pi$ is a bijective mapping and $\pi^{-1}=\pi'$.

By the definitions of $\lplus$ and $[\cdot,\cdot]$, it is obvious that $\pi$ is a Lie algebra isomorphism.
\hfill $\Box$

The following properties are available for classical $gl(n,\F)$ \cite{hum72}, \cite{wan13}.  Using Proposition \ref{p4.5.6}, it is easy to verify that they are also available for $gl(\F)$.

\begin{prp}\label{p4.5.601} Let $g\subset \gl(\F)$ be a Lie sub-algebra.
\begin{enumerate}
\item If $g$ is nilpotent then it is solvable.
\item If $g$ is solvable (or nilpotent) then so is its sub-algebra, its homomorphic image.
\item If $h\subset g$ is an idea of $g$ and $h$ and $g/h$ are solvable, then $g$ is also solvable.
\end{enumerate}
\end{prp}

\begin{dfn}\label{d4.5.7} Let $\A\in \gl(\F)$. The adjoint representation $\ad_{\A}:~\gl(\F)\ra \gl(\F)$ is defined as
\begin{align}\label{4.5.7}
\ad_{\A}\B=[\A,\B].
\end{align}
\end{dfn}

To see (\ref{4.5.7}) is well defined, we have to prove that
\begin{align}\label{4.5.8}
\ad_{\A}\B=\left<\ad_AB\right>,\quad A\in\A,~\B\in \B.
\end{align}
It follows from the consistence of $\ltimes$ and $\lplus$ ($\lminus$) with $\sim_{\ell}$ immediately.

\begin{exa}\label{e4.5.701} Consider $\A\in \Sigma$. Assume $\A$ is nilpotent, that is, there is a $k>0$ such that $\A^k=0$. Then $\ad_{\A}$ is also nilpotent.

Note that $A^k=0$ if and only if $(A\otimes I_s)^k=0$. Similarly, $\ad_A^k=0$ if and only if $\ad^k_{A\otimes I_s}=0$. Hence, we need only to show that $\ad_{A}$ is nilpotent, where $A\in \A$ and $A\in {\cal M}_{n\times n}$. Using the definition that
$$
\ad_AB=AB-BA,
$$
a straightforward computation shows that
$$
\ad^m_A B=\dsum_{i=0}^m(-1)^i\binom{m}{i}A^{m-i}BA^i,\quad \forall B\in {\cal M}_{n\times n}.
$$
As $m=2k-1$, it is clear that $\ad^m_AB=0$, $\forall B\in {\cal M}_{n\times n}$. It follows that
$$
\ad^{2k-1}_A=0.
$$
\end{exa}

\begin{dfn}\label{d4.5.8}
\begin{enumerate}
\item Let $A,~B\in \gl(n,\F)$. Then the Killing form $(\cdot,\cdot)_K: \gl(n,\F)\times \gl(n,\F)\ra \F$ is defined as (We refer to \cite{hum72} for original definition. The following definition is with a mild modification.)
\begin{align}\label{4.5.801}
(A,B)_K:=\Tr(\ad_A \ad_B).
\end{align}

\item Assume $\A,\B\in \gl(\F)$. The killing form $(\cdot,\cdot)_K:\gl(\F)\times \gl(\F)\ra \F$ is defined as
\begin{align}\label{4.5.9}
(\A,\B)_K:=\Tr\left(\ad_{\A}\ltimes \ad_{\B}\right).
\end{align}
\end{enumerate}
\end{dfn}

To see the killing form is well defined, we also need to prove
\begin{align}\label{4.5.10}
(\A,\B)_K=(A,B)_K,\quad A\in\A,~B\in \B.
\end{align}
Similar to (\ref{4.5.8}), it can be verified by a straightforward calculation.

Because of the equations (\ref{4.5.8}) and (\ref{4.5.10}), the following properties of finite dimensional Lie algebras \cite{wan13} can easily be extended to $\gl(\F)$:

\begin{prp}\label{p4.5.9} Consider $g=\gl(\F)$. Let $\A$, $\A_1$, $\A_2$, $\B$, $\E\in g$, $c_1,~c_2\in \F$. Then
\begin{enumerate}
\item
\begin{align}\label{4.5.11}
(\A,\B)_K=(\B,\A)_K.
\end{align}
\item
\begin{align}\label{4.5.12}
\begin{array}{l}
\left(c_1\A_1\lplus c_2\A_2,\B\right)_K\\
=c_1\left(\A_1,\B\right)_K\lplus c_2\left(\A_2,B\right)_K,\\
~~~~~~~~~~c_1,\;c_2\in \F.
\end{array}
\end{align}
\item
\begin{align}\label{4.5.13}
\left(\ad_{\A}\B,~\E\right)_K\lplus  \left(\B,~\ad_{\A}\E\right)_K=0.
\end{align}
\item Let $h\subset g$ be an idea of $g$, and $\A,~\B\in h$. Then
\begin{align}\label{4.5.14}
\left(\A,~\B\right)_K =\left(\A,~\B\right)_K^h.
\end{align}
The right hand side means the Killing form on the ideal $h$.
\item
A sub-algebra $\xi\subset g$ is semi-simple, if and only if, its Killing form is non-degenerated.
\end{enumerate}
\end{prp}

The Engel theorem can easily be extended to $gl(\F)$:

\begin{thm}\label{t4.5.10} Let $\{0\}\neq g\subset \gl(\F)$ be a bundled Lie sub-algebra. Assume each $\A\in g$ is nilpotent, (i.e., for each $\A \in g$ there exists a $k>0$ such that $\A^k=\left<A^k\right>=0$).
\begin{enumerate}
\item  If $g$ is finitely generated, then there exists a vector $X\neq 0$ (of suitable dimension) such that
$$
G\ltimes X=0,\quad \forall G\in g.
$$
\item $g$ is nilpotent.
\end{enumerate}
\end{thm}

\begin{dfn}[\cite{hum72}]\label{d4.5.11} Let $V$ be an $n$-dimensional vector space. A flag in $V$ is a chain of subspaces
$$
0=V_0\subset V_1\subset V_2\subset \cdots\subset V_n=V,
$$
with $\dim(V_i)=i$. Let $A\in End(V)$ be an endomorphism of $V$. $A$ is said to stabilize this flag if
$$
AV_i\subset V_i,\quad i=1,\cdots,n.
$$
\end{dfn}

Lie theorem can be extended to $\gl(\F)$ as follows.

\begin{thm}\label{t4.5.12} Assume $g\subset gl(\F)$ is a solvable Lie sub-algebra. Then for any $p>0$ there is a flag of ideals
$0={\cal I}_0\subset {\cal I}_1\subset \cdots\subset {\cal I}_p$, such that the truncated $g^{[\cdot, p]}$ stabilizes the flag.
\end{thm}

\begin{cor}\label{c4.5.13} Assume $g\subset gl(\F)$ is a Lie sub-algebra. $g$ is solvable if and only if, $D(g)$ is nilpotent.
\end{cor}

\begin{exa}\label{e4.5.14} Consider the bundled Lie sub-algebras $t(\F)$ and $n(\F)$. It is easy to verify the following:
\begin{enumerate}
\item $t(\F)$ is solvable;
\item $n(\F)$ is nilpotent.
\end{enumerate}
\end{exa}

Even though $gl(\F)$ is an infinite dimensional Lie algebra, it has almost all the properties of finite dimensional Lie algebras. This claim can be verified one by one easily. The reason is: $\gl(\F)$ is essentially a union of finite dimensional Lie algebras.

\section{Lie Group on Nonsingular M-equivalence Space}

\subsection{Bundled Lie Group}

Consider $\Sigma:=\Sigma_1$, we define a subset
\begin{align}\label{5.7.1.1}
GL(\F):=\left\{\A\in \Sigma\;|\; \Dt(\A)\neq 0\right\}.
\end{align}

We emphasize the fact that $GL(\F)$ is an open subset of $\Sigma$. For an open subset of bundled manifolds we have the following result.

\begin{prp}\label{p5.7.1.1} Let $M$ be a bundled manifold, and $N$ an open subset of $M$. Then $N$ is also a bundled manifold.
\end{prp}

\noindent{Proof}. It is enough to construct an open cover of $N$. Starting from the open cover of $M$, which is denoted as
$${\cal C}=\{U_{\lambda}\;|\; \lambda\in \Lambda\},$$
we construct
$${\cal C}_N:=\{U_{\lambda}\cap N\;|\; U_{\lambda}\in {\cal C},~\lambda\in \Lambda\}$$
Then we can prove that it is $C^r$ ($C^{\infty}$ or $C^{\omega}$) comparable as long as ${\cal C}$ is.
Verifying other conditions is trivial.
\hfill $\Box$

\begin{cor}\label{c5.7.1.2} $GL(\F)$ is a bundled manifold.
\end{cor}

Note that, Proposition \ref{p2.4.5} shows that $\ltimes:GL(\F)\times GL(\F)\ra GL(\F)$ is well define.

\begin{dfn}\label{d5.7.1.3} A topological space $G$ is a bundled Lie group, if
\begin{enumerate}
\item it is a bundled analytic manifold;
\item it is a group;
\item the product $A\times B\ra AB$  and the inverse mapping $A\ra A^{-1}$ are analytic.
\end{enumerate}
\end{dfn}

The following result is an immediate consequence of the definition.

\begin{thm}\label{t5.7.1.4}
 $GL(\F)$ is a bundled Lie group.
\end{thm}

\noindent\textit{Proof}.
 We already known that $GL(\F)$ is a bundled analytic manifold. We first prove $GL(\F)$ is a group.
It is ready to verify that $\left<1\right>$ is the identity. Moreover, $\A^{-1}=\left<A^{-1}\right>$. The conclusion follows.

Using a simple coordinate chart, it is obvious that the inverse and product are two analytic mappings.
\hfill $\Box$

\subsection{Relationship with $gl(\F)$}

Denote by
$$
{\cal W}:=\left\{A\in {\cal M}_1\;|\; \det(A)\neq 0\right\},
$$
and
$$
{\cal W}_s:={\cal W}\cap {\cal M}_{s\times s},\quad s=1,2,\cdots.
$$

Consider the bundle: $\left({\cal M}_1, \PR, \Sigma\right)$, where the map $\PR:A\mapsto \A$ is the national projection. It has a natural sub-bundle: $\left({\cal W},\PR, GL(\F)\right)$  via the following bundle morphism as
\begin{align}\label{5.7.2.1}
\begin{array}{ccc}
{\cal W}&\xrightarrow{~~~\pi~~~}&{\cal M}_{1}\\
\!\!\!\!\!\!\!\!\!\PR\downarrow&~&\!\!\!\!\!\!\!\!\!\!\PR\downarrow\\
GL(\F)&\xrightarrow{~~~\pi'~~~}&\Sigma
\end{array}
\end{align}

In fact, the projection leads to Lie group homomorphism.

\begin{thm}\label{t5.7.2.0}
\begin{enumerate}
\item With natural group and differential structures, ${\cal W}_s=GL(s,\F)$ is a Lie group.

\item  Consider the projection $\PR$. Restrict it to each leaf yields
\begin{align}\label{5.7.2.101}
\PR\big|_{{\cal W}_s}: GL(s,\F)\ra GL(\F).
\end{align}
Then $\PR\big|_{{\cal W}_s}$ is a Lie group homomorphism.
\item
Set the image set as $\PR(GL(s,\F)):=\Psi_s$. Then $\Psi_s < GL(\F)$ is a Lie sub-group. Moreover,
\begin{align}\label{5.7.2.102}
\PR\big|_{{\cal W}_s}:GL(s,\F)\ra \Psi_s
\end{align}
is a Lie group isomorphism.
\end{enumerate}
\end{thm}
%
%

\begin{dfn}\label{d5.7.2.2} A vector field $\left<\xi\right>\in V\left(GL(\F)\right)$, (where for each $P\in GL(\F)$, $\left<\xi(P)\right>\in T_P\left(GL(\F)\right)$), is called a left-invariant vector field, if for any $\A\in GL(\F)$
$$
\begin{array}{ccl}
\left(L_{\A}\right)_*\left(\left<\xi\right>(P)\right)&:=&\left<\left(L_A\right)_*(\xi(P))\right>\\
~&=&\left<\xi(AP)\right>=\left<\xi\right>(AP).
\end{array}
$$
\end{dfn}

Then it is easy to verify the following relationship between $GL(\F)$ and $gl(\F)$:

\begin{thm}\label{t5.7.2.3} The corresponding Lie algebra of the bundled Lie group $GL(\F)$ is $gl(\F)$ in the following natural sense:
$$
gl(\F)\simeq T_{\left<1\right>}\left(GL(\F)\right)\xrightarrow{\left(L_{\A}\right)_*}  T_{\A}\left(GL(\F)\right).
$$
That is, $gl(\F)$ is a Lie algebra isomorphic to the Lie algebra consists of the vectors on the tangent space of $GL(\F)$ at identity. Then these vectors generate the left-invariant vector fields which form the tangent space at any $\A\in GL(\F)$.
\end{thm}

Let $\PR:{\cal M}_{n\times n}\ra \Sigma$ be the natural mapping $A\mapsto \A$. Then we have the following commutative picture:
\begin{align}\label{5.7.2.3}
\begin{array}{ccc}
gl(n,\F)&\xrightarrow{~~~\exp~~~}&GL(n,\F)\\
\!\!\!\!\!\!\PR\downarrow&~&\!\!\!\!\!\!\PR\downarrow\\
gl(\F)&\xrightarrow{~~~\exp~~~}&GL(\F)
\end{array}
\end{align}
where $n=1,2,\cdots$. Recall (\ref{2.2.15}), we know that the exponential mapping $\exp$ is well defined $\forall \left<X\right>\in gl(\F)$.

Graph (\ref{5.7.2.3}) also shows the relationship between $gl(\F)$ and $GL(\F)$, which is a generalization of the relationship between $gl(n,\F)$ and $GL(n,\F)$.

\subsection{Lie Subgroups of $GL(\F)$}

It has been discussed that $gl(\F)$ has some useful Lie sub-algebras. It is obvious that $GL(\F)$ has some Lie sub-groups, corresponding to those sub-algebras of $gl(\F)$. They are briefly discussed as follows.

\begin{itemize}

\vskip 2mm
\item

Bundled orthogonal Lie sub-group

\vskip 2mm

\begin{dfn}\label{d5.7.3.1} $\A\in GL(\F)$ is said to be orthogonal, if $A^T=A^{-1}$.
\end{dfn}
It is also easy to verify the following:

\begin{prp}\label{p5.7.3.2} Assume $\A$ and $\B$ are orthogonal, then so is $\A\ltimes \B$.
\end{prp}

We, therefore, can define the following bundled Lie sub-group of $GL(\F)$ as follows.

\begin{dfn}\label{d5.7.3.3}
$$
O(\F):=\left\{\A\in GL(\F)\;\big|\;  \A^T=\A^{-1}\right\}
$$
is called the bundled orthogonal group.
\end{dfn}

It is easy to verify the following proposition:

\begin{prp}\label{p5.7.3.4} Consider the bundled orthogonal group.

\begin{enumerate}

\item

$O(\F)$ is a Lie sub-group of $GL(\F)$, i.e., $O(\F)< GL(\F)$.

\item
$$SO(\F):=\left\{\A\in O(\F)\;\big|\; \det(\A)=1\right\}.
$$
Then $SO(\F)<O(\F)<GL(\F)$.

\item The Lie algebra for both $O(\F)$ and $SO(\F)$  is $o(\F)$.
\end{enumerate}
\end{prp}

\vskip 2mm

\item Bundled special linear group

\vskip 2mm

\begin{dfn}\label{d5.7.3.5}
$$
SL(\F):=\left\{\A\in GL(\F)\;\big|\;  \det(\A)=1\right\}
$$
is called the bundled special linear group.
\end{dfn}

Similar to the case of orthogonal algebra, it is easy to verify the following:

\begin{prp}\label{t5.7.3.6}  Consider the bundled special linear group.

\begin{enumerate}
\item $SL(\F)$ is a Lie sub-group of $GL(\F)$, i.e., $SL(\F)<GL(\F)$.
\item The Lie algebra of $SL(\F)$ is $sl(\F)$.
\end{enumerate}
\end{prp}

\vskip 2mm
\item Bundled upper triangular group
\vskip 2mm

\begin{dfn}\label{d5.7.3.7}
$$
T(\F):=\left\{\A\in GL(\F)\;\big|\;  \A~\mbox{is upper triangular}\right\}
$$
is called the bundled upper triangular group.
\end{dfn}

\begin{prp}\label{t5.7.3.8} Consider the bundled upper triangular group.

\begin{enumerate}
\item $T(\F)$ is a Lie sub-group of $GL(\F)$, i.e., $T(\F)<GL(\F)$.
\item The Lie algebra of $T(\F)$ is $t(\F)$.
\end{enumerate}
\end{prp}

\vskip 2mm

\item Bundled special upper triangular group

\vskip 2mm

\begin{dfn}\label{d5.7.3.9}
$$
N(\F):=\left\{\A\in T(\F)\;\big|\;  \det(\A)=1\right\}
$$
is called the bundled special upper triangular group.
\end{dfn}

\begin{prp}\label{t5.7.3.10} Consider the bundled special upper triangular group.

\begin{enumerate}
\item $N(\F)$ is a Lie sub-group of $T(\F)$, i.e., $N(\F)<T(\F)<GL(\F)$.
\item The Lie algebra of $N(\F)$ is $n(\F)$.
\end{enumerate}
\end{prp}

\vskip 2mm

\item Bundled symplectic group

\vskip 2mm

\begin{dfn}\label{d5.7.3.11}
$$
\begin{array}{ccl}
SP(\F)&:=&\left\{\A\in GL(\F)\;\big|\; A_1\in {\cal M}_{2n}\right.\\
~&~&~~~~~\left.\mbox{satisfies (\ref{5.7.3.01})},\; n\in \N \right\},
\end{array}
$$
is called the bundled symplectic group.
\end{dfn}
\begin{align}\label{5.7.3.01}
\A^T\left<J\right>\A=\left<J\right>,
\end{align}
where $J$ is defined in (\ref{2.3.8}).
\end{itemize}

\begin{prp}\label{t5.7.3.12} Consider the bundled symplectic group.

\begin{enumerate}
\item $SP(\F)$ is a Lie sub-group of $GL(\F)$, i.e., $SP(\F)<GL(\F)$.
\item The Lie algebra of $SP(\F)$ is $sp(\F)$.
\end{enumerate}
\end{prp}

\subsection{Symmetric Group}

Let ${\bf S}_k$ be the $k$-th order symmetric group. Denote
$$
{\bf S}:=\bigcup_{k=1}^{\infty}{\bf S}_k.
$$

\begin{dfn}\label{d5.7.3.1}
A matrix $A\in {\cal M}_k$ is called a permutation matrix, if
$\Col(A)\subset \D_k$ and $\Col(A^T)\subset \D_k$. The set of $k\times k$ permutation matrices is
denoted by ${\cal P}_k$.
\end{dfn}

\begin{prp}\label{p5.7.3.2} Consider the set of permutation matrix.
\begin{enumerate}
\item If $P\in {\cal P}_k$, then
\begin{align}\label{5.7.3.1}
P^T=P^{-1}.
\end{align}
\item
$$
{\cal P}_k< O(k,\R)< GL(k,\R).
$$
\end{enumerate}
\end{prp}

Let $\sigma\in {\bf S}_k$. Define a permutation matrix
$M_{\sigma}\in {\cal P}_k$ as $M_{\sigma}:=\left[m_{i,j}\right]$, where
\begin{align}\label{5.7.3.2}
m_{i,j}=\begin{cases}
1,& \sigma(j)=i\\
0,& \mbox{otherwise}.
\end{cases}
\end{align}

The following proposition is easily verifiable.

\begin{prp}\label{p5.7.3.3}
Define $\pi: {\bf S}_k\ra {\cal P}_k$, where $\pi(\sigma):=M_{\sigma}\in {\cal P}_k$ is constructed by (\ref{5.7.3.2}). Then $\pi$ is an isomorphism.
\end{prp}

Assume $\sigma,~\lambda\in {\bf S}_k$. Then Proposition \ref{p5.7.3.3} leads to
\begin{align}\label{5.7.3.3}
M_{\sigma\circ \lambda}=M_{\sigma}M_{\lambda}.
\end{align}

Next, assume $\sigma\in {\bf S}_m$, $\lambda\in {\bf S}_n$, we try to generalize (\ref{5.7.3.3}).

\begin{dfn}\label{d5.7.3.4} Assume $\sigma\in {\bf S}_m$, $\lambda\in {\bf S}_n$. The (left) STP of $\sigma$ and $\lambda$ is defined by
\begin{align}\label{5.7.3.4}
M_{\sigma\ltimes \lambda}=M_{\sigma}\ltimes M_{\lambda}\in {\cal P}_t,
\end{align}
where $t=m\vee n$. That is,
\begin{align}\label{5.7.3.5}
\sigma\ltimes \lambda:=\pi^{-1}\left(M_{\sigma}\ltimes M_{\lambda}\right)\in {\bf S}_t.
\end{align}
\end{dfn}

Similarly, we can define the right STP of $\sigma$ and $\lambda$.

Now, it is clear that $(S,\ltimes)<({\cal M}_1,\ltimes)$ is a sub-monoid. To get a bundled Lie
subgroup structure, we consider the quotient space
$$
{\cal D}:=\left({\bf S}, \ltimes\right)/\sim_{\ell}.
$$

Then we have the following:

\begin{thm}\label{t5.7.3.5} ${\cal D}$ is a discrete bundled sub-Lie group of $GL(\F)$.
\end{thm}

${\cal D}$ may be used to investigate the permutation of uncertain number of elements.


\section{V-equivalence}

This section considers the vector equivalence (V-equivalence) of vectors. Since many results are parallel and simpler than M-equivalence, some detailed discussions are omitted.

\subsection{Equivalence of Vectors of Different Dimensions}

Consider the set of vectors on field $\F$. Denote it as
\begin{align}\label{6.2.1}
{\cal V}:=\bigcup_{i=1}^{\infty}{\cal V}_i,
\end{align}
where ${\cal V}_i$ is the $i$-dimensional vector space, which is a subset of ${\cal V}$.

Our purpose is to build a vector space structure on ${\cal V}$, precisely speaking, on equivalence classes of ${\cal V}$.  To this end, we first propose an equivalence relation.

\begin{dfn}\label{d6.2.1}
\begin{enumerate}
\item Let $X,~Y\in {\cal V}$. $X$ and $Y$ are said to be V-equivalent, denoted by $X\lra Y$, if there exist two one-vectors ${\bf 1}_s$ and ${\bf 1}_t$ such that
\begin{align}\label{6.2.2}
X\otimes {\bf 1}_s=Y\otimes {\bf 1}_t.
\end{align}
\item The equivalent class is denoted as
$$
\lceil X \rceil:=\left\{Y\;\big|\;Y\lra X\right\}.
$$
\item In an equivalent class $\lceil X\rceil$ a partial order ($\leqslant$) is defined as: $X\leqslant Y$, if there exists a one-vector ${\bf 1}_s$ such that $X\otimes {\bf 1}_s=Y$. $X_1\in \lceil X\rceil$ is irreducible, if there are no $Y$ and ${\bf 1}_s$, $s>1$, such that $X_1=Y\otimes {\bf 1}_s$.
\end{enumerate}
\end{dfn}

\begin{rem}\label{r6.2.2}
\begin{enumerate}
\item The equivalence defined above can be seen as the left equivalence. Formally, we set $\lra:=\lra_l$.
\item The right equivalence can be defined as follows: Let $X,~Y\in {\cal V}$. $X$ and $Y$ are said to be right equivalent, denoted by $X\lra_r Y$, if there exist two one-vectors ${\bf 1}_s$ and ${\bf 1}_t$ such that
\begin{align}\label{6.2.3}
{\bf 1}_s\otimes X={\bf 1}_t\otimes Y.
\end{align}
\item The right equivalent class of $X$ is denoted as $\lceil X\rceil_r$. Of course, we have $\lceil X\rceil_l=\lceil X\rceil$.
\item Hereafter, the left equivalence is considered as the defaulted one. That is, we always assume $\lra=\lra_l$. But with some obvious modifications one sees easily that the arguments/results in the sequel are also valid for right equivalence.
\end{enumerate}
\end{rem}

The following properties of the M-equivalence are also true for V-equivalence. The proofs are also similar to the M-equivalence. Therefore, they are omitted.

\begin{thm}\label{t6.2.3}
\begin{enumerate}
\item If $X\lra Y$, then there exists a vector $\Gamma$ such that
\begin{align}\label{6.2.4}
X=\Gamma\otimes \I_{\b},\quad Y=\Gamma\otimes \I_{\a}.
\end{align}
\item In each class $\lceil X\rceil $ there exists a unique $X_1\in \lceil X\rceil $, such that $X_1$ is  irreducible. \end{enumerate}
\end{thm}

\begin{rem}\label{r6.2.4}
\begin{enumerate}
\item  If $X=Y\otimes {\bf 1}_s$, then $Y$ is called a divisor of $X$ and $X$ is called a multiple of $Y$. This relation determined the order $
Y\leqslant X$.
\item  If (\ref{6.2.4}) holds and $\a,~\b$ are co-prime, then the $\Gamma$ satisfying (\ref{6.2.4}) is called the greatest common divisor of $X$ and $Y$. Moreover, $\Gamma=gcd(X,Y)$ is unique.
\item  If (\ref{6.2.2}) holds and $s,~t$ are co-prime, then
\begin{align}\label{6.2.5}
\Xi:=X\otimes \I_{s}=Y\otimes \I_{t}
\end{align}
is called the least common multiple  of $X$ and $Y$. Moreover, $\Xi=\lcm(X,Y)$ is unique.

\item Consider an equivalent class $\lceil X\rceil$, denote the unique irreducible element by $X_1$, which is called the root element. All the elements in $\lceil X\rceil$ can be expressed as
 \begin{align}\label{6.2.6}
 X_i=X_1\otimes \I_i,\quad i=1,2,\cdots.
 \end{align}
$X_i$ is called the $i$-th element of $\lceil X\rceil$. Hence, an equivalent class $\lceil X\rceil$ is a well ordered sequence as:
$$
\lceil X\rceil=\left\{X_1,~X_2,~X_3,\cdots\right\}.
$$
\end{enumerate}
\end{rem}

We also have a lattice structure on $\lceil X\rceil =\{X_1,X_2,\cdots\}$:

\begin{prp}\label{p6.2.5} $(\lceil X\rceil,~\leq)$ is a lattice.
\end{prp}

\noindent\textit{Proof.} It is easy to verify that for $U,~V\in \X$,
$$
\sup(U,~V)=lcm(U,~V);\quad \inf(U,~V)=gcd(U,~V).
$$
The conclusion follows.
\hfill $\Box$

\begin{prp}\label{p6.2.6} Let $A_1\in {\cal M}$ and $X_1\in {\cal V}$ be both irreducible. Then the two lattices $\A$ and $\X$, generated by $A_1$ and $X_1$ are isomorphic. Precisely,
 \begin{align}\label{6.2.7}
\left(\A,\prec\right)\approxeq \left(\lceil X\rceil,\leqslant \right).
\end{align}
The isomorphism is: $\varphi:A_s=A_1\otimes I_s\mapsto X_s=X_1\otimes \I_s$.
\end{prp}
%
%
%

Next we investigate the lattice structure on ${\cal V}$. Consider ${\cal V}_i$ and ${\cal V}_j$, where $i|j$ and $k=j/i$. Then we have
${\cal V}_i\otimes {\bf 1}_k\subset {\cal V}_j$. This order is denoted by
$$
{\cal V}_i\sqsubseteq {\cal V}_j.
$$
Using this order, the following result is obvious.

\begin{prp}\label{p6.2.801}
$({\cal V}, \sqsubseteq)$ is a lattice with
\begin{align}\label{6.2.8}
\begin{array}{l}
\sup ({\cal V}_i,{\cal V}_j)={\cal V}_{i\vee j},\\
\inf ({\cal V}_i,{\cal V}_j)={\cal V}_{i\wedge j}.
\end{array}
\end{align}
\end{prp}

\begin{prp}\label{p6.2.802} Assume $\lceil X \rceil$ has its irreducible vector $X_1\in \R^i$. Define $\varphi: \lceil X \rceil\ra {\cal V}$ as $\varphi (X_1):={\cal V}_i$.
Then $\varphi: (\lceil X \rceil,\leq ) \ra ({\cal V},\sqsubseteq)$ is a lattice homomorphism.
\begin{align}\label{6.2.9}
 (\lceil X \rceil,\leq ) \approx ({\cal V},\sqsubseteq).
\end{align}
\end{prp}

We also have the following isomorphic relation:

\begin{prp}\label{p6.2.8} Define $\varphi: {\cal V}\ra {\cal M}_{\mu}$ as:
\begin{align}\label{6.2.10}
\varphi\left({\cal V}_i\right):= {\cal M}_{\mu}^i.
\end{align}
Then $\varphi$ is a lattice isomorphism.
\end{prp}

\begin{dfn}\label{d6.2.9}
\begin{enumerate}
\item  Let $p\in \N$. The $p$-lower truncated vector space is defined as
\begin{align}\label{6.2.11}
{\cal V}^{[p,\cdot]}:=\bigcup_{\{s\;\big|\; p|s\}}{\cal V}_s.
\end{align}
\item The quotient space of ${\cal V}/\lra$ is denoted as
\begin{align}\label{6.2.12}
\Omega_{{\cal V}}:=\left\{\lceil X\rceil\;\big|\;X\in {\cal V}\right\}.
\end{align}
\item The subspace
\begin{align}\label{6.2.13}
\Omega_{{\cal V}}^i:={\cal V}^{[i,\cdot]}/\lra
=\left\{\lceil X\rceil \;\big|\; X_1\in {\cal V}_i\right\}.
\end{align}
\end{enumerate}
\end{dfn}

Consider an equivalence $\lceil X\rceil$. Let $X_1\in \lceil X\rceil$ be the irreducible element and $\dim(X_1)=r$. Then we define a mapping $\psi: \lceil X\rceil \ra {\cal V}$ as
\begin{align}\label{6.2.14}
\psi(X_i)={\cal V}_{ir}, \quad i=1,2,\cdots
\end{align}
Similar to matrix case, we have
\begin{prp}\label{p6.2.10}
Let $\psi: \lceil X\rceil \ra {\cal V}$ be defined by (\ref{6.2.14}). Then $\psi: (\lceil X\rceil,\leq)\ra ({\cal V},~\sqsubseteq)$ is a lattice homomorphism.
\end{prp}

\subsection{Vector Space Structure on V-equivalence Space}

To begin with, we define an addition between vectors of different dimensions.

\begin{dfn}\label{d6.3.1} Let $X\in {\cal V}_p$ and $Y\in {\cal V}_q$. $t=p\vee q$. Then
\begin{enumerate}
\item the vector addition $\lvplus: {\cal V}\times {\cal V}\ra {\cal V}$ of $X$ and $Y$ is defined as
\begin{align}\label{6.3.1}
X\lvplus Y:=\left( X\otimes \I_{t/p}\right)+ \left( Y\otimes \I_{t/q}\right);
\end{align}
\item the subtraction $\lvminus$ is defined as
\begin{align}\label{6.3.2}
X\lvminus Y:=X\lvplus (-Y).
\end{align}
\end{enumerate}
\end{dfn}

\begin{rem}\label{r6.3.101}  Let $X\in {\cal V}^T_p$ and $Y\in {\cal V}^T_q$ be two row vectors. Then we define
\begin{enumerate}
\item the vector addition as
\begin{align}\label{6.3.201}
X\lvplus Y:=\left( X\otimes \I^T_{t/p}\right)+ \left( Y\otimes \I^T_{t/q}\right);
\end{align}
\item the subtraction  as
\begin{align}\label{6.3.202}
X\lvminus Y:=X\lvplus (-Y).
\end{align}
\end{enumerate}
\end{rem}

\begin{prp}\label{p6.3.2} The vector addition $\lvplus$ is consistent with the equivalence $\lra$. That is,
if $X\lra \tilde{X}$ and $Y\lra \tilde{Y}$, then $X\lvplus Y\lra \tilde{X}\lvplus\tilde{Y}$.
\end{prp}

\noindent\textit{Proof}. Since $X\lra \tilde{X}$, according to Theorem \ref{t6.2.3}, there exists $\Gamma$, say, $\Gamma\in {\cal V}_p$, such that
$$
X=\Gamma\otimes \I_{\a},\quad \tilde{X}=\Gamma\otimes \I_{\b}.
$$
Similarly, there exists $\Pi$, say, $\Pi\in {\cal V}_q$, such that
$$
Y=\Pi\otimes \I_{s},\quad \tilde{Y}=\Pi\otimes \I_{t}.
$$
Let $\xi=p\vee q$, $\eta=p\a\vee sq$, and $\eta=\xi \ell$. Then
$$
\begin{array}{l}
X\lvplus Y=\left(\Gamma\otimes \I_{\a}\right)\lvplus \left(\Pi\otimes \I_{s}\right)\\
~=\left[\Gamma\otimes \I_{\a}\otimes \I_{\eta/(\a p)}\right]+
\left[\Pi\otimes \I_{s}\otimes \I_{\eta/(s q)}\right]\\
~=\left[\Gamma\otimes \I_{\eta/p}\right]+
\left[\Pi\otimes \I_{\eta/q}\right]\\
~=\left(\left[\Gamma\otimes \I_{\xi/p}\right]+
\left[\Pi\otimes \I_{\xi/q}\right]\right)\otimes \I_{\ell}\\
~=(\Gamma\lvplus\Pi)\otimes \I_{\ell}.
\end{array}
$$
Hence $X\lvplus Y\lra \Gamma \lvplus \Pi$. Similarly, we can show that $\tilde{X}\lvplus\tilde{Y}\lra \Gamma \lvplus \Pi$.
The conclusion follows.
\hfill $\Box$

\begin{cor}\label{c6.3.3} The vector addition $\lvplus$ (or subtraction $\lvminus$) is well defined on the quotient space $\Omega_{{\cal V}}$, as well as $\Omega_{{\cal V}}^i$, $i=1,2,\cdots$. That is,
\begin{align}\label{6.3.3}
\begin{array}{l}
\X \lvplus \Y:=\lceil X \lvplus Y \rceil,\\
~~~~X,~Y\in {\cal V}, \left(\mbox{or}~ \X,~\Y\in \Omega_{{\cal V}}\right).
\end{array}
\end{align}
\end{cor}

Let $\X\in \Omega_{{\cal V}}$ (or $\X \in \Omega_{{\cal V}}^i$). Then we define a scale product
\begin{align}\label{6.3.4}
a\X:=\lceil aX \rceil,\quad a\in \F.
\end{align}

Using (\ref{6.3.3}) and (\ref{6.3.4}), one sees easily that $\Omega_{{\cal V}}$ becomes a vector space:

\begin{thm}\label{t6.3.4} Using the vector addition defined in (\ref{6.3.3}) and the scale product defined in (\ref{6.3.4}), we have
\begin{enumerate}
\item $\Omega_{{\cal V}}$ is a vector space over $\F$;
\item $\Omega_{{\cal V}}^i$, $i=1,2,\cdots$ are the subspaces of $\Omega_{{\cal V}}$;
\item if $i|j$, then $\Omega_{{\cal V}}^i$ is a subspace of $\Omega_{{\cal V}}^j$.
\end{enumerate}
\end{thm}

Let $E\subset {\cal V}$ be a set of vectors. Then
\begin{align}\label{6.3.5}
\lceil E \rceil :=\left\{\X\;\big|\; X\in E \right\}.
\end{align}

The following proposition shows that the vector equivalence keeps the space-subspace relationship unchanged.

\begin{prp}\label{p6.3.5} Assume $E\subset {\cal V}_i$ is a subspace of ${\cal V}_i$. Then $\lceil E \rceil \subset \Omega_{{\cal V}}^i$ is a subspace of $\Omega_{{\cal V}}^i$.
\end{prp}

Definition \ref{d6.3.1} can be translated to its corresponding right one as follows:

\begin{dfn}\label{d6.3.6} Let $X\in {\cal V}_p$ and $Y\in {\cal V}_q$. $t=p\vee q$. Then
\begin{enumerate}
\item the (right) vector addition $\rvplus: {\cal V}\times {\cal V}\ra {\cal V}$ of $X$ and $Y$ is defined as
\begin{align}\label{6.3.6}
X\rvplus Y:=\left(  \I_{t/p} \otimes X\right)+ \left( \I_{t/q} \otimes Y\right);
\end{align}
\item the subtraction $\rvminus$ is defined as
\begin{align}\label{6.3.7}
X\rvminus Y:=X\rvplus (-Y).
\end{align}
\end{enumerate}
\end{dfn}

Then all the arguments in this subsection and the following several subsections can be stated in a parallel way for right addition and corresponding linear spaces.

\subsection{Inner Product and Linear Mappings}

\begin{dfn}\label{d6.0.1} Let $X\in {\cal V}_m\subset {\cal V}$, $Y\in {\cal V}_n\subset {\cal V}$, and $t=m\vee n$. Then the weighted inner product is defied as
\begin{align}\label{6.0.1}
\left<X\;\big|\;Y\right>_W:=\frac{1}{t}\left<X\otimes {\bf 1}_{t/m}\;\big|\;Y\otimes {\bf 1}_{t/n}\right>,
\end{align}
where $\left<X\;\big|\;Y\right>$ is the conventional inner product. Say, if $\F=\R$, $\left<X\;\big|\;Y\right>=X^TY$, and if $\F=\C$, $\left<X\;\big|\;Y\right>=\bar{X}^TY$.
\end{dfn}

\begin{dfn}\label{d6.0.2} Let $\X, \Y\in \Omega_{{\cal V}}$. Their inner product is defined as
\begin{align}\label{6.0.2}
\left<\X\;\big|\;\Y\right>:=\left<X\;\big|\;Y\right>_W, \quad X\in \X, Y\in \Y.
\end{align}
\end{dfn}

It is easy to verify the following proposition, which assures that the Definition \ref{d6.0.2} is reasonable.

\begin{prp}\label{p6.0.3} Equation (\ref{6.0.2}) is well defined. That is, it is independent of the choice of $X$ and $Y$.
\end{prp}

Since $\Omega_{{\cal V}}$ is a vector space, using (\ref{6.0.2}) as an inner product on $\Omega_{{\cal V}}$, then we have the follows:

\begin{prp} The $\Omega_{{\cal V}}$ with inner product defined by (\ref{6.0.2}) is an inner product space.
It is not a Hilbert space.
\end{prp}

\noindent\textit{Proof}. The first part is obvious. As for the second part, we construct a sequence as
$$
\begin{cases}
X_1=a\in \F\\
X_{i+1}=X_i\otimes \I_2+\frac{1}{2^{i+1}}\left(\d_{2^{i+1}}^1-\d_{2^{i+1}}^2\right),\quad i=1,2,\cdots.
\end{cases}
$$
Then we can prove that $\{X_i\}$ is a Cauchy sequence and it does not converge to any $X\in {\cal V}$.
\hfill $\Box$

Given a vector  $X\in {\cal V}$, then it determined a linear mapping $\varphi_X:~{\cal V}\ra \F$ via inner product as
$$
\varphi_X:Y\mapsto \left<X\;\big|\; Y\right>_W.
$$
Similarly, $\X\in \Omega_{{\cal V}}$ can also determine a linear mapping $\varphi_{\X}:~\Omega_{{\cal V}} \ra \F$ as
$$
\varphi_{\X}:\Y\mapsto \left<\X\;\big|\; \Y\right>.
$$
Unfortunately, the inverse is not true. Because  $\Omega_{{\cal V}}$ is an infinite dimensional vector space but any vector in $\Sigma$ has only finite dimensional representatives.

Next, we consider the linear mappings on ${\cal V}$.
It is well known that assume $A\in {\cal M}_{m\times n}$ and $X\in {\cal V}_n$. Then the product $\times:~{\cal M}_{m\times n}\times {\cal V}_n\ra {\cal V}_m$, defined as $(A,~X)\mapsto AX$, for given $A$ is a linear mapping. We intend to generalize such a linear mapping to arbitrary matrix and arbitrary vector.

\begin{dfn}\label{d6.4.1} Let $A\in {\cal M}_{m\times n}$, $X\in {\cal V}_p$, and $t=n\vee p$.  Then the vector product, denoted by $\lvtimes$, is defined as
\begin{align}\label{6.4.1}
A \lvtimes X:=\left(A\otimes I_{t/n}\right)\left(X\otimes \I_{t/p}\right).
\end{align}
\end{dfn}

\begin{rem}\label{r6.4.101}
\begin{enumerate}
\item The vector product defined in (\ref{6.4.1}) is the left vector product. Of course, we can define right vector product, denoted by $\rvtimes$, defined as
 \begin{align}\label{6.4.101}
A \rvtimes X:=\left( I_{t/n}\otimes A\right)\left( \I_{t/p}\otimes X\right).
\end{align}
\item In fact, $\lvtimes$ is a combination of $\sim_{\ell}$ of matrices and $\lra_{\ell}$ of vectors, and $\rvtimes$ is a combination of $\sim_{r}$ of matrices and $\lra_{r}$ of vectors. Of course, we may define two more vector products by combinations of   $\sim_{l}$ with $\lra_{r}$ and  $\sim_{r}$ with $\lra_{l}$ respectively.

\item
Note that when $n=p$, $A\lvtimes X=AX$. That is, the linear mapping defined in (\ref{6.4.1}) is a generalization of conventional linear mapping. It is also true for other vector products.
\item
To avoid similar but tedious arguments, hereafter the default vector product is $\lvtimes$.
\end{enumerate}
\end{rem}

The following proposition is easily verifiable.

\begin{prp}\label{p6.4.2} Consider the vector product
$$\lvtimes:~{\cal M}\times {\cal V}\ra {\cal V}.
$$
\begin{enumerate}
\item It is linear with respect to the second variable, precisely,
\begin{align}\label{6.4.2}
A\lvtimes (a X\lvplus b Y)=aA\lvtimes X\lvplus bA\lvtimes Y,\quad a,b\in \F.
\end{align}
\item Assume both $A,~B\in {\cal M}_{\mu}$, then the vector product is also linear with respect to the first variable, precisely,
\begin{align}\label{6.4.3}
(aA \lplus bB) \lvtimes X=aA \lvtimes X \lvplus bB \lvtimes X.
\end{align}
\end{enumerate}
\end{prp}

The following proposition shows that the vector product is consistent with both M-equivalence and V-equivalence.

\begin{prp}\label{p6.4.3}
Assume $A\sim B$ and $X\lra Y$. Then
\begin{align}\label{6.4.2}
A\lvtimes X \lra B\lvtimes Y.
\end{align}
\end{prp}

\noindent\textit{Proof}. Assume
$A=\Lambda\otimes I_s$, $B=\Lambda\otimes I_{\a}$; $X=\Gamma\otimes \I_t$, $Y=\Gamma\otimes \I_{\b}$, where $\Lambda\in {\cal M}_{n\times p}$ and $\Gamma\in {\cal V}_q$. Denote $\xi=p\vee q$, $\eta=ps\vee qt$, and $\eta=k\xi$.
Then we have
$$
\begin{array}{l}
A\lvtimes X=\left(\Lambda\otimes I_s\right)\lvtimes \left(\Gamma\otimes \I_t  \right)\\
~=\left(\Lambda\otimes I_s\otimes I_{\eta/ps}\right)\left(\Gamma\otimes \I_t \otimes \I_{\eta/qt} \right)\\
~=\left(\Lambda\otimes I_{\xi/p}\otimes I_{k}\right)\left(\Gamma\otimes \I_{\xi/q} \otimes \I_{k} \right)\\
~=\left[\left(\Lambda\otimes I_{\xi/p}\right)\left(\Gamma\otimes \I_{\xi/q} \right)\right]\otimes
\left[I_k\I_k\right]\\
~=\left(\Lambda \lvtimes \Gamma\right)\otimes \I_k.
\end{array}
$$
Hence
$$
A\lvtimes X\lra \Lambda\lvtimes \Gamma.
$$
Similarly, we have
$$
B\lvtimes Y\lra \Lambda\lvtimes \Gamma.
$$
Equation (\ref{6.4.2}) follows.
\hfill $\Box$

The above propositions have an immediate consequence as follows:

\begin{cor}\label{p6.4.4} The vector product $\lvtimes: {\cal M}\times {\cal V}\ra {\cal V}$ can be extended to $\lvtimes: \Sigma_{{\cal M}}\times \Omega_{{\cal V}} \ra \Omega_{{\cal V}}$. Particularly, each $\A\in \Sigma_{{\cal M}}$ determines a linear mapping on the vector space $\Omega_{{\cal V}}$.
\end{cor}

Next, we extend the vector equivalence to matrices.

\begin{dfn}\label{d6.4.5}
\begin{enumerate}
\item Assume $V,~W \in {\cal M}_{\cdot \times n}$. $V$ and $W$ are said to be vector equivalent, denoted by $V\lra W$, if there exist $\I_s$ and $\I_t$ such that
\begin{align}\label{6.4.3}
V\otimes \I_s=W\otimes \I_t.
\end{align}
\item The equivalence class of $V$ is denoted as
\begin{align}\label{6.4.4}
\lceil V\rceil :=\left\{W\;\big|\; W\lra V\right\}.
\end{align}
\end{enumerate}
\end{dfn}

\begin{rem}\label{r6.4.6}
\begin{enumerate}
\item It is clear the equivalence defined by (\ref{6.4.3}) is left vector equivalence $\lra_{\ell}$. The right vector equivalence $\lra_r$ can be defined similarly. Moreover, the right vector equivalence class is denoted by $\lceil \cdot \rceil_r$.
\item A matrix can be considered as a linear mapping, and correspondingly the M-equivalence is considered. A matrix can also be considered as a vector subspace generated by its columns. In this way its V-equivalence, defined by (\ref{6.4.3})--(\ref{6.4.4}), makes sense. Denote the equivalence class as
$$
\Omega_{{\cal M}}:={\cal M}/\lra,
$$
and
$$
\Omega_{{\cal M}}^{n}:={\cal M}_{\cdot\times n}/\lra.
$$
\end{enumerate}
\end{rem}

All the results about vector space ${\cal V}$ can be extended to ${\cal M}_{\cdot\times n}$ for $n\in \N$. They are briefly summarized as follows:

\begin{prp}\label{p6.4.7}
\begin{enumerate}
\item Assume $V\in {\cal M}_{r\times n}$, $W\in {\cal M}_{s\times n}$, $V\lra W$ and $r|s$. Then we said $V$ is a divisor of $W$ or $W$ is a multiple of $V$. Moreover, an order  determined by this is denoted as $V\sqsubseteq W$.
\item $(\lceil V\rceil,\sqsubseteq)$ is a lattice.
\item Assume $r|s$, then an order is given in ${\cal M}_{\cdot \times n}$ as
$$
{\cal M}_{r\times n}\sqsubseteq {\cal M}_{s\times n}.
$$
\item $\left({\cal M}_{\cdot \times n},\sqsubseteq \right)$ is a lattice.
\item Assume $V\in {\cal M}_{r\times n}$ is irreducible. then
$\varphi: \lceil V\rceil \ra {\cal M}_{\cdot \times n}$ defined as $V_i\mapsto {\cal M}_{ir\otimes n}$ is a lattice homomorphism.
\end{enumerate}
\end{prp}

The operators $\lvplus$ and $\lvtimes$ can also be defined in a similar way as for vectors:

\begin{dfn}\label{d6.4.8}
\begin{enumerate}
\item Assume $V~\in {\cal M}_{p\times n}$ and $W\in {\cal M}_{q\times n}$, $t=p\vee q$. The vector addition is defined as
$$
V\lvplus W:=\left(V\otimes \I_{t/p}\right)+\left(W\otimes \I_{t/q}\right).
$$
\item Let $A\in {\cal M}_{m\times n}$, $V\in {\cal M}_{p\times q}$, and $t=n\vee p$.  Then the vector product of $A$ and $V$ is defined as
\begin{align}\label{6.4.5}
A \lvtimes V:=\left(A\otimes I_{t/n}\right)\left(V\otimes \I_{t/p}\right).
\end{align}
\end{enumerate}
\end{dfn}

We can denote the equivalence class of ${\cal M}_{\cdot\times n}$ as
$$
\Sigma_{{\cal M}}^n:={\cal M}_{\cdot\times n}/\lra.
$$

For vector product and vector addition we also have the distributive law with respect to these two operators for matrix case.

\begin{prp}\label{p6.4.9} $\lvtimes: {\cal M}\times {\cal M}\ra {\cal M}$ is distributive. Precisely, Let $A,~B\in {\cal M}_{\mu}$ and $C\in {\cal M}_{n\times \b}$, $D\in {\cal M}_{p\times \b}$. Then
\begin{align}\label{6.7.11}
(aA\lplus bB)\lvtimes C=a A\lvtimes C \lvplus b B\lvtimes C,
\quad a,b\in \F;\\
\label{6.7.12}
A\lvtimes (aC\lvplus bD)=a A\lvtimes C \lvplus b A\lvtimes D,
\quad a,b\in \F.
\end{align}
\end{prp}

\noindent\textit{Proof}. We prove (\ref{6.7.11}). The proof of (\ref{6.7.12}) is similar.
Denote $m\vee n=s$, $m\vee p=r$, $n\vee p)=t$, and $m\vee n \vee p=\xi$.
Then for (\ref{6.7.11}) we have
$$
\begin{array}{l}
LHS=\left(aA\otimes I_{s/m}+bB\otimes I_{s/n}\right)\lvtimes C\\
~=\left[\left(aA\otimes I_{s/m}+bB\otimes I_{s/n}\right)\otimes I_{\xi/s}\right]
 \left[C\otimes {\bf 1}_{\xi/p}\right]\\
~=a\left(A\otimes I_{\xi/m}\right)\left(C\otimes {\bf 1}_{\xi/p}\right)+b\left(B\otimes I_{\xi/n}\right)\left(C\otimes {\bf 1}_{\xi/p}\right)\\
~=a\left[\left(A\otimes I_{s/m}\right)\left(C\otimes {\bf 1}_{s/p}\right)\right]\otimes {\bf 1}_{\xi/s}\\
~+b\left[\left(B\otimes I_{t/n}\right)(C\otimes {\bf 1}_{t/p})\right]\otimes {\bf 1}_{\xi/t}\\
~=a\left[\left(A\otimes I_{s/m}\right)\left(C\otimes {\bf 1}_{s/p}\right)\right]
\lvplus b\left[\left(B\otimes I_{t/n}\right)(C\otimes {\bf 1}_{t/p})\right]\\
~=RHS.
\end{array}
$$
\hfill $\Box$

The distributive law can be extended to equivalence spaces as follows.

\begin{cor}\label{c6.4.10} $\lvtimes: \Sigma_{{\cal M}}\times \Sigma_{{\cal M}}^n\ra \Sigma_{{\cal M}}^n$ is distributive. Precisely, Let $\A,~\B\in \Sigma_{\mu}$ and $\lceil C\rceil, \lceil D\rceil \in \Sigma_{{\cal M}}^{\b}$. Then
\begin{align}\label{6.7.13}
\begin{array}{ccl}
(a\A\lplus b\B)\lvtimes \lceil C\rceil&=&a \A\lvtimes \lceil C\rceil \lvplus b \B\lvtimes \lceil C\rceil ,\\
~&~& a,b\in \F.
\end{array}\\
\label{6.7.14}
\begin{array}{ccl}
\A\lvtimes (a\lceil C\rceil \lvplus b\lceil D\rceil )&=&a \A\lvtimes \lceil C\rceil \lvplus b \A\lvtimes \lceil D\rceil,\\
~&~& a,b\in \F.
\end{array}
\end{align}
\end{cor}

\subsection{Type-1 Invariant Subspace}

Given $A\in {\cal M}_{\mu}^i$, we seek a subspace $S\subset {\cal V}$ which is $A$ invariant.

\begin{dfn}\label{d6.5.1} Let $S\subset {\cal V}$ be a  vector subspace. If
$$A\lvtimes S\subset S,
$$ $S$ is called an $A$-invariant subspace.  Moreover, if $S\subset {\cal V}_t$, it is called the type-1 invariant subspace; Otherwise, it is called the type-2 invariant subspace.
\end{dfn}

This subsection considers the type-1 invariant subspace only.

\begin{prp}\label{p6.5.2} Let $A\in {\cal M}_{\mu}$ and $S={\cal V}_t$. Then $S$ is $A$-invariant, if and only if,
\begin{itemize}
\item[(i)]
\begin{align}\label{6.5.1}
\mu_y=1;
\end{align}
\item[(ii)] $A\in {\cal M}^i_{\mu}$, where $i$ satisfies
\begin{align}\label{6.5.2}
i\mu_x\vee t=t \mu_x.
\end{align}
\end{itemize}
\end{prp}

\noindent\textit{Proof}. (Necessity) Assume $\xi=i\mu_x\vee t$, by definition we have
$$
\begin{array}{l}
A\lvtimes X=\left(A\otimes I_{\xi/i\mu_x}\right)\left(X\otimes \I_{\xi/t}\right)\\
~\in {\cal V}_t,\quad X\in {\cal V}_t.
\end{array}
$$
Hence we have
\begin{align}\label{6.5.3}
i\mu_y\left(\frac{\xi}{i\mu_x}\right)=t.
\end{align}
It follows from (\ref{6.5.3}) that
\begin{align}\label{6.5.4}
\frac{\xi}{t}=\frac{\mu_x}{\mu_y}.
\end{align}
Since $\mu_x$ and $\mu_y$ are co-prime and the left hand side of (\ref{6.5.4}) is an integer, we have $\mu_y=1$. It follows that
$$
\xi=t\mu_x.
$$

(Sufficiency)  Assume (\ref{6.5.2}) holds, then (\ref{6.5.3}) holds. It follows that
$$
A\lvtimes X \in {\cal V}_t,\quad \mbox{when}~ X\in {\cal V}_t.
$$
\hfill $\Box$

Assume $A\in {\cal M}_{\mu}$ and $S={\cal V}_t$, a natural question is: Can we find $A$ such that $S$ is $A$-invariant? According to Proposition \ref{p6.5.2}, we know that it is necessary that $\mu_y=1$. Let $k_1,\cdots,k_{\ell}$ be the prime divisors of $\mu_x\wedge ~t$, then we have
\begin{align}\label{6.5.5}
\mu_x=k_1^{\a_1}\cdots k_{\ell}^{\a_{\ell}}p;~t=k_1^{\b_1}\cdots k_{\ell}^{\b_{\ell}}q,
\end{align}
where $p,~q$ are co-prime and $k_i\nmid p$, $k_i\nmid q$, $\forall i$.

Now it is obvious that to meet (\ref{6.5.2}) it is necessary and sufficient that
\begin{align}\label{6.5.6}
i=k_1^{\b_1}\cdots k_{\ell}^{\b_{\ell}}\lambda,
\end{align}
where $\lambda|(pq)$.

Summarizing the above argument we have the following result.

\begin{prp}\label{p6.5.3} Assume $A\in {\cal M}_{\mu}$, $S={\cal V}_t$. $S$ is $A$-invariant, if and only if, (i) $\mu_y=1$, and (ii)
$A\in {\cal M}_{\mu}^i$, where $i$ satisfies (\ref{6.5.6}).
\end{prp}

\begin{rem}\label{r6.5.4} Assume $A\in {\cal M}_1$, $S={\cal V}_t$. Using Proposition \ref{p6.5.3}, it is clear that $S$ is $A$-invariant, if and only if, $A\in {\cal M}_1^i$ with
$$
i\in \{\ell\;\big|\;\ell|t\}.
$$
Particularly, when $i=1$ then $A=a$ is a number. So it is a scale product of $S$, i.e, for $V\in S$ we have $a \lvtimes V=aV\in S$.
When $i=t$ then $A\in {\cal M}_{t\times t}$. So we have $A\lvtimes V=AV\in S$. This is the classical linear mapping on ${\cal V}_t$. We call these two linear mapping the standard linear mapping.
The following example shows that there are lots of non-standard linear mappings.
\end{rem}

\begin{exa}\label{e6.5.401}
\begin{enumerate}
\item Assume $\mu=0.5$, $t=6$. Then $\mu_y=1$, $\mu_x=2$. $\mu_x\wedge t=2$. Using (\ref{6.5.5}), $\mu_x=2\times 1$ and $t=2\times 3$. That is, $p=1$, $q=3$.  According to (\ref{6.5.6}) $i=2^{\b_1}\lambda$, where $\b_1=1$, $\lambda=1$ or $\lambda=3$. According to Proposition \ref{p6.5.3}, ${\cal V}_6$ is $A$-invariant, if and only if, $A\in {\cal M}_{0.5}^2$ or $A\in {\cal M}_{0.5}^6$.
\item Next, we give a numerical example: Assume
\begin{align}\label{6.5.7}
A=\begin{bmatrix}
1&-1&0&0\\
0&0&1&0
\end{bmatrix}\in {\cal M}_{0.5}^2.
\end{align}
Then $\R^6$ is $A$-invariant space. For instance, let
\begin{align}\label{6.5.8}
X=\begin{bmatrix}1+i,&2,&1-i,&0,&0,&0
\end{bmatrix}^T\in \C^6.
\end{align}
Then
\begin{align}\label{6.5.9}
A\lvtimes X=iX.
\end{align}
\end{enumerate}

\end{exa}

Motivated by (\ref{6.5.9}), we give the following definition.

\begin{dfn} \label{d6.5.5} Assume $A\in {\cal M}$ and $X\in {\cal V}$. If
\begin{align}\label{6.5.10}
A\lvtimes X=\a X, \quad \a\in \F, X\neq 0,
\end{align}
then $\a$ is called an eigenvalue of $A$, and $X$ is called an eigenvector of $A$ with respect to $\a$.
\end{dfn}

\begin{exa}\label{e6.5.6} Recall Example \ref{e6.5.401}. In fact it is easy to verify that matrix $A$ in (\ref{6.5.7}), as a linear mapping on $\R^6$, has $6$ eigenvalues: Precisely,
$$
\sigma(A)=\left\{i,-i,0,0,0,1\right\}.
$$
Correspondingly, the first eigenvector is $X$, defined in (\ref{6.5.8}). The other $5$ eigenvectors are
$(1-i,2,1+i,0,0,0)^T$, $(1,1,1,0,0,0)^T$, $(0,0,0,0,1,0)^T$ (this is a root vector), $(0,0,0,0,0,1)^T$, and $(0,0,0,1,1,1)^T$, respectively.
\end{exa}

Assume $S={\cal V}_t$ is $A$-invariant with respect to $\lvtimes$, that is,
\begin{align}\label{6.5.11}
A\lvtimes S\subset S.
\end{align}
Then the restriction $A|_S$ is a linear mapping on $S$. It follows that there exists a matrix, denoted by $A|_{t}\in {\cal M}_{t\times t}$, such that $A|_S$ is equivalent to $A|_{t}$. We state it as a proposition.

\begin{prp}\label{p6.5.7} Assume $A\in {\cal M}$ and $S={\cal V}_t$ is $A$-invariant. Then there exists a unique matrix $A|_{t}\in {\cal M}_{t\times t}$, such that $A|_S$ is equivalent to $A|_{t}$. Precisely,
\begin{align}\label{6.5.12}
A \lvtimes X=A|_{t}X, \quad \forall X\in S.
\end{align}
$A|_{t}$ is called the realization of $A$ on $S={\cal V}_{t}$.
\end{prp}

\begin{rem}\label{r6.5.8}
\begin{enumerate}
\item To calculate $A|_{t}$ from $A$ is easy. In fact, it is clear that
\begin{align}\label{6.5.13}
\Col_i\left(A|_{t}\right)=A\lvtimes \d_t^i, \quad i=1,\cdots, t.
\end{align}
\item Consider Example \ref{e6.5.401} again.
$$
\begin{array}{l}
\Col_1(A|_{6})=A\lvtimes \d_6^1\\
~=\left[
\begin{bmatrix}
1&-1&0&0\\
0&0&1&0\end{bmatrix}\otimes I_3
\right]\left[\d_6^1\otimes \I_2\right]\\
~=(1,1,0,0,0,0)^T.
\end{array}
$$
Similarly, we can calculate all other columns. Finally, we have
$$A|_{6}=
\begin{bmatrix}
     1&    -1&     0&     0&     0&     0\\
     1&     0&    -1&     0&     0&     0\\
     0&     1&    -1&     0&     0&     0\\
     0&     0&     0&     1&     0&     0\\
     0&     0&     0&     1&     0&     0\\
     0&     0&     0&     0&     1&     0
\end{bmatrix}
$$
\end{enumerate}
\end{rem}

Finally, given a matrix $A\in {\cal M}_{\mu}^i$, we would like to know whether it has (at least one) type-1 invariant subspace $S={\cal V}_t$?

Then as a corollary of Proposition \ref{p6.5.2}, we can prove the following result.

\begin{cor}\label{c6.5.9}Assume $A\in {\cal M}_{\mu}^i$, then
$A$ has (at least one) type-1 invariant subspace $S={\cal V}_t$, if and only if,
\begin{align}\label{6.5.14}
\mu_y=1.
\end{align}
\end{cor}

\noindent\textit{Proof.}
According to Proposition \ref{p6.5.2},  $\mu_y=1$ is obvious necessary. We prove it is also sufficient.
Assume $i$ is factorized into its prime factors as
\begin{align}\label{6.5.15}
i=\prod_{j=1}^{n}i_j^{k_j},
\end{align}
and correspondingly, $\mu_x$ is factorized as
\begin{align}\label{6.5.16}
\mu_x=\prod_{j=1}^{n}i_j^{r_j}p,
\end{align}
where $p$ is co-prime with $i$; $t$ is factorized as
\begin{align}\label{6.5.17}
t=\prod_{j=1}^{n}i_j^{t_j}q,
\end{align}
where $q$ is co-prime with $i$.

Using Proposition \ref{p6.5.2} again, we have only to prove that there exists at least one $t$ satisfying
(\ref{6.5.2}). Calculate that
$$
\begin{array}{ccl}
i\mu_x\vee t&=&\prod_{j=1}^n i_j^{\max(r_j+k_j,~t_j)} (p\vee q);\\
\mu_x t&=&\prod_{j=1}^n i_j^{r_j+t_j}pq.
\end{array}
$$
To meet (\ref{6.5.2}) a necessary condition is: $p$ and $q$ are co-prime. Next, fix $j$, we consider two cases: (i) $t_j>k_j+r_j$: Then on the LHS (left hand side) of (\ref{6.5.2}) we have factor $i_j^{t_j}$ and on the RHS of (\ref{6.5.2}) we have factor $i_j^{r_j+t_j}$. Hence, as long as $r_j=0$, we can choose $t_j>k_j$ to meet (\ref{6.5.2}). (ii) $t_j<k_j+r_j$: Then on the LHS we have factor $i_j^{k_j+r_j}$ and on the RHS  we have factor $i_j^{r_j+t_j}$. Hence, as long as $t_j=k_j$, (\ref{6.5.2}) is satisfied.
\hfill $\Box$

Using above notations, we also have the following result

\begin{cor}\label{c6.5.10}
Assume $\mu_y=1$. Then ${\cal V}_t$ is $A$-invariant, if and only if,
(i) for $t_j=0$, the corresponding $r_j\geq k_j$; (ii) for $t_j>0$, the corresponding $r_j=k_j$.
\end{cor}

\begin{exa}\label{e6.5.11} Recall Example \ref{e6.5.401} again.

\begin{enumerate}
\item Since the matrix $A$ defined in (\ref{6.5.7}) is in ${\cal M}_{0.5}^2$. We have $i=2$, ${\mu}_y=1$, ${\mu}_x=2=ip$ and hence $p=1$. According to Corollary \ref{c6.5.10}, $S={\cal V}_t$ is $A$-invariant, if and only if, $t=iq=2q$ and $q$ is co-prime with $i=2$. Hence
$$
{\cal V}_{2(2n+1)},\quad n=0,1,2,\cdots,
$$
are type-1 invariant subspaces of $A$.

\item Assume $q=5$. Then the restriction is
$$A|_{\R^{10}}=A|_{10}=\begin{bmatrix}
     1&     0&    -1&     0&     0&     0&     0&     0&     0&     0\\
     1&     0&     0&    -1&     0&     0&     0&     0&     0&     0\\
     0&     1&     0&    -1&     0&     0&     0&     0&     0&     0\\
     0&     1&     0&     0&    -1&     0&     0&     0&     0&     0\\
     0&     0&     1&     0&    -1&     0&     0&     0&     0&     0\\
     0&     0&     0&     0&     0&     1&     0&     0&     0&     0\\
     0&     0&     0&     0&     0&     1&     0&     0&     0&     0\\
     0&     0&     0&     0&     0&     0&     1&     0&     0&     0\\
     0&     0&     0&     0&     0&     0&     1&     0&     0&     0\\
     0&     0&     0&     0&     0&     0&     0&     1&     0&     0
\end{bmatrix}
$$
The eigenvalues are
$$
\sigma(A|_{\R^{10}})=\{-1,i,-i,0,1,0,0,0,0,1\}.
$$
The corresponding eigenvectors are:
$$
\begin{array}{l}
E_1=(0,1,0,1,2,0,0,0,0,0)^T\\
E_2=(0.3162+0.1054i,0.5270,0.4216 - 0.2108i,\\
~0.3162 - 0.4216i, 0.1054 - 0.3162i,0,0,0,0,0)^T\\
E_3=(0.3162-0.1054i,0.5270,0.4216 + 0.2108i,\\
~0.3162 + 0.4216i, 0.1054 + 0.3162i,0,0,0,0,0)^T\\
E_4=(1,1,1,1,1,0,0,0,0,0)^T\\
E_5=(2,1,0,1,0,0,0,0,0,0)^T\\
R_6=(0,0,0,0,0,0,1,0,0,0)^T\\
R_7=(0,0,0,0,0,0,1,1,0,0)^T\\
E_8=(0,0,0,0,0,0,0,0,0,1)^T\\
E_9=(0,0,0,0,0,0,0,0,1,0)^T\\
E_{10}=(0,0,0,0,0,1,1,1,1,1)^T.
\end{array}
$$
Note that $R_6$, $R_7$ are two root vectors. That is,
$$
A\lvtimes R_6=R_7;\quad A\lvtimes R_7=E_8.
$$
\end{enumerate}
\end{exa}

\begin{rem}\label{r6.5.11} In fact, the set of type-1 $A$-invariant subspaces depends only on the shape of $A$. Hence we define
$$
{\cal I}_{\mu}^i:=\left\{{\cal V}_t\;|\; {\cal V}_t ~\mbox{is}~ A\in {\cal M}_{\mu}^i ~\mbox{invariant} \right\}.
$$
We can also briefly call ${\cal I}_{\mu}^i$ the set of ${\cal M}_{\mu}^i$-invariant subspaces.
\end{rem}

Before ending this subsection we consider a general type-1 invariant subspace $S\subset {\cal V}_t$. The following proposition is obvious.

\begin{prp}\label{p6.5.12} If $S\subset {\cal V}_t$ is an $A$-invariant subspace, then ${\cal V}_t$ is also an $A$-invariant subspace.
\end{prp}

Because of Proposition \ref{p6.5.12} searching $S$ becomes a classical problem. Because we can first find a matrix $P\in {\cal M}_{t\times t}$, which is equivalent to $A|_{{\cal V}_t}$. Then $S$ must be a classical invariant subspace of $P$.

\subsection{Type-2 Invariant Subspace}

Denote the set of Type-1 $A$-invariant subspaces as
$$
{\cal I}_A:=\left\{{\cal V}_s\;\big|\; {\cal V}_s ~\mbox{is}~A\mbox{-invariant}\right\}.
$$
Assume $A\in {\cal M}_{\mu}^i$. Then ${\cal I}_A={\cal I}_{\mu}^i$.

To assure ${\cal I}_A\neq \emptyset$, through this subsection we assume $\mu_y=1$. We give a definition for this.

\begin{dfn}\label{d6.6.0}
\begin{enumerate}
\item Assume $A\in {\cal M}_{\mu}$. $A$ is said to be a bounded  operator (or briefly, $A$ is bounded,) if $\mu_y=1$.

\item A sequence $\left\{{\cal V}_i\;\big|\; i=1,2,\cdots\right\}$, is called an $A$ generated sequence if
$$
A\lvtimes {\cal V}_i\subset {\cal V}_{i+1},\quad i=1,2,\cdots.
$$
\item A finite sequence  $\left\{{\cal V}_i\;\big|\; i=1,2,\cdots,p\right\}$, is called an $A$ generated loop if ${\cal V}_p={\cal V}_1$.
\end{enumerate}
\end{dfn}

\begin{lem}\label{l6.6.1} Assume there is an $A$ generated loop ${\cal V}_p$, ${\cal V}_{q_1}$, $\cdots$, ${\cal V}_{q_r}$, ${\cal V}_p$, as depicted in (\ref{6.6.1}).
\begin{align}\label{6.6.1}
{\cal V}_p\xrightarrow{\mbox{A}} {\cal V}_{q_1}\xrightarrow{\mbox{A}}\cdots \xrightarrow{\mbox{A}} {\cal V}_{q_r}\xrightarrow{\mbox{A}} {\cal V}_{p}
\end{align}
Then
\begin{align}\label{6.6.2}
q_j=p,\quad j=1,\cdots,r.
\end{align}
\end{lem}

\noindent\textit{Proof}. Assume $A\in {\cal M}_{\mu}^i$. According to Definition \ref{d6.4.1}, we have the following dimension relationship:
\begin{align}\label{6.6.3000}
\begin{array}{lll}
s_0:=i\mu_x\vee p&\Rightarrow&q_1=\mu s_0;\\
s_1:=i\mu_x\vee q_1&\Rightarrow&q_2=\mu s_1;\\
~&\vdots&~\\
s_{r-1}:=i\mu_x\vee q_{r-1}&\Rightarrow&q_r=\mu s_{r-1};\\
s_r:=i\mu_x\vee q_r&\Rightarrow&p=\mu s_r.\\
\end{array}
\end{align}
Next, we
$$
\begin{array}{llll}
\mbox{set}& s_0:=t_0p&\mbox{then}& q_1=\mu t_0p;\\
\mbox{set}& s_1:=t_1q_1&\mbox{then}& q_2=\mu^2t_1 t_0p;\\
~&\vdots&~&~\\
\mbox{set}& s_{r-1}:=t_{r-1}q_{r-1}&\mbox{then}& q_r=\mu^{r}t_{r-1}\cdots t_1 t_0p;\\
\mbox{set}& s_r:=t_rq_r&\mbox{then}& p=\mu^{r+1}t_r\cdots t_1 t_0p.\\
\end{array}
$$
We conclude that
\begin{align}\label{6.6.4}
\mu^{r+1}t_rt_{r-1}\cdots t_0=1.
\end{align}
Equivalently, we have
$$
\frac{\mu_x^{r+1}}{\mu_y^{r+1}}=t_rt_{r-1}\cdots t_0.
$$
It follows that
$$
\mu_y=1.
$$

Define
$$
s_r=i \mu_x\vee q_r:=i\mu_x \xi,
$$
where $\xi\in \N$. Then from the last equation of (\ref{6.6.3000}) we have
\begin{align}\label{6.6.5}
p=i\xi.
\end{align}
That is,
$$
\mu(i\mu_x\vee q_r)=i\xi.
$$
Using (\ref{6.6.4}), and the expression
$$
q_r=\mu^{r}t_{r-1}\cdots t_1 t_0p,
$$
we have
$$
\left(\mu_x \vee \frac{\mu_x \xi}{t_r}\right)=\xi \mu_x.
$$
From above it is clear that
\begin{align}\label{6.6.6}
\begin{array}{l}
t_r\big|\mu_x\\
\mu_x\wedge\xi=1.
\end{array}
\end{align}

Next, using last two equation in (\ref{6.6.3000}), we have
$$
\mu\left(i\mu_x\vee \mu\lcm(i\mu_x,q_{r-1})\right)=p=i\xi.
$$
Similar to the above argument, we have
$$
\xi\mu_x \big| \left(\mu_x \vee \mu(\mu \vee \frac{q_{r-1}}{i})\right).
$$
Hence
$$
\xi\mu_x \big| \left(\mu_x \vee \mu(\mu \vee \frac{\mu_x^2}{t_rt_{r-1}}\xi)\right).
$$
To meet this requirement, it is necessary that
$$
t_rt_{r-1}\big|\mu_x^2.
$$
Continuing this process, finally we have
\begin{align}\label{6.6.7}
t_rt_{r-1}\cdots t_s\big|\mu_x^{r-s+1},\quad s=r-1,r-2,\cdots,0.
\end{align}
Combining (\ref{6.6.4}) with (\ref{6.6.7}) yields that
$$
t_s=\mu_x, \quad s=0,1,\cdots,r.
$$
That is, $q_1=q_2=\cdots=q_r=p$.
\hfill $\Box$
%
%

\begin{thm}\label{t6.6.2} A finite dimensional subspace $S\subset {\cal V}$ is $A$-invariant, if and only if, $S$ has the following structure:
\begin{align}\label{6.6.9}
S=\oplus_{i=1}^{\ell}S^i,
\end{align}
where
$$
A\lvtimes S^{i}\subset S^{i+1}, \quad i=1,\cdots,\ell-1.
$$
\end{thm}

\noindent\textit{Proof}. Since $S$ is of finite dimension, there are only finite ${\cal V}_{t_i}$ such that
$$
S^j:=S\cap {\cal V}_{t_j}\neq \{0\}.
$$
Now for each $0\neq X_0\in S^j\subset {\cal V}_{t_j}$ we construct $X_1:=A\lvtimes X\in V_{t_r}$ for certain $t_r$. Note that $S$ is $A$-invariant, if for all $X_0\in S^j$, we have $t_r=t_j$, then this $S^j=S^{\ell}$ is the end element in the sequence. Otherwise, we can find a successor $S^{r}=S\cap{\cal V}_{t_r}$. Note that since there are only finite $S^j$, according to Lemma \ref{l6.6.1}, starting from $X_0\in S$ there are only finite sequence of differen $S^j$ till it reach an $A$-invariant $S^{\ell}$ (equivalently, $A$-invariant ${\cal V}_{t_{\ell}}$). The claim follows.  \hfill $\Box$

\subsection{Higher Order Linear Mapping}

\begin{dfn}\label{d6.7.1} Let $A\in {\cal M}_{\mu}^i$, ${\cal V}_t$ is $A$-invariant subspace of ${\cal V}$. That is,
$A:{\cal V}_t\ra {\cal V}_t$ is a linear mapping. The higher order linear mapping of $A$, is defined as
\begin{align}\label{6.7.1}
\begin{cases}
A^{[1]}\lvtimes X:=A\lvtimes X,\quad X\in {\cal V}_t\\
A^{[k+1]}\lvtimes X:=A\lvtimes\left(A^{[k]}\lvtimes X\right),\quad k\geq 1.
\end{cases}
\end{align}
\end{dfn}

\begin{dfn}\label{d6.7.2} Let $X\in {\cal V}$. The $A$-sequence of $X$ is the sequence $\left\{X_i\right\}$, where
$$
\begin{cases}
X_0=X,\\
X_{i+1}=A\lvtimes X_i,\quad i=0,1,2,\cdots.
\end{cases}
$$
\end{dfn}

Using notations (\ref{6.5.15})--(\ref{6.5.17}), we have the following result.

\begin{lem}\label{l6.7.3} Assume $A\in {\cal M}_{\mu}^i$ is bounded. $X\in {\cal V}_t$,
where $i$, $\mu= \frac{1}{\mu_x}$, and $t$ are described by (\ref{6.5.15})--(\ref{6.5.17}).
Then $A\lvtimes X\in {\cal V}_s\in {\cal I}_A$, if and only if,
for each $0<j<n$, one of the following is true:
\begin{align}\label{6.7.2}
r_j=0;
\end{align}
or
\begin{align}\label{6.7.3}
t_j\leq k_j+r_j.
\end{align}
\end{lem}

\noindent\textit{Proof}.
Since
$$
\begin{array}{l}
i\mu_x\vee t=\left(\prod_{j=1}^n i_j^{\max(k_j+r_j,t_j)}q \vee \prod_{j=1}^ni_j^{t_j}q\right)\\
~=\prod_{j=1}^ni_j^{\max(k_j+r_j,t_j)}(p\vee q),
\end{array}
$$
we have
$$
\begin{array}{l}
s=(i\mu_x\vee t)\mu\\
~=\prod_{j=1}^ni_j^{\max(k_j+r_j,t_j)-r_j}(p\vee q)/p.
\end{array}
$$
Then we can calculate that
$$
i\mu_x\vee s=\prod_{j=1}^ni_j^{\max\left(k_j+r_j, max(k_j+r_j,t_j)-r_j\right)}(p\vee q),
$$
and
$$
\mu_x s=\prod_{j=1}^ni_j^{\max\left(k_j+r_j, t_j\right)}(p\vee q).
$$
Note that ${\cal V}_s\in {\cal I}_A$ is $A$-invariant. Using Proposition \ref{p6.5.2},   (\ref{6.5.2}) leads to
\begin{align}\label{6.7.4}
\max\left(k_j+r_j,\max(k_j+r_j,t_j)-r_j\right)=\max\left(k_j+r_j,t_j\right).
\end{align}
\begin{enumerate}
\item[Case 1]: $t_j>k_j+2r_j$,
(\ref{6.7.4}) leads to $r_j=0$. Hence, we have
\begin{align}\label{6.7.5}
r_j=0~\mbox{and}~t_j>k_j.
\end{align}
\item[Case 2]: $k_j+r_j\leq t_j\leq k_j+2r_j$,
which leads to $k_j+r_j=t_j$.
\item[Case 3]: $t_j<k_j+r_j$, which assures (\ref{6.7.4}).
Combining Case 2 and Case 3 yields
\begin{align}\label{6.7.6}
t_j\leq k_j+r_j.
\end{align}
\end{enumerate}
Note that when $r_j=0$, if $t_j\leq k_j$ we have (\ref{6.7.6}). Hence  $t_j\leq k_j$ is also allowed. The conclusion follows.
\hfill $\Box$

The following result is important.

\begin{thm}\label{t6.7.4} Let $A\in {\cal M}_{\mu}^i$ be bounded. Then for any $X\in {\cal V}_{t^0}$ the $A$-sequence of $X$ will enter a ${\cal V}_t\in {\cal I}_A$ at finite steps.
\end{thm}

\noindent\textit{Proof}.  Assume $X_1:= A\lvtimes X \in {\cal V}_{\tilde{t^1}}$.  Using notations (\ref{6.5.15})--(\ref{6.5.17}), it is easy to calculate that after one step the $j$-th index of $t$ becomes
$$
t^1_j=\max\left(k_j+r_j,t^0_j)\right)-r_j.
$$
Assume $r_j=0$, this component already meets the requirement of Lemma \ref{l6.7.3}. Assume for some $j$, $r_j>0$ and $t^0_j>k_j+r_j$, then we have
$$
t^1_j=t^0_j-r_j<t^0_j.
$$
Hence after finite times, say $k$, the $j$-th index of $X_k$, denoted by $t^k_j$, satisfies
\begin{align}\label{6.7.7}
t^k_j=t^0_j-kr_j,
\end{align}
will satisfy (\ref{6.7.3}), and as long as (\ref{6.7.3}) holds, $t^s_j=t^k_j$  $\forall s>k$. Hence, after finite steps either (\ref{6.7.2}) or (\ref{6.7.3}) (or both) is satisfied. Then at the next step the sequence enters into ${\cal V}_t\in {\cal I}_A$.
\hfill $\Box$

\begin{dfn}\label{d6.7.5} Given a polynomial
\begin{align}\label{6.7.8}
p(x)=x^n+c_{n-1}x^{n-1}+\cdots+c_1x+c_0,
\end{align}
a matrix $A\in {\cal M}$ and a vector $X\in {\cal  V}$.
\begin{enumerate}
\item $p(x)$ is called an $A$-annihilator of $X$, if
\begin{align}\label{6.7.8}
\begin{array}{ccl}
p(A)X&:=&A^{[n]}\lvtimes X \lvplus c_{n-1} A^{[n-1]}\lvtimes X\\
~&~&\lvplus \cdots \lvplus c_1A\lvtimes X\lvplus c_0=0.
\end{array}
\end{align}
\item Assume $q(x)$ is the $A$-annihilator of $X$ with minimum degree, then $q(x)$ is called the minimum
$A$-annihilator of $X$.
\end{enumerate}
\end{dfn}

\begin{rem}\label{r6.7.401} Theorem \ref{t6.7.4} shows why $A\in {\cal M}_{\mu}$ with $\mu_y=1$ is called a bounded operator. In fact, it is necessary and sufficient for $A$ to have a finite dimensional invariant subspace of either type-1 or type-2. We also know that if $A$ has type-2 invariant subspace, it also has type-1 invariant subspace. If $\mu_y\neq 1$, $A$ is called an unbounded operator.
\end{rem}

The following result is obvious:

\begin{prp}\label{p6.7.6} The minimum $A$-annihilator of $X$ divides any  $A$-annihilator of $X$.
\end{prp}

The following result is an immediate consequence of Theorem \ref{t6.7.4}.

\begin{cor}\label{c6.7.7} Assume $A$ is bounded, then for any $X\in {\cal V}$ there exists at least one
$A$-annihilator of $X$.
\end{cor}

\noindent\textit{Proof.} According to Theorem \ref{t6.7.4}, there is a finite $k$ such that $A^{[k]}X\in {\cal V}_{s}$ with ${\cal V}_s$ being $A$-invariant. Now in ${\cal V}_s$ assume the minimum annihilator polynomial for $A^{[k]}X$ is $q(x)$, then $p(x)=x^kq(x)$ is an $A$-annihilator of $X$.
\hfill $\Box$

\begin{exa}\label{e6.7.8}
\begin{enumerate}
\item Assume $A\in {\cal M}_{2/3}^1$. Since $\mu_y=2\neq 1$, we know any $X$ does not have its $A$-annihilator.

Now assume $X_0\in {\cal V}_k$, where $k=3^sp$, and $3,~p$ are co-prime. Then it is easy to see that the $A$ sequence of $X_0$ has the dimensions, $\dim(X_i):=d_i$, which are: $d_1=2\times 3^{s-1}p$, $d_2=2^2\times 3^{s-2}p$ $\cdots$ $d_s=2^sp$, $d_{s+1}=2^{s+1}p$, $d_{s+2}=2^{s+2}p$, $\cdots$. It can not reach a ${\cal V}_t\in {\cal I}_A$.

\item Given
$$
A=\begin{bmatrix}
1&0&1&1\\
0&1&0&1
\end{bmatrix};\quad X=\begin{bmatrix}
1\\0\\0
\end{bmatrix}.
$$
We try to find the minimum $A$-annihilator of $X$. Set $X_0=X$. It is easy to see that
$$
X_1=A\lvtimes X_0\in {\cal V}_6\in {\cal I}_A.
$$
Hence, we can find the annihilator of $X$ in the space of $\R^6$.
Calculating
$$
\begin{array}{lll}
X_1=\begin{bmatrix}1\\1\\0\\0\\0\\0\end{bmatrix};&
X_2=\begin{bmatrix}1\\1\\1\\1\\0\\0\end{bmatrix};&
X_3=\begin{bmatrix}2\\2\\1\\0\\0\\1\end{bmatrix};\\
X_4=\begin{bmatrix}2\\1\\1\\2\\1\\1\end{bmatrix};&
X_5=\begin{bmatrix}3\\3\\1\\-1\\-1\\0\end{bmatrix};&
X_6=\begin{bmatrix}3\\2\\2\\4\\2\\2\end{bmatrix},
\end{array}
$$
it is easy to verify that $X_1,~X_2,~X_3,~X_4,~X_5$ are linearly independent. Moreover,
$$
X_6=X_1+X_2+X_3+X_4-X_5.
$$
The minimum $A$-annihilator of $X=X_0$ follows as
$$
p(x)=x^6+x^5-x^4-x^3-x^2-x.
$$
\end{enumerate}
\end{exa}

\subsection{Invariant Subspace on V-equivalence Space}

Recall that
$$
\Omega_{{\cal M}}:={\cal M}/\lra;\quad \Omega_{{\cal M}}^{n}:={\cal M}_{\cdot\times n}/\lra.
$$
We extend the vector product to the equivalence spaces.

\begin{dfn}\label{d6.8.1} Let $A\in {\cal M}$ and  $B\in {\cal M}_{\cdot\times q}$. Then we define $\lvtimes: \Sigma_{{\cal M}}\times \Omega_{{\cal M}}^q\ra \Omega_{{\cal M}}^q$ as
\begin{align}\label{6.8.1}
\A\lvtimes \lceil B\rceil:=\lceil A\lvtimes B \rceil.
\end{align}
\end{dfn}

The following proposition shows that (\ref{6.8.1}) is well defined.

\begin{prp}\label{p6.8.2} (\ref{6.8.1}) is independent of the choice of $A$ and $B$.
\end{prp}

\noindent\textit{Proof}. Assume $A_1\in \A$ is irreducible and $A_1\in {\cal M}_{m\times n}$; $B_1\in \lceil B\rceil$ is also irreducible and $B_1\in {\cal M}_{p\times q}$. $A_i=A_1\otimes I_{i}$ and $B_j=B_1\otimes \I_j$.
Set $s=n\vee p$, $t=ni\vee pj$, and $s\xi=t$. Then
$$
\begin{array}{l}
\left(A_i\lvtimes B_j\right)\\
=\left(A_i\otimes I_{t/ni}\right)\left(B_j\otimes \I_{t/pj}\right)\\
=\left(A_1\otimes I_{i}\otimes I_{t/ni}\right)\left(B_1\otimes \I_{t}\otimes \I_{t/pj}\right)\\
=\left(A_1\otimes \otimes I_{t/n}\right)\left(B_1\otimes \I_{t/p}\right)\\
=(A_1\lvtimes B_1)\otimes \I_{\xi}\lra (A_1\lvtimes B_1).
\end{array}
$$
\hfill $\Box$

Precisely speaking, because ${\cal V}$ is not a vector space, ``invariant subspace" is not a rigorous subspace. But $\Omega_{{\cal V}}:={\cal V}/\lra $ is a vector space. It is easy to see that the results about ${\cal V}$, and ${\cal M}_{\cdot\times n}$ can be extended to $\Omega_{{\cal V}}$ and $\Omega_{{\cal M}}$. For instance, $\A\in \Sigma_{\mu}$ is a bounded operator on $\Omega_{{\cal V}}$ if and only if, $\mu_y=1$.

\subsection{Generalized Linear System}

\begin{dfn}[\cite{kopnet,liu08}]\label{d7.8.1} Let $S$ be a semigroup, $X$ a Hausdorff space. A mapping $\varphi:S\times X\ra X$ is called a topological dynamics, if
\begin{enumerate}
\item[(1)]
\begin{align}\label{7.8.1}
\varphi(s_1,\varphi(s_2,x))=\varphi(s_1s_2,x),\quad s_1,s_2\in S,\; x\in X.
\end{align}
\item[(2)]
\begin{align}\label{7.8.2}
\varphi(e,x)=x,\quad x\in X,
\end{align}
where $e\in S$ is the identity of $S$.
\item[(3)] For each $s\in S$, $\varphi_s:X\ra X$ is continuous.
\end{enumerate}
\end{dfn}

\begin{thm}\label{t7.8.2} Let $S={\cal M}$ and $X={\cal V}$. Then $\lvtimes: S\times X\ra X$ is a topological dynamics. Precisely, let $\F=\R$, and $A\in {\cal M}$, then
\begin{align}\label{7.8.3}
\begin{array}{l}
x(t+1):=A\lvtimes x(t),\\
x(0)=x_0\in {\cal V},
\end{array}
\end{align}
is a topological dynamics, which is called a generalized linear system.
\end{thm}

To prove this theorem we have to prove that the conditions (1)--(3) are satisfied. It is not difficult to see that (2) and (3) are satisfied. So we need to prove (\ref{7.8.1}). We state it as the following lemma.

\begin{lem}\label{l7.8.3} For any two matrices $A,~B\in {\cal M}$ and any vector $X\in {\cal V}$, it holds that
\begin{align}\label{7.8.4}
(A\ltimes B)\lvtimes X=A\lvtimes (B\lvtimes X).
\end{align}
\end{lem}

\noindent\textit{Proof.}
Assume $A\in {\cal M}_{m\times n}$, $B\in {\cal M}_{p\times q}$, $X\in \R^r$. Then
\begin{align}\label{a1.1}
\begin{array}{l}
\left(A\ltimes B\right)\lvtimes X\\
=\left[(A\otimes I_{t/n})(B\otimes I_{t/p})\right]\lvtimes X\\
=\left\{\left[(A\otimes I_{t/n})(B\otimes I_{t/p})\right]\otimes I_{\frac{sp}{qt}}\right\}\left(X\otimes {\bf 1}_{s/r}\right)\\
=\left(A\otimes I_{\frac{sp}{nq}}\right)\left(B\otimes I_{s/q}\right)\left(X\otimes {\bf 1}_{s/r}\right)\\
=\left(A\otimes I_{\frac{sp}{nq}}\right)\left\{\left[\left(B\otimes I_{\ell/q}\right)\left(X\otimes {\bf 1}_{\ell/r}\right)\right]\otimes {\bf 1}_{\phi}\right\}\\
=\left(A\otimes I_{\frac{sp}{nq}}\right)\left[\left(B\lvtimes X\right)\otimes {\bf 1}_{\phi}\right],
\end{array}
\end{align}
where
$$
\begin{array}{l}
t=n\vee p\\
s=\left(\frac{qt}{p}\right)\vee r\\
\ell=q\vee r\\
s=\ell\phi.
\end{array}
$$
Note that $B\lvtimes X\in \R^{\frac{p\ell}{q}}$. By definition if
\begin{align}\label{a1.2}
\frac{n\vee \frac{p\ell}{q}}{n}=\frac{sp}{nq}
\end{align}
and
\begin{align}\label{a1.3}
\frac{n\vee \frac{p\ell}{q}}{\frac{p\ell}{q}}=\phi,
\end{align}
then (\ref{a1.1}) becomes
$$
A\lvtimes \left(B\lvtimes X\right),
$$
and we are done.

It is clear that both (\ref{a1.2}) and (\ref{a1.3}) are equivalent to
\begin{align}\label{a1.4}
n\vee \frac{p\ell}{q}=\frac{sp}{q}.
\end{align}
Hence as long as (\ref{a1.4}) holds, we are done.

In the following we prove (\ref{a1.4}).

Since $t=n\vee p$, assume
$$
n\wedge p=u,
$$
then
$$
n=\a u,\quad p=\b u, \quad t=\a \b u, \quad \a\wedge \b=1.
$$
Next,
$$
s=\frac{qt}{p}\vee r=q\a\vee r.
$$
Assume $d_1,d_2,\cdots,d_s$ are the set of prime factors of $q\wedge r$. Then we can express $q$ and $r$ as follows:
$$
\begin{array}{l}
q=d_1^{Q_1}d_2^{Q_2}\cdots d_s^{Q_s}q_0,\\
r=d_1^{R_1}d_2^{R_2}\cdots d_s^{R_s}r_0,\\
\end{array}
$$
where $q_0\wedge r_0=1$, $q_0\wedge d_i=1$, and $r_0\wedge d_i=1$, $Q_i>0$, $R_i>0$, $i=1,\cdots,s$.

Assume $\epsilon_1,\cdots, \epsilon_t$ are the common prime factors of $\a$ and $r$, which are distinct from $\{d_1,\cdots,d_s\}$.
Then we can express $\a$ and $r$ respectively as follows:
$$
\begin{array}{l}
\a=d_1^{A_1}d_2^{A_2}\cdots d_s^{A_s}\epsilon_1^{B_1}\epsilon_2^{B_2}\cdots \epsilon_t^{B_t}\a_0,\\
r=d_1^{R_1} d_2^{R_2}\cdots d_s^{R_s}\epsilon_1^{C_1}\epsilon_2^{C_2}\cdots \epsilon_t^{C_t} r_{00},\\
\end{array}
$$
where $B_i> 0$, $C_i>0$, $i=1,\cdots,t$, $A_j\geq 0$, $j=1,\cdots,s$. $A_j$ might be zero, if $\a$ has no such a factor. Now $\a_0\wedge r_{00}=1$, $\a_0\wedge d_j=1$, $r_{00}\wedge d_j=1$, $j=1,\cdots,s$, $\a_0\wedge \epsilon_i=1$, $r_{00}\wedge \epsilon_i=1$, $i=1,\cdots,t$.

Now we can calculate the left hand side of (\ref{a1.4}):
$$
\ell=q\vee r= d_1^{\max(Q_1, R_1)}\cdots d_s^{\max(Q_s, R_s)}\epsilon_1^{C_1}\cdots \epsilon_t^{C_t} q_0r_{00};
$$
$$
\frac{\ell}{q}=d_1^{\max(Q_1, R_1)-Q_1}\cdots d_s^{\max(Q_s, R_s)-Q_s}\epsilon_1^{C_1}\cdots \epsilon_t^{C_t}r_{00};
$$
$$
\frac{p\ell}{q}=d_1^{\max(Q_1, R_1)-Q_1}\cdots d_s^{\max(Q_s, R_s)-Q_s}\epsilon_1^{C_1}\cdots \epsilon_t^{C_t}r_{00}\b u.
$$
Hence, we have
\begin{align}\label{a1.5}
\begin{array}{l}
n\vee\frac{p\ell}{q}=\left(u\a_0\d_1^{A_1}\cdots\d_s^{A_s}\epsilon_1^{B_1}\cdots \epsilon_t^{B_t} \right)\\
~~\vee \left(d_1^{\max(Q_1, R_1)-Q_1}\cdots d_s^{\max(Q_s, R_s)-Q_s}\epsilon_1^{C_1}\cdots \epsilon_t^{C_t}r_{00}\b u\right)\\
~=u\a_0r_{00}\b d_1^{\max(\max(Q_1, R_1)-Q_1, A_1)}\cdots \\
~~d_s^{\max(\max(Q_s, R_s)-Q_s, A_s)}\epsilon_1^{\max(C_1, B_1)} \cdots \epsilon_t^{\max(C_t, B_t)}.
\end{array}
\end{align}

Next, we calculate the right hand side of (\ref{a1.4}):
$$
q\a=d_1^{A_1+Q_1}\cdots d_s^{A_s+Q_s}\epsilon_1^{B_1}\cdots \epsilon_t^{B_t}q_0\a_0;
$$
$$
\begin{array}{l}
s=(q\a)\vee r=d_1^{\max(A_1+Q_1, R_1)}\cdots d_s^{\max(A_s+Q_s, R_s)}\\
~~\epsilon_1^{\max(B_1, C_1)}\cdots \epsilon_t^{\max(B_t, C_t)}q_0\a_0 r_{00}.
\end{array}
$$
Then we have
\begin{align}\label{a1.6}
\begin{array}{l}
\frac{sp}{q}=d_1^{\max(A_1+Q_1, R_1)-Q_1}\cdots d_s^{\max(A_s+Q_s, R_s)-Q_s}\\
~~\epsilon_1^{\max(B_1, C_1)}\cdots \epsilon_t^{\max(B_t, C_t)}\b u\a_0 r_{00}.
\end{array}
\end{align}
Comparing (\ref{a1.5}) with (\ref{a1.6}), one sees easily that to prove (\ref{a1.4}) it is enough to prove
\begin{align}\label{a1.7}
\begin{array}{ccl}
\max(A_i+Q_i, R_i)-Q_i&=&\max(\max(Q_i, R_i)-Q_i, A_i),\\
~&~&\quad i=1,\cdots,s.
\end{array}
\end{align}
Note that the left hand side of (\ref{a1.7}) equals to
$$
\max(A_i, R_i-Q_i).
$$
Then it is easy to verify that
when $R_i\geq Q_i$ both sides of (\ref{a1.7}) equal to $\max(A_i,R_i-Q_i)$, and when $R_i<Q_i$, both sides equal to $A_i$. The proof is completed.
\hfill $\Box$

\begin{rem}\label{r7.8.4}
\begin{enumerate}
\item The invariance subspace has been discussed in previous subsections. Then it is clear that the general linear system (\ref{7.8.3}) is dimension bounded, which means there exists an $n$ such that $\dim(x(t)<n$, $\forall t$, if and only if, $A$ is bounded.
\item It is easy to extend $(\ref{7.8.3})$ to equivalence space $\Sigma_{{\cal M}}$ and the vector space $\Omega_{{\cal V}}$.
\end{enumerate}
\end{rem}


\section{Concluding Remarks}

Matrix theory is one of the most fundamental and useful tools in modern science and technology. But one of the major weaknesses is its dimension restriction. To overcome this barrier,  the purpose of this paper is to set up a framework for an almost dimension-free matrix theory.

First, we review the STP ($\ltimes$), which extends the conventional matrix product to overall matrices ${\cal M}$. The related monoid structure for $\left({\cal M},~\ltimes\right)$ is obtained. The M-equivalence $\sim$ is proposed. A lattice structure over each equivalence class is obtained. The equivalence space, as the quotient space ${\cal M}/\sim$ is introduced and discussed.

Second, the set of overall matrices is partitioned into subspaces as ${\cal M}=\bigcup_{\mu\in \Q_+}{\cal M}_{\mu}$. The STA ($\lplus$) is proposed. Under this addition the quotient spaces $\Sigma_{\mu}={\cal M}_{\mu}/\sim$
become  vector spaces. Certain geometric and algebraic structures are revealed. Including topological structure, inner product structure, differential manifold structure, etc.

Particularly, when $\mu=1$, (corresponding to square matrices) we have extended the Lie algebra and Lie group theory to ${\cal M}_1$. A fiber bundle structure, called the discrete bundle,  is proposed for ${\cal M}_{\mu}$ and the extended Lie group and Lie algebra.

Finally, the set of overall vectors ${\cal V}$ are considered as a universal vector space, based on the vector equivalence $\lra$.  A matrix $A$ of any dimension can be considered either a linear mapping on ${\cal V}$, or a subspace generated by its columns.
The $A$-invariant subspace is discussed in details. Many key concepts such as eigenvalue/eigenvector, characteristic polynomial of a matrix have been extended from square matrices to no-square matrices.

It was said by Asimov that ``Only in mathematics is there no significant correction - only existence. Once the Greeks had developed the deductive method, they were correct in what they did, correct for all the time."  \cite{tan13} All extensions we did in this paper consist with the classical ones. That is, when the dimension restrictions required by the classical matrix theory are satisfied the new operators proposed in this paper coincide with the classical ones.

There are many questions remain for further discussion. For instance, is it possible to construct an equivalence over ${\cal M}$, which is consistent with certain matrix product, such that the quotient space becomes a vector space?

The followings are some possible equivalences on ${\cal M}$.
\begin{enumerate}
\item  Equivalence 1:

\begin{dfn}\label{d7.1.1} Let $A,~B\in {\cal M}$. $A$ is said to be equivalent to $B$, denoted by $A \simeq B$, if there exist $I_i$, $I_j$, $\I_s$, $\I_t$ such that
\begin{align}\label{7.1.1}
\I_\a^T\otimes A\otimes I_i=\I_\b^T\otimes B\otimes I_j.
\end{align}
\end{dfn}

It is easy to verify that $\simeq$ is an equivalence relation. Moreover, similar to M-equivalence or vector equivalence, we have the following result:

\begin{thm}\label{t7.1.2} Assume $A\simeq B$, then there exists a $\Lambda$, such that
\begin{align}\label{7.1.2}
\begin{array}{l}
A=\I^T_p\otimes \Lambda \otimes I_s,\\
B=\I^T_q\otimes \Lambda \otimes I_t.
\end{array}
\end{align}
\end{thm}
Hence the lattice structure similar to M-equivalence exists.

It may be considered as a combination of M- and V-equivalences. Unfortunately, (i) it is not consistent with STP ($\ltimes$); (ii) the quotient space is not a vector space.

\item Equivalence 2:

\begin{dfn}\label{d7.1.3} Let $A,~B\in {\cal M}$. $A$ is said to be  equivalent to $B$, if there exist $\I_i$, $\I_j$, $\I_s$, $\I_t$ such that
\begin{align}\label{7.1.3}
\I_\a^T\otimes A\otimes \I_i=\I_\b^T\otimes B\otimes \I_j.
\end{align}
\end{dfn}

The lattice structure can also be determined in a similar way. Moreover, the quotient space is a vector space. Unfortunately, a proper product, which is consistent with the equivalence, is unknown.

\end{enumerate}

Further geometric/algebraic structures may be investigated.

\begin{enumerate}
\item More geometric structure on equivalence space could be interesting. For instance, a Riemannian geometric structure or a Symplectic geometric structure may be posed on the equivalence space.

\item Under the STP and $\lplus$ ($\lminus$), for any $A\in {\cal M}_{\mu}$ and $B\in {\cal M}_{\lambda}$,
\begin{align}\label{7.1.4}
[A,B]=A\ltimes B\lminus B\ltimes A
\end{align}
is well defined. Moreover, the three requirements (\ref{4.3.1})--(\ref{4.3.3}) in Definition \ref{d4.3.1} can also be satisfied (under obvious modification). Exploring the properties of this generalized Lie algebra is challenging and interesting.

\end{enumerate}

One may be more interested in its applications. For instance, can we use the extended structure proposed in this paper to the analysis and control of certain dynamic systems? Particularly, we may consider
the following special cases:
\begin{enumerate}

\item Consider a dynamic system
$$
\dot{x}=A(t)x,\quad x\in \R^n,
$$
where $A(t)$ satisfies
$$
\dot{A}_t=V(x),
$$
where $V(x)$ is a vector field on ${\cal M}_{n\times n}$. What can we say about this system?
Is it possible to extend this system to the equivalence space $\Sigma$?

\item A dimension-varying dynamic control system as
$$
\begin{cases}
x(t+1)=A\lvtimes x(t)\lvplus B\lvtimes u(t)\\
y(t)=C\lvtimes x(t),
\end{cases}
$$
where $x(t)\in {\cal V}$.

What can we say about this, say, controllability? observability etc.?

\end{enumerate}

In one word, this paper could be the beginning of investigating dimension-free matrix theory and its applications.

\vskip 2mm

\centerline{\bf Acknowledgment}

The author would like to thank the anonymous reviewers for their valuable suggestions, comments, and detailed typo corrections.

\end{document}